\documentclass[11pt,a4]{article}
\usepackage[T2A]{fontenc}
\usepackage[cp866]{inputenc}
\usepackage[mathscr]{eucal}
\usepackage{amssymb}
\usepackage{amsmath}
\usepackage{fancyhea}
\renewcommand{\leq}{\leqslant}

\renewcommand{\C}{{\mathbb C}}
\newcommand{\N}{{\mathbb N}}
\newcommand{\R}{{\mathbb R}}

\renewcommand{\k}{\rule{0.7em}{0.7em}}
\begin{document}
\sloppy
\headrulewidth = 2pt
\pagestyle{fancy}
\lhead[\scriptsize V.N.Gorbuzov]
{\scriptsize V.N.Gorbuzov}
\rhead[\it \scriptsize Integral equivalence of multidimensional differential systems]
{\it \scriptsize Integral equivalence of multidimensional differential systems}
\headrulewidth=0.25pt

{\normalsize 

\thispagestyle{empty}

\mbox{}
\\[-3.5ex]
\centerline{
{
\bf 
INTEGRAL EQUIVALENCE OF MULTIDIMENSIONAL 
}
}
\\[0.75ex]
\mbox{}\hfill
{
\bf 
DIFFERENTIAL SYSTEMS\,}\footnote[1]{
The basic results of this paper have been published in the monographs
"Integrals of systems of differential equations", Grodno, 2006 [1] and 
"Integrals of systems of total differential equations", Grodno, 2005 [2],
and the journal "Vestnik of the Yanka Kupala Grodno State Univ.", 
2005, Ser. 2, No. 2, 10-29 [3]. 
}
\hfill\mbox{}
\\[2.75ex]
\centerline{
\bf 
V.N. Gorbuzov
}
\\[1.75ex]
\centerline{
\it 
Department of Mathematics and Computer Science, 
Yanka Kupala Grodno State University,
}
\\[1ex]
\centerline{
\it 
Ozeshko 22, Grodno, 230023, Belarus
}
\\[1ex]
\centerline{
E-mail: gorbuzov@grsu.by
}
\\[5ex]
\centerline{{\large\bf Abstract}}
\\[1ex]
\indent
The bases of the theory of integrals 
for multidimensional differential systems are stated. 
The integral equivalence of total differential systems, 
linear homogeneous systems of partial differential equations, and 
Pfaff systems of equations is established.
\\[2ex]
\\[0.75ex]
\\[3.5ex]
{\large\bf Contents} 
\\[1ex]
{\bf  Introduction}
                                                                                                                    \hfill\ 2
\\[0.5ex]
{\bf  
1. First integrals of  total differential system}
                                                                                                                    \hfill\ 5
\\[0.5ex]
\mbox{}\hspace{1em}
1.1. First integral                                                                                             \dotfill\ 5
\\[0.5ex]
\mbox{}\hspace{1em}
1.2. Basis of first integrals                                                     \dotfill\ 7
\\[0.5ex]
\mbox{}\hspace{1em}
1.3. Dimension of  ba\-sis of first integrals for completely solvable systems 
                                                                                                                         \dotfill\ 8
\\[1ex]
\noindent
{\bf  
2. First integrals for linear homogeneous system of partial differential 
\\
\mbox{}\hspace{1em}\!
equations}      \hfill\ 10
\\[0.5ex]
\mbox{}\hspace{1em}
2.1. Basis of first integrals 
                                                                                                                                   \dotfill\ 10
\\[0.5ex]
\mbox{}\hspace{1em}
2.2. Incomplete system
                                                                                                                                   \dotfill\ 11
\\[0.5ex]
\mbox{}\hspace{1em}
2.3. Complete system
                                                                                                                                   \dotfill\ 13
\\[0.5ex]
\mbox{}\hspace{1em}
2.4. Dimension of integral basis
                                                                                                                                   \dotfill\ 17
\\[1ex]
\noindent
{\bf  
3. Dimension of  integral basis for not completely solvable 
\\
\mbox{}\hspace{1em}\!
total differential system}                                                                                         \hfill\ 21
\\[1ex]
\noindent
{\bf  
4. First integrals for Pfaff system of equations}                                                              \hfill\ 26
\\[0.5ex]
\mbox{}\hspace{1em}
4.1. Integrally equivalent Pfaff systems of equations 
                                                                                                                                   \dotfill\ 26
\\[0.5ex]
\mbox{}\hspace{1em}
4.2. Integral basis                                                                                                   \dotfill\ 26
\\[0.5ex]
\mbox{}\hspace{1em}
4.3. Existence criterion of first integral
                                                                                                                                   \dotfill\ 28
\\[0.5ex]
\mbox{}\hspace{1em}
4.4. Integral equivalence with linear homogeneous system 
\\
\mbox{}\hspace{2.9em}
of partial differential equations                                                                               \dotfill\ 29
\\[0.5ex]
\mbox{}\hspace{1em}
4.5. Transformation of Pfaff system of equations by known first integrals
                                                                                                                                   \dotfill\ 31
\\[0.5ex]
\mbox{}\hspace{1em}
4.6. Closed systems                                                                                                  \dotfill\ 32
\\[0.5ex]
\mbox{}\hspace{1em}
4.7. Interpretation of closure in terms of differential forms                               \dotfill\ 33
\\[0.5ex]
\mbox{}\hspace{1em}
4.8. Nonclosed  systems                                                                                          \dotfill\ 37
\\[0.5ex]
\mbox{}\hspace{1em}
4.9. Integral equivalence with  total differential system                                        \dotfill\ 39
\\[0.5ex]
{\bf  References}                                                                                                  \dotfill\ 40

\newpage

\mbox{}
\\[-1.75ex]
\centerline{
\large\bf  
Introduction
}
\\[1.25ex]
\indent
Basis of the general theory of differential equations are 
the theory of solutions and the theory of integrals. 
The functional-analytical research of integrals is most deeply developed 
for ordinary differential systems and 
linear systems of partial differential equations.
\vspace{0.15ex}

The initial problem of the integration in quadratures for differential equations
has led to necessity to develop methods of analytical and qualitative researches both 
for solutions and integrals of differential systems.
So, J. Liouville [4, 5] considered the problem of the integration in quadratures for 
the Riccati equation. His investigation gave the classical problems about a form of solutions and 
about development of methods for finding solutions of given forms.
\vspace{0.15ex}

For the first time the problem of building a general integral from first integrals was considered by J. Jacobi in [6 -- 8].
He also introduced the notion of a last multiplier (also known in publications as Jacobi's last multiplier) 
and he used this notion of a last multiplier to solve the problem of finding a general integral.
\vspace{0.15ex}

The profound researches, which are the base of the theory of integrals, are due to  
F.G. Minding [9],  A.V. Letnikov [10], and A.N. Korkine [11]. 
\vspace{0.15ex}

One of  such problems is the Darboux problem [12] about building of a general integral by known partial integrals 
and about the  form of a general integral for the ordinary differential equation of the first order. 
The Darboux problem for ordinary differential systems, total differential systems, and 
systems of partial differential equations was considered in the monographs [1, 2, 13, 14] 
and in the articles [15 -- 46].
\vspace{0.15ex}

The method of last multiplier was developed by S. Lie. 
He has given the new interpretation of this method and has created
the theory of infinitesimal transformations in [47, 48].
In his papers were selected the fundamental approaches of the integration.
These approaches are the base of the theory of closed differential systems.
Methods of the integration were considered in  [49] 
with regard to the uniform positions of the group analysis. 
This approach is gave the theoretical base for the classification [50, 51]
of these methods.
\vspace{0.15ex}

At the beginning of the twentieth century the interest to global researches for 
differential systems has decreased. 
But in the middle of the twentieth century the interest in it has ap\-pe\-a\-r\-ed again.
Let us note only the monographies:
N.M. Gjunter [52], E. Cartan [53, 54], N.G. Che\-bo\-ta\-rev [55],
L. Eis\-en\-hart [56], P.K. Rashevskii [57], H. Cartan [58], 
A.S. Galiullin [59], L.V. Ov\-sian\-ni\-kov [60], 
N.P. Erugin [61], 
E.A. Barbashin [62], A.M. Samoilenko [63], 
\linebreak
P. Olver [64], A. Goriely [65]. 
It hap\-pen\-ed for the reason that these results find applications in the 
ma\-the\-ma\-ti\-cal physics [66 -- 71].
\vspace{0.15ex}

The intensive development of the qualitative theory of differential equations 
(the base of this theory was founded by H. Poincar\'{e} [72])
also has led to solving some problems for the theory of integrals. 
So, the analytic structure of integrals and of integrating multipliers 
in a neighbourhood of a center was obtained by A.M. Liapunov [73] and 
N.A. Sakharnikov [74]  respectively.
Investigation of limit cycles for differential systems on the base of partial integrals and 
integrating multipliers was considered by M.V. Dolov [75 -- 77].
Behaviour of trajectories for ordinary autonomous differential systems of the second order 
having the integral curves of the concrete forms and with the special qualitative properties 
\vspace{0.5ex}
was researched\! in\! [16; 78 -- 90].

The subject of our investigation is a system of total differential equations 
\\[2ex]
\mbox{}\hfill                                                
$
dx = X(t,x)\,dt, 
$
\hfill(TD)
\\[2ex]
where\vspace{0.5ex}
$\!t\in\R^m,\, x\in\R^n,\, m\leq n,\, 
dt\!={\rm colon}\,(dt^{}_1,\ldots,dt^{}_m\!),\, 
dx\!={\rm colon}\,(dx^{}_1,\ldots,dx^{}_n\!),\,X\!(t,x)=\!\!$ 
$=\|X^{}_{ij}(t,x)\|$ is an $n\times m$\vspace{0.5ex}
matrix with entries $X^{}_{ij}\colon \Pi\to\R,\ 
i=1,\ldots,n,\ j=1,\ldots,m,\ \Pi$ is a domain of the extended space $\R^{m+n};$
a linear homogeneous system of partial differential equations
\\[0.75ex]
\mbox{}\hfill                                                
$
{\frak L}^{}_j (x)\,y=0,
\quad 
j=1,\ldots,m,
$
\hfill{\rm(\!$\partial$\!)}
\\[1.75ex]
where $y\in\R,\ x\in\R^n,\ m\leq n,$
the linear differential ope\-ra\-tors of first order
\\[2ex]       
\mbox{}\hfill                                                
$
\displaystyle
{\frak L}^{}_j(x)=
\sum\limits_{i=1}^n u^{}_{ji}(x)\,
\partial^{}_{x^{}_i}$ 
\ for all $x\in G,
\quad 
j=1,\ldots,m,
$
\hfill(0.1)
\\[2ex]
with coordinates $u^{}_{ji}\colon G\to\R,\ j=1,\ldots,m,\ i=1,\ldots,n,$
a domain $G\subset\R^n;$ and a Pfaff system of equations 
\\[1.5ex]
\mbox{}\hfill                                                 
$
\omega^{}_j(x)=0,
\quad 
j=1,\ldots,m,
$
\hfill(Pf)
\\[2ex]
where $x\in\R^n,\ m\leq n,$ the linear differential forms
\\[2ex]
\mbox{}\hfill                                                  
$
\displaystyle
\omega^{}_j(x)=
\sum\limits_{i=1}^n
w^{}_{ji}(x)\,dx^{}_i$
\ for all $x\in G,
\quad 
j=1,\ldots,m,
$
\hfill(0.2)
\\[2.25ex]
with coefficients $w^{}_{ji}\colon G\to\R,\ j=1,\ldots,m,\ i=1,\ldots,n,$
a domain $G\subset\R^n.$ 

We recall that by domain we mean open arcwise connected set.

Let the linear differential ope\-ra\-tors (0.1) of the system of partial differential equations (\!$\partial$\!) and
the linear differential forms (0.2) of the Pfaff system of equations (Pf) be not linearly bound
on the domain $G$ [91, pp. 105 -- 115].
Note that the operators ${\frak L}^{}_j\,,\ j=1,\ldots,m$
(the 1-forms $\omega^{}_j\,,\ j=1,\ldots,m)$ are called 
{\it linearly bound} on the domain $G$ if 
these operators (1-forms) are 
linearly dependent in any point of the domain $G$
[60, pp. 113 -- 114].

The system (TD) is induced the $m$ linear differential ope\-ra\-tors of first order
\\[2ex]
\mbox{}\hfill                                                   
$                                       
\displaystyle
{\frak X}^{}_j(t,x)=
\partial^{}_{t^{}_j}+ 
\sum\limits_{i=1}^{n}
X^{}_{ij}(t,x)\,
\partial^{}_{x^{}_i}$
\ for all $(t,x)\in\Pi,
\quad 
j=1,\ldots,m.
$
\hfill(0.3)
\\[2ex]
We'll say that these operators are 
{\it operators of differentiation by virtue of sys\-tem} (TD).

Under the condition $m=1,$ we have the system (TD) is an
ordinary differential system of $n$ order.

The system (TD) is said to be {\it completely solvable on a domain} 
$\Pi^{\;\!\prime}\subset \Pi$ if for any point $(t^{}_0\,,x^{}_0)\in \Pi^{\;\!\prime}$ 
there exists a unique solution of the Cauchy problem with initial data $(t^{}_0\,,x^{}_0)$
[1, p. 17]. If $\Pi^{\;\!\prime}= \Pi,$ then we say that the system (TD) is 
{\it completely solvable}.

Suppose $X\in C^1(\Pi),$ i.e.,  the functions  
$X^{}_{ij}\,,\ i=1,\ldots,n,\ j=1,\ldots,m,$ are continuously differentiable on the 
domain $\Pi.$ Then the system (TD) is completely solvable if the {\it Frobenius theorem} 
 [1, pp.~17~-- 25; 58, pp. 290 -- 297; 92, pp. 309 -- 311; 93, p. 21] is true. 
\vspace{0.5ex}

{\bf Theorem 0.1}({\sl the Frobenius theorem}).                                    
{\it
The system {\rm (TD)} with $X\in C^1(\Pi)$ is completely solvable if and only if
the Frobenius conditions hold}
\\[1.5ex]
\mbox{}\hfill                                                   
$
\displaystyle
\partial^{}_{t^{}_j}
X^{}_{i\zeta}(t,x)+
\sum\limits_{\xi=1}^n
X^{}_{\xi j}(t,x)
\partial^{}_{x^{}_{\xi}}
X^{}_{i\zeta}(t,x)=
\partial^{}_{t^{}_{\zeta}}
X^{}_{ij}(t,x)+
\sum\limits_{\xi=1}^n
X^{}_{\xi\zeta}(t,x)
\partial^{}_{x^{}_{\xi}}
X^{}_{ij}(t,x)
\hfill
$
\\
\mbox{}\hfill {\rm (0.4)}
\\
\mbox{}\hfill
for all $(t,x)\in\Pi,\ \ 
i=1,\ldots,n,\ \
j=1,\ldots,m,\ \
\zeta=1,\ldots,m.
\hfill
$
\\[2ex]
\indent
Using the operators (0.3), we obtain the Frobenius conditions (0.4)
are represented via Poisson brackets as the system of the operator identities
\\[2ex]
\mbox{}\hfill                                                 
$
\bigl[{\frak X}^{}_j(t,x),
{\frak X}^{}_\zeta(t,x)\bigr]={\frak O}$
\ for all $(t,x)\in\Pi,
\ \ \
j=1,\ldots,m,\ \,
\zeta=1,\ldots,m,
$
\hfill (0.5)
\\[2ex]
where ${\frak O}$ is the null operator.

The system (TD) is the Pfaff system of equations 
\\[2ex]
\mbox{}\hfill                                                 
$
\eta^{}_i (t,x)=0,
\quad 
i=1,\ldots,n,
$
\hfill(0.6)
\\[2ex]
with the linear differential forms
\\[2ex]
\mbox{}\hfill                                                 
$                                       
\displaystyle
\eta^{}_i (t,x)=
dx^{}_i- 
\sum\limits_{j=1}^{m}
X^{}_{ij}(t,x)\,dt^{}_j$
\ for all 
$(t,x)\in\Pi,\ \ i=1,\ldots,n.
$
\hfill(0.7)
\\[2ex]
\indent
The Frobenius conditions (0.4) for completely solvability 
of system (TD) (the Pfaff system of equations (0.6)) are represented via 
differential 1-forms (0.7) as the system of  exterior differential identities [58, pp. 290 -- 297]
\\[2ex]
\mbox{}\hfill                                                  
$
d\eta^{}_i(t,x)\wedge
\Bigl(\,\mathop{\wedge}\limits_{\xi=1}^{n}
\eta^{}_{\xi}(t,x)\Bigr)= 0$
\ for all $(t,x)\in\Pi,
\ \  i=1,\ldots,n.
$
\hfill(0.8)
\\[2.25ex]
\indent
Consider an autonomous total differential system
\\[2ex]
\mbox{}\hfill                                                
$
dx=X(x)\,dt, 
$
\hfill(ATD)
\\[2.25ex]
where\vspace{0.5ex}
$t\in\R^m,\ x\in\R^n,\ m\leq n,\ 
dt={\rm colon}\,(dt^{}_1,\ldots,dt^{}_m),\ 
dx={\rm colon}\,(dx^{}_1,\ldots,dx^{}_n),$ 
the entries of the $n\times m$ matrix\vspace{0.5ex}
$X(x)\!=\!\|X^{}_{ij}(x)\|\!$ are  $X^{}_{ij}\colon G\to\R,\, i=1,\ldots,n,\, j=1,\ldots,m,\
G\!$ is a domain of the phase space $\!\R^n.\!\!$
We say that the linear differential ope\-ra\-tors of first order
\\[2ex]       
\mbox{}\hfill                                                
$
\displaystyle
{\frak X}^{}_j(t,x)=
\partial^{}_{t^{}_j}+
\sum\limits_{i=1}^n 
X^{}_{ij}(x)\,
\partial^{}_{x^{}_i}$ 
\ for all $(t,x)\in \R^m\times G,
\quad  j=1,\ldots,m,
$
\hfill (0.9)
\\[2ex]
induced by this system are 
{\it nonautonomous operators} of differentiation by virtue of sys\-tem (ATD).
The linear differential ope\-ra\-tors of first order
\\[2ex]       
\mbox{}\hfill                                                
$
\displaystyle
{\frak A}^{}_j(x)=
\sum\limits_{i=1}^n 
X^{}_{ij}(x)\,
\partial^{}_{x^{}_i}$ 
\ for all $x\in G,
\quad j=1,\ldots,m, 
$
\hfill(0.10)
\\[2ex]
are called {\it autonomous operators} of differentiation by virtue of sys\-tem (ATD).

Let $\!X\!\in\! C^1(G).\!$ Then the Frobenius conditions (0.5) for completely solvability 
of sys\-tem (ATD) are the identities [1, pp. 112 -- 113]
\\[2ex]
\mbox{}\hfill                                                 
$
\bigl[{\frak A}^{}_j(x),
{\frak A}^{}_\zeta(x)\bigr]={\frak O}$
\ for all $x\in G, 
\quad 
j=1,\ldots,m,\ \
\zeta=1,\ldots,m.
$
\hfill (0.11)
\\[2ex]
\indent
The system (\!$\partial$\!) with
\vspace{0.5ex} 
${\frak L}^{}_j\in C^1(G),\ j=1,\ldots,m$
(i.e., the coordinates $u^{}_{ji}\,,\ j=1,\ldots,m,$ $i=1,\ldots,n,$
are continuously differentiable on the domain $G)$ 
is said to be {\it complete} [92, p. 521; 53, p. 117] if
Poisson brackets of any two differential operators (0.1) are
the linear combination of operators (0.1)
\\[2ex]
\mbox{}\hfill                                           
$
\displaystyle
\bigl[{\frak L}^{}_j(x),
{\frak L}^{}_{\zeta}(x)\bigr]=
\sum\limits_{\nu=1}^m 
A^{}_{j\zeta\nu}(x)
{\frak L}^{}_{\nu}(x)$
\ for all $x\in G,
\quad
j=1,\ldots,m,\ \
\zeta=1,\ldots,m,
$
\hfill (0.12)
\\[2.25ex]
where the coefficients $A^{}_{j\zeta\nu}\in C^1(G),\ 
j=1,\ldots,m,\  \zeta=1,\ldots,m,\  \nu=1,\ldots,m.$
\vspace{0.5ex}

If Poisson brackets of operators (0.1) are symmetric, i.e.,
\\[2ex]
\mbox{}\hfill                                           
$
\bigl[{\frak L}^{}_j(x),
{\frak L}^{}_{\zeta}(x)\bigr]=
\bigl[{\frak L}^{}_{\zeta}(x),
{\frak L}^{}_j(x)\bigr]$
\ for all $x\in G,
\quad
j=1,\ldots,m,\ \
\zeta=1,\ldots,m,
$
\hfill (0.13)
\\[2ex]
then we say that the linear homogeneous system of partial differential equations (\!$\partial$\!) with
${\frak L}^{}_j\in C^1(G),\ j=1,\ldots,m,$ is  {\it ja\-co\-bian} [52, p. 62; 92, p. 523].
\vspace{0.5ex}

The symmetry (0.13) of the Poisson brackets 
of  operators (0.1) is equivalent to
\\[2ex]
\mbox{}\hfill                                           
$
\bigl[{\frak L}^{}_j(x),
{\frak L}^{}_{\zeta}(x)\bigr]=
{\frak O}$
\ for all $x\in G,
\quad
j=1,\ldots,m,\ \
\zeta=1,\ldots,m.
$
\hfill(0.14)
\\[2ex]
\indent
The identities (0.14)
\vspace{0.5ex} 
are the identities (0.12) with 
$A^{}_{j\zeta\nu}(x)=0$ for all $x\in G,\,
j\!=1,\ldots,m,$ $\zeta=1,\ldots,m,\ \nu=1,\ldots,m.$
\vspace{0.5ex} 
Therefore the jacobian system (\!$\partial$\!) is complete [52, p. 62].

A differential system
\\[2ex]
\mbox{}\hfill                                           
$
\partial^{}_{x^{}_j} y={\frak M}^{}_j(x) y,
\quad 
j=1,\ldots,m,
$
\hfill{\rm(N\!$\partial$\!)}
\\[2ex]
where $y\in\R,\ x\in\R^n, \ m\leq n ,$ 
the linear differential ope\-ra\-tors of first order
\\[2ex]
\mbox{}\hfill                                            
$
\displaystyle
{\frak M}^{}_j(x)=
\sum\limits_{s=m+1}^n
u^{}_{js}(x)\,\partial^{}_{x^{}_s}$
\ for all $x\in G,
\quad 
j=1,\ldots,m,
$
\hfill(0.15)
\\[2ex]
is called a {\it normal} linear homogeneous 
system of partial differential equations [52, p. 64].
\vspace{0.5ex}

A complete normal partial system is jacobian [52, p. 65; 1, pp. 38 -- 40].
\vspace{0.5ex}

Problems of the theory of integrals for total differential systems were considered 
in process of necessity at the decision of adjacent problems. 
There are first of all the problems of orbits topology for completely solvable autonomous 
total differential systems were studied [93, 94]. 
Directly problems of the theory of integrals for multidimensional 
differential systems (TD), (\!$\partial$\!), and (Pf) are considered in [1, 2].
\\[5.75ex]
\centerline{
\large\bf  
1. First integrals of  total differential system 
}
\\[1.75ex]
\indent
{\bf  1.1. First integral}
\\[0.75ex]
\indent
Suppose the system (TD) has the matrix $X\in C(\Pi),$ 
i.e.,  the entries $X^{}_{ij}\,,\ i=1,\ldots,n,\linebreak 
j=1,\ldots,m,$ 
of the matrix $X$ are continuous functions on the domain $\Pi.$
\vspace{0.75ex}

{\bf Definition 1.1.}
{\it
A scalar function $F\in C^1(\Pi^{\;\!\prime})$
is said to be a 
\textit{\textbf{first integral}} on a domain $\Pi^{\;\!\prime}\subset \Pi$ 
of system {\rm (TD)} with $X\in C(\Pi)$
if the differential of the function $F$ by virtue of system {\rm (TD)} 
vanishes}
on the domain $\Pi^{\;\!\prime}\colon$
\\[2ex]
\mbox{}\hfill                                           
$
dF(t,x)_{\displaystyle |_{\rm(TD)}}=0$
\ for all $(t,x)\in\Pi^\prime.
$
\hfill(1.1)
\\[2.5ex]
\indent
The differential of the function $F$ by virtue of system {\rm (TD)} is
\\[2.25ex]
\mbox{}\hfill
$
\displaystyle
dF(t,x)_{\displaystyle |_{\rm(TD)}}=
\sum_{j=1}^m\,
\partial^{}_{\scriptstyle t^{}_j}
F(t,x)\,
dt^{}_j+
\sum_{i=1}^n\,
\partial^{}_{\scriptstyle x^{}_i}
F(t,x)\,
{dx^{}_i}_{\displaystyle |_{\rm(TD)}}=
\hfill
$
\\[2ex]
\mbox{}\hfill
$
\displaystyle
=\sum_{j=1}^m\,
\biggl(
\partial^{}_{\scriptstyle t^{}_j}
F(t,x)+
\sum_{i=1}^n\,
X^{}_{\scriptstyle ij}(t,x)\,
\partial^{}_{\scriptstyle x^{}_i} 
F(t,x)\biggr) 
dt^{}_j=
\sum_{j=1}^m\,
{\frak X}^{}_j 
F(t,x)\,dt^{}_j$
\ for all $(t,x)\in\Pi^\prime,
\hfill
$
\\[2.5ex]
where ${\mathfrak X}^{}_j\,,\ j=1,\ldots,m,$ are the operators (0.3).
From (1.1) it follows that
\\[2ex]
\mbox{}\hfill                                              
$
{\frak X}^{}_j F(t,x)=0$
\ for all $(t,x)\in\Pi^\prime,
\quad 
j=1,\ldots,m.
$
\hfill(1.2)
\\[2ex]
\indent
The connection between the identity (1.1) and the system of identities (1.2) is 
the jus\-ti\-fi\-ca\-ti\-on of that the linear differential operators (0.3) have been named 
the operators of differentiation by virtue of system (TD).
\vspace{0.75ex}

{\bf Example 1.1.}                                                   
The total differential system 
\\[2ex]
\mbox{}\hfill                                              
$
\begin{array}{c}
dx^{}_1=(x^{}_1 t_1^{-1}+\,t^{}_1 x^{}_2)\,dt^{}_1+\,
t^{}_1 x^{}_2\,dt^{}_2,
\\[2.25ex]
dx^{}_2=({}-1-x^{}_1 t_1^{-1}+\;\!
x_1^2\,t_1^{-2}+\;\!x_2^2)\,dt^{}_1\, +\,
({}-1-x^{}_1 t_1^{-1}+
x_1^2 t_1^{-2}+x_2^2-x^{}_2 x^{}_3)\,dt^{}_2,
\\[2.25ex]
dx^{}_3=x^{}_2(x^{}_3\,dt^{}_1+\,
(x^{}_2+x^{}_3)\,dt^{}_2)
\end{array}
$
\hfill (1.3)
\\[2ex]
induces the linear differential ope\-ra\-tors of first order
\\[2ex]
\mbox{}\hfill                                                
$
{\frak X}^{}_1(t,x)=
\partial^{}_{t^{}_1}+
(x^{}_1 t_1^{-1}+\,
t^{}_1 x^{}_2)\,
\partial^{}_{x^{}_1}+({}-1-x^{}_1 t_1^{-1}+x_1^2\;\! t_1^{-2}+x_2^2)\,
\partial^{}_{x^{}_2}+\,
x^{}_2 x^{}_3\,\partial^{}_{x^{}_3}
$
\hfill (1.4)
\\[1.25ex]
and
\\[1.25ex]
\mbox{}\hfill                                                
$
{\frak X}^{}_2(t,x)=
\partial^{}_{t^{}_2}+\,
t^{}_1 x^{}_2\,\partial^{}_{x^{}_1}+\,
({}-1-x^{}_1 t_1^{-1}+x_1^2 t_1^{-2}+\,
x_2^2-x^{}_2 x^{}_3)\,
\partial^{}_{x^{}_2}+\,
x^{}_2(x^{}_2+x^{}_3)\,
\partial^{}_{x^{}_3}
$
\hfill (1.5)
\\[2.25ex]
on the set $D=\{(t,x)\colon t^{}_1\ne0\}\subset \R^5.$
\vspace{0.5ex}

The operations of operators (1.4) and (1.5) on the scalar function
\\[2ex]
\mbox{}\hfill                                                
$
F\colon (t,x)\to
(1-x_1^2 t_1^{-2}-x_2^2-x_3^2)
\exp({}-2x^{}_1 t_1^{-1})$
\ for all $(t,x)\in D
$
\hfill(1.6)
\\[2.25ex]
are identically equal to zero on the set
$D\colon\  {\frak X}^{}_1 F(t,x)={\frak X}^{}_2 F(t,x)=0$ for all $(t,x)\in D.$
\vspace{0.35ex}

By definition 1.1, the function (1.6) is a first integral on any domain $\Pi\subset D$ 
of the total differential system (1.3).
\vspace{0.75ex}

The direct corollary of Definition 1.1 is
{\sl the criterion of the existence of a first integral for 
a completely solvable total differential system}. 
\vspace{0.35ex}

{\bf Theorem 1.1.}                                            
{\it
A scalar function $F\in C^1(\Pi^\prime)$\vspace{0.35ex} 
is a first integral on the domain $\Pi^\prime\subset\Pi$ of 
the completely solvable on the domain $\Pi^\prime$ system
{\rm(TD)} with $X\in C(\Pi)$ if and only if
the function $F\colon\Pi^\prime\to\R$ is constant along any solution
$x\colon t\to x(t)$ for all $t\in {\rm T}^\prime$ of system {\rm(TD)}, 
where the domain ${\rm T}^{\prime}\subset\R^m$\vspace{0.5ex} is such that 
$(t,x(t))\in\Pi^\prime$ for all $t\in {\rm T}^\prime,$ i.e.,
$F(t,x(t))=C$ for all $t\in {\rm T}^\prime,\ C={\rm const}.$
}
\vspace{0.5ex}

If the system (TD) is not completely solvable, 
then the system (TD) can have first integrals even in the case 
when the system (TD) does not have solutions.
\vspace{0.75ex}

{\bf Example 1.2.}                                                   
The autonomous linear system of total differential equations 
\\[2ex]
\mbox{}\hfill                                                 
$
dx^{}_1=x^{}_1\,dt^{}_1+
3x^{}_1\,dt^{}_2\,,
\quad \ 
dx^{}_2=(1+x^{}_1+2x^{}_2)\,dt^{}_1+
(x^{}_1+3x^{}_2)\,dt^{}_2
$
\hfill (1.7)
\\[2ex]
has the linear differential ope\-ra\-tors of first order
\\[2.25ex]
\mbox{}\hfill
$
{\frak X}^{}_1(t,x)=
\partial^{}_{t^{}_1}+ {\frak A}^{}_1(x)$
\ for all $(t,x)\in\R^4,
\qquad
{\frak X}^{}_2(t,x)=
\partial^{}_{t^{}_2}+{\frak A}^{}_2(x)$
\ for all $(t,x)\in\R^4,
\hfill
$
\\[2ex]
where
\\[1ex]
\mbox{}\hfill
$
{\frak A}^{}_1(x)=
x^{}_1\,\partial^{}_{x^{}_1}+ (1+x^{}_1+2x^{}_2)\,\partial^{}_{x^{}_2}$
\ for all $x\in\R^2,
\hfill
$
\\[2.5ex]
\mbox{}\hfill
$
{\frak A}^{}_2(x) =
3x^{}_1\,\partial^{}_{x^{}_1}+(x^{}_1+3x^{}_2)\,\partial^{}_{x^{}_2}
$
\ for all $x\in\R^2.
\hfill
$
\\[2.25ex]
\indent
The Poisson bracket for the autonomous operators of differentiation by virtue of 
the total differential sys\-tem (1.7) 
\\[2ex]
\mbox{}\hfill
$
\bigl[{\frak A}^{}_1(x),
{\frak A}^{}_2(x)\bigr]=
(3-x^{}_1)\,\partial_{x^{}_2}$
\ for all $x\in\R^2
\hfill
$
\\[2ex]
is not the null operator on any two-dimensional domain from the phase plane $\R^2.$
By the Frobenius theorem (Theorem 0.1), the system (1.7) 
is not completely solvable.

The operations of the nonautonomous operators of differentiation by virtue of sys\-tem (1.7)
on the holomorphic scalar function
\\[2ex]
\mbox{}\hfill                                                
$
F\colon (t,x)\to
x^{}_1\exp({}-(t^{}_1+3t^{}_2))$
\ for all $(t,x)\in\R^4
$
\hfill(1.8)
\\[2.25ex]
are identically equal to zero on the space 
\vspace{0.35ex}
$\R^4\colon\ {\frak X}^{}_1 F(t,x)\!={\frak X}^{}_2 F(t,x)\!=0$ for all $(t,x)\in\R^4.$
By de\-fi\-ni\-ti\-on 1.1, the function (1.8) is a first integral on the space $\R^4$ of system (1.7).

Let us prove that the system (1.7) has no solutions. 
\vspace{0.5ex}

Since $x^{}_1\colon t\to C\exp(t^{}_1+3t^{}_2)$ for all $t\in\R^2,$
\vspace{0.35ex}
we see that the first equation of system (1.7) is 
an identity on the plane $\R^2$ and 
the second equation of system (1.7) is 
\\[2ex]
\mbox{}\hfill                                                 
$
dx^{}_2=P(t,x^{}_2)\,dt^{}_1+
Q(t,x^{}_2)\,dt^{}_2\,,
$
\hfill(1.9)
\\[1.25ex]
where
\\[1ex]
\mbox{}\hfill
$
P\colon (t,x^{}_2)\to
1+2x^{}_2+C\exp(t^{}_1+3t^{}_2)$
\ for all $(t,x^{}_2)\in\R^3,
\hfill
$
\\[2.5ex]
\mbox{}\hfill
$
Q\colon (t,x^{}_2)\to
3x^{}_2+C\exp(t^{}_1+3t^{}_2)$
\ for all $(t,x^{}_2)\in\R^3,
\hfill
$
\\[1.75ex]
$C$ is a constant from the field $\R.$

If the equation (1.9) possess the solution
$x^{}_2\colon t\to x^{}_2(t)$ for all $t\in{\rm T},$ 
where ${\rm T}$ is some domain of the plane $\R^2,$
then the following conditions hold
\\[2ex]
\mbox{}\hfill
$
\partial^{}_{t^{}_1}\,x^{}_2(t)=P(t,x^{}_2(t))$
\ for all $t\in{\rm T},
\qquad 
\partial^{}_{t^{}_2}\,x^{}_2(t)=Q(t,x^{}_2(t))$
\ for all $t\in{\rm T}.
\hfill
$
\\[2.25ex]
\indent
But such the function $x^{}_2\colon{\rm T}\to\R$ is not exist.
It follows that the mixed derivatives of the second order
\\[2ex]
\mbox{}\hfill
$
\partial^{}_{t^{}_1 t^{}_2}\,x^{}_2(t)=
\partial^{}_{t^{}_2}\,P(t,x^{}_2(t))=
6x^{}_2(t)+5C\exp(t^{}_1+3t^{}_2)$
\ for all $t\in{\rm T},
\hfill
$
\\[2.25ex]
\mbox{}\hfill
$
\partial^{}_{t^{}_2 t^{}_1}\,x^{}_2(t)=
\partial^{}_{t^{}_1}\,Q(t,x^{}_2(t))=
3+6x^{}_2(t)+4C\exp(t^{}_1+3t^{}_2)$
\ for all $t\in{\rm T}
\hfill
$
\\[2ex]
are not coincide neither at any $C$ from the field $\R$ nor on any domain of the plane $\R^2.$
 
Thus the not completely solvable total differential system (1.7) has the first integral (1.8), 
but the system (1.7) doesn't have solutions.
\\[2.5ex]
\indent
{\bf  1.2. Basis of first integrals}
\\[0.5ex]
\indent
Consider the set of scalar functions 
\\[1.5ex]
\mbox{}\hfill         
$
F^{}_s\colon (t, x)\to F^{}_s(t, x)$
for all $\!(t,x)\in\Pi^\prime,
\ \ \
F^{}_s\subset C^1(\Pi),\
s=1,\ldots,k,\  
\Pi^\prime\subset\Pi\subset\R^{m+n},
$
\hfill (1.10)
\\[1.5ex]
and form the vector function
\\[1.5ex]
\mbox{}\hfill              
$
F\colon(t,x)\to \bigl(F^{}_1(t,x),\ldots, F^{}_k(t,x)\bigr)$
\ for all $(t,x)\in\Pi^\prime
$
\hfill (1.11)
\\[1.5ex]
with range ${\rm E}F\subset \R^k.$
\vspace{0.5ex}

{\bf Theorem 1.2.}
{\it
Suppose the functions {\rm (1.10)} are first integrals 
on the domain $\Pi^{\;\!\prime}\subset \Pi$ of system {\rm (TD)} with $X\in C(\Pi).$ 
Then the function
\\[1.5ex]
\mbox{}\hfill
$
\Psi\colon(t,x)\to \Phi\bigl(F^{}_1(t,x),\ldots, F^{}_k(t,x)\bigr)$
\ for all $(t,x)\in\Pi^\prime, 
$
\hfill {\rm (1.12)}
\\[1.75ex]
where arbitrary scalar function $\Phi\in C^1_{}({\rm E}F),$ 
is a first integral on the domain $\Pi^{\;\!\prime}$ of the 
total differential system {\rm (TD)}.
}

{\sl Proof}. From Definition 1.1 for the first integrals (1.10), we have
\\[1.5ex]
\mbox{}\hfill
$
{\frak X}^{}_j\,F^{}_s(t,x)=0$
\ for all $(t,x)\in\Pi^\prime,
\ \  
j=1,\ldots,m,\ 
s=1,\ldots,k.
\hfill
$
\\[1.5ex]
\indent
Then, for any scalar function $\Phi\!\in\! C^1({\rm E}\,F)$ on the range $\!{\rm E}\,F\!$ of the 
vector function (1.11), we obtain 
\\[1.5ex]
\mbox{}\hfill
$
\displaystyle
{\frak X}^{}_j\,
\Phi\bigl(F(t,x)\bigr)=
\sum\limits_{s=1}^{k}\,
\partial^{}_{F^{}_s}\,
\Phi(F)_{\displaystyle 
|_{\scriptstyle F=F(t,x)}}\,
{\frak X}^{}_j\,F^{}_s(t,x)=0$
\ for all $(t,x)\in\Pi^\prime,
\ \ 
j=1,\ldots,m.
\hfill
$
\\[1.5ex]
\indent
By Definition 1.1, the function (1.12) is a first integral on the domain $\!\Pi^\prime\!$
\vspace{0.5ex}
of  system (TD).\k

This theorem expresses functional ambiguity of a first integral for a 
system of total differential equations: 
{\it 
if a scalar function $F^{}_1\!\in\! C^1(\Pi^\prime)$ 
is a first integral on the domain $\Pi^\prime\!\subset\!\Pi$ of system 
{\rm(TD)} with $X\in C(\Pi),$ then the function 
$\Psi^{}_1\colon(t,x)\to \Phi\bigl(F^{}_1(t,x)\bigr)$
for all $(t,x)\in\Pi^\prime,\!\!$ 
where arbitrary scalar function $\!\Phi\!\in\!\! C^1_{}\!({\rm E}F_1\!),\!\!$
\vspace{0.5ex}
is a first integral on the domain $\!\Pi^\prime\!\!$ of system {\rm(\!TD\!)}\!.}

Thus the priority of first integrals functionally independent on some domain is installed.
Therefore for a system of total differential equations  we have the problems about the existence and 
the number of functionally independent first integrals.
\vspace{0.5ex}

{\bf Definition 1.2.}
{\it
A set of the functionally independent on the domain $\Pi^{\;\!\prime}\subset \Pi$
first integrals {\rm (1.10)} of system {\rm (TD)} with
$X\in C(\Pi)$ is called a 
\textit{\textbf{basis of first integrals}} 
on the domain $\Pi^{\;\!\prime}$ of system {\rm (TD)} if
for any first integral $\Psi$ on the domain $\Pi^{\;\!\prime}$ of system {\rm (TD)},
we have $\Psi(t,x) = \Phi(F_1^{}(t,x),\ldots, F_k^{}(t,x))$ for all $(t,x) \in \Pi^{\;\!\prime},$
where $\Phi$ is some function of class $C^{1}_{}({\rm E}F),\ 
{\rm E}F$ is the range of the vector function  {\rm (1.11)}.
The number k is said to be the \textit{\textbf{dimension}} of 
ba\-sis of first integrals on the domain $\Pi^{\;\!\prime}$ of system {\rm (TD)}.
}
\vspace{0.5ex}

A basis of first integrals we'll name also as an {\sl integral basis}.
\vspace{0.5ex}

{\bf Definition 1.3.}
{\it
We'll say that two systems of total differential equations
are \textit{\textbf{integrally equivalent}} on some domain if  on this domain
each first integral of the first system is a first integral of the second system 
and on the contrary each first integral of the second system 
is a first integral of the first system.
}
\vspace{0.5ex}

The integrally equivalent on the domain $\Pi^{\;\!\prime}$ 
total differential systems have the same integral basis on this domain.
\\[2.5ex]
\indent
{\bf                                                     
1.3. Dimension of  ba\-sis of first integrals for completely solvable systems}
\\[1ex]
\indent
Suppose the system (TD) has the matrix $X\in C^\infty(\Pi),$ 
i.e.,  the entries $X^{}_{ij}\,,\ i=1,\ldots,n,\linebreak 
j=1,\ldots,m,$ 
of the matrix $X$ are holomorphic functions on the domain $\Pi.$
${\rm T}$ is a domain from the space $\R^m$
such that for any solution $x\colon t\to x(t)$ for all $t\in{\rm T}$ of 
the completely solvable system {\rm(TD)}, it follows that
$(t,x(t))\in\Pi$ for all $t\in{\rm T}.$ 
Then, by the Cauchy theorem (see, for example, [93, p. 26]),
for any point $t^{}_0\in{\rm T}$ there exists a neighbourhood 
${\rm U}^{}_0$ such that the solution 
$x\colon t\to x(t)$ for all $t\in {\rm U}^{}_0$ of the 
completely solvable system {\rm(TD)} with $X\in C^\infty(\Pi)$ is
a holomorphic function.
Moreover, solutions of the completely solvable total differential system {\rm(TD)} with $X\in C^\infty(\Pi)$ 
depends holomorphically on parameters and initial data [93, pp.~39~-- 41].

By $x\colon t\to x(t;(t^{}_0\,,x^{}_0))$ for all $t\in{\rm T}$ we denote
the solution $x\colon t\to x(t)$ for all $t\in{\rm T}$ satisfying the initial condition 
$x(t^{}_0)=x^{}_0.$
\vspace{0.75ex}

{\bf Lemma 1.1.}                                     
{\it
Let the vector function 
\vspace{0.25ex}
$x\colon t\to  x(t;(t^{}_0\,,x^{}_0))$ for all $t, t^{}_0\in{\rm T}^{\;\!\prime}$ be
a solution on some simply connected domain
${\rm T}^{\,\prime}\subset{\rm T}$ 
\vspace{0.25ex}
of the completely solvable system {\rm(TD)} with $X\in C^\infty(\Pi).$
\vspace{0.35ex}
Then for any $t^{}_\ast\in{\rm T}^\prime$ there exists the solution
$x\colon t\to x(t;(t^{}_\ast\,,x^{}_\ast))$ for all $t\in{\rm T}^\prime$
of system {\rm(TD)}, where $x^{}_\ast=x(t^{}_\ast;(t^{}_0\,,x^{}_0)),$ such that
\vspace{0.75ex}
$x(t^{}_0;(t^{}_\ast\,,x(t^{}_{\ast};(t^{}_0\,,x^{}_0))))=x^{}_0\,.$
}

{\sl Proof.}\! Let $\widetilde{x},\ \widehat{x}$ be the solutions of 
the corresponding Cauchy's problems of system {\rm(\!TD\!)}:
\\[2ex]
\mbox{}\hfill
$
\widetilde{x}\colon t\to x(t;(t^{}_0\,,x^{}_0))$
\ for all $t\in{\rm T}^\prime
\hfill
$
\\[0.5ex]
and
\\[1ex]
\mbox{}\hfill
$
\widehat{x}\colon t\to
x(t;(t^{}_{\ast}\,,
x(t^{}_{\ast};(t^{}_0\,,x^{}_0))))$
\ for all $t\in{\rm T}^\prime.
\hfill
$
\\[2.25ex]
The values $\widetilde{x}(t^{}_{\ast})=\widehat{x}(t^{}_{\ast}).$
Then, by the Cauchy theorem, we have
$\widetilde{x}(t)=\widehat{x}(t)$ for all $t\in{\rm T}^\prime,$
where a simply connected domain ${\rm T}^\prime$ such that 
the points $t^{}_0$ and $t^{}_{\ast}$ belongs to ${\rm T}^\prime.$
\vspace{0.5ex}

Thus $\widehat{x}(t^{}_0)=\widetilde{x}(t^{}_0)=x^{}_0,$ i.e.,
$x(t^{}_0;(t^{}_{\ast}\,,x(t^{}_{\ast};(t^{}_0\,,x^{}_0))))=x^{}_0.\ \k$
\vspace{1.25ex}

{\bf Lemma 1.2.}\!                                  
\vspace{0.25ex}
{\it
Suppose the completely solvable system \!{\rm(TD)}\! with $\!\!X\!\in\! C^\infty(\Pi)\!$ in 
a ne\-ig\-h\-bo\-ur\-hood of the point $\!(t^{}_0,x^{}_0)\!\in\!\Pi\!$ satisfies the conditions of 
the Cauchy theorem. Then this system has $n$ 
functionally independent on some neighbourhood of the point $(t^{}_0\,,x^{}_0)$
\vspace{0.5ex}
first in\-teg\-rals.
}

{\sl Proof.} 
\vspace{0.35ex}
Let $x\colon t\to x(t;(t^{}_0\,,x^{}_0))$ for all $t\in{\rm T}^\prime$
be a solution of system {\rm(TD)} on a simply connected domain
${\rm T}^\prime\ni t^{}_0\,.$ The function
\\[2ex]
\mbox{}\hfill
$
F\colon(t,x)\to x(t^{}_0;(t,x))$
\ for all $(t,x)\in{\rm U}^{}_{00}\,,
\hfill
$
\\[2ex]
where ${\rm U}^{}_{00}$ is some neighbourhood of the point $(t^{}_0\,,x^{}_0),$ 
under
\\[2ex]
\mbox{}\hfill
$
\displaystyle
x(t;(t^{}_0\,,x^{}_0))=x^{}_0+\sum\limits_{k=1}^{\infty}\, 
a^{}_k(t-t^{}_0)^k$
\ for all $t\in{\rm U}^{}_0\,,
\hfill
$
\\[2.25ex]
where ${\rm U}^{}_0$ is some neighbourhood of the point $t^{}_0\,,$ is such that
\\[2ex]
\mbox{}\hfill
$
\displaystyle
F(t,x)=x(t^{}_0;(t,x))=x+\sum\limits_{k=1}^{\infty}\,
a^{}_k(t^{}_0-t)^k$
\ for all $(t,x)\in{\rm U}^{}_{00}\,.
\hfill
$
\\[2ex]
\indent
Since solutions of the Cauchy problem depends holomorphically on initial data, 
we see that the function $F$ is holomorphic on the neighbourhood ${\rm U}^{}_{00}.$
\vspace{0.5ex}
The Jacobi matrix at the point $(t^{}_0\,,x^{}_0)$ 
with respect to $x$ is the identity matrix, i.e., 
$\partial^{}_x\,F(t,x)_{\displaystyle 
|_{\scriptstyle (t^{}_0\,,x^{}_0)}}={\rm E}.$
Therefore Jacobian
${\rm det}\,\partial^{}_x\,F(t,x)\ne 0$ 
for all $(t,x)\in{\rm U}^{}_{00}\,.$
\vspace{0.5ex}

This implies that the coordinate functions 
\vspace{0.25ex}
$F^{}_i\colon{\rm U}^{}_{00}\to\R,\ i=1,\ldots,n,$ 
of the vector function $F$ are functionally independent of $x$ on 
the neighbourhood ${\rm U}^{}_{00}.$
\vspace{0.35ex}

By Lemma 1.1,
\\[1.25ex]
\mbox{}\hfill
$
F(t,x(t;(t^{}_0\,,x^{}_{\ast})))=
x(t^{}_0;(t,x(t;(t^{}_0\,,x^{}_{\ast}))))=
x^{}_{\ast}\,.
\hfill
$
\\[1.75ex]
Hence the function $F$ is a constant vector
along a solution of system {\rm(TD)}.
Then, by The\-o\-rem~1.1, the scalar functions
\\[1.5ex]
\mbox{}\hfill                                        
$
F^{}_i\colon(t,x)\to F^{}_i(t,x)$
\ for all $(t,x)\in{\rm U}^{}_{00}\,,
\ \
i=1,\ldots,n,
$
\hfill(1.13)
\\[1.75ex]
are first integrals on the neighbourhood ${\rm U}^{}_{00}$ of system {\rm(TD)}.\ \k
\vspace{0.75ex}

{\bf Lemma 1.3.}                                      
{\it
Suppose the completely solvable system {\rm(TD)} with $\!X\!\in C^\infty(\Pi)$
has $n$ fun\-c\-ti\-o\-nal\-ly independent on a neighbourhood ${\rm U}^{}_{00}$
\vspace{0.25ex}
of the point $(t^{}_0\,,x^{}_0)\in\Pi$ first integrals {\rm(1.13)}.
Then for any first integral $\Psi\colon{\rm U}^{}_{00}\to\R$ of system {\rm (TD)},
\vspace{0.25ex}
we have $\Psi(t,x)=\Phi(F(t,x))$ for all $(t,x)\in{\rm U}^{}_{00}\,,$
\vspace{0.35ex}
where $\Phi$ is some scalar holomorphic function on the range ${\rm E}F$ of 
the vector function
$F\colon (t,x)\to \bigl(F^{}_1(t,x),\ldots, F^{}_n(t,x)\bigr)$ for all $(t,x)\in{\rm U}^{}_{00}.$
}
\vspace{0.5ex}

{\sl Proof.}
The first integrals (1.13) of system (TD) has the form
$F^{}_i\colon(t,x)\to x^{}_i(t^{}_0;(t,x))$ for all $(t,x)\in{\rm U}^{}_{00}\,,\ i=1,\ldots,n.$
\vspace{0.35ex}
Jacobian ${\rm det}\,\partial^{}_x\, F(t,x)\ne 0$ for all $(t,x)\in {\rm U}^{}_{00}\,.$
Therefore the function $F$ at fixed $t$ has the inverse function $S$ and
\\[1.75ex]
\mbox{}\hfill                                             
$
F(t,S(t,x))=x$
\ for all $(t,x)\in{\rm U}^{}_{00}\,.
$
\hfill(1.14)
\\[2ex]
\indent
The function $\Phi\colon(t,x)\to\Psi(t,S(t,x))$ for all $(t,x)\in{\rm U}^{}_{00}$ 
with the first integrals (1.13) is connected by the identity
\\[1.5ex]
\mbox{}\hfill
$
\Psi(t,x)=\Phi(t,F(t,x))$
\ for all $(t,x)\in{\rm U}^{}_{00}\,.
\hfill
$
\\[1.75ex]
\indent
Let us prove that the function $\Phi$ is independent of $t$ on 
the neighbourhood ${\rm U}^{}_{00}\colon$
\\[2ex]
\mbox{}\hfill
$
\partial^{}_{t^{}_j}\Phi(t,x)=0$
\ for all $(t,x)\in{\rm U}^{}_{00}\,,
\ \
j=1,\ldots,m.
\hfill
$
\\[2ex]
\indent
Differentiating the identity (1.14) with respect to $t,$ we get
\\[2ex]
\mbox{}\hfill
$
\partial^{}_{t^{}_j} F(t,S(t,x))+
\partial^{}_x\,F(t,S(t,x))\,
\partial^{}_{t^{}_j} S(t,x)=0$
\ for all $(t,x)\in{\rm U}^{}_{00}\,,
\ \
j=1,\ldots,m.
\hfill
$
\\[2ex]
\indent
The functions (1.13) are first integrals of system (TD). Hence,
\\[1.75ex]
\mbox{}\hfill
$
\partial^{}_{t^{}_j} F(t,x)=
{}-\partial^{}_x\,F(t,x)\,
X^{{}^{\scriptstyle j}}(t,x)$
\ for all $(t,x)\in{\rm U}^{}_{00}\,,
\ \
j=1,\ldots,m,
\hfill
$
\\[2ex]
where the vector functions 
\vspace{0.5ex} 
$\!X^{{}^{\scriptstyle j}}\!\colon\! (t,x)\!\to\!
\bigl(X^{}_{1j}(t,x),\ldots, X^{}_{nj}(t,x)\bigr)\!$ for all $\!(t,x)\!\in\!\Pi,\,  j\!=\!1,\ldots,m.\!\!$

Then
\\[1.5ex]
\mbox{}\hfill
$
\partial^{}_x\,F(t,S(t,x))
\bigl(\partial^{}_{t^{}_j}S(t,x)-
X^{{}^{\scriptstyle j}}
(t,S(t,x))\bigr)=0$ 
\ for all $(t,x)\in{\rm U}^{}_{00}\,,
\ \
j=1,\ldots,m.
\hfill
$
\\[1.75ex]
\indent
Therefore,
\\[1.25ex]
\mbox{}\hfill
$
\partial^{}_{t^{}_j}S(t,x)=X^{{}^{\scriptstyle j}}(t,S(t,x))$
\ for all $(t,x)\in{\rm U}^{}_{00},
\ \
j=1,\ldots,m,
\hfill
$
\\[2ex]
since $\partial^{}_x\,F$ is a nonsingular matrix on ${\rm U}^{}_{00}\,.$
\vspace{0.5ex}

Taking into account 
the function $\Psi$ is a first integral of system (TD), we obtain
\\[2ex]
\mbox{}\hfill
$
\partial^{}_{t^{}_j}\Phi(t,x)=
\partial^{}_{t^{}_j}\Psi(t,S(t,x))+
\partial^{}_x\,\Psi(t,S(t,x))\,
\partial^{}_{t^{}_j}S(t,x)={}
\hfill
$
\\[1.75ex]
\mbox{}\hfill
$
=\partial^{}_{t^{}_j}\Psi(t,S(t,x))+
\partial^{}_x\,\Psi(t,S(t,x))\,
X^{{}^{\scriptstyle j}}(t,S(t,x))=0$
\ for all $(t,x)\in{\rm U}^{}_{00}\,,
\ 
j=1,\ldots,m. \ \k
\hfill
$
\\[1.75ex]
\indent
From lemmas 1.2 and 1.3, we have the  following 
\vspace{0.75ex}

{\bf Theorem 1.3.}                     
{\it
The completely solvable system {\rm(TD)} with $X\in C^\infty(\Pi)$
on a neighbourhood of any point of the domain $\Pi$ 
has a basis of first integrals of dimension $n.$
}
\vspace{0.75ex}

{\bf Example 1.3.}                                    
The completely solvable autonomous total differential system 
\\[1.25ex]
\mbox{}\hfill                                        
$
dx^{}_1=dt^{}_1,
\quad\
dx^{}_2=dt^{}_2,
\quad\
dx^{}_3=
\partial^{}_{x^{}_1}
g(x^{}_1,x^{}_2)\,dt^{}_1+ 
\partial^{}_{x^{}_2}
g(x^{}_1,x^{}_2)\,dt^{}_2,
$
\hfill (1.15)
\\[2ex]
where the scalar function $g$ is holomorphic on a domain
\vspace{0.25ex} 
${\mathscr D}\subset \R^2,$ 
has the basis of first integrals on the simply connected domain
$\Pi^\prime=\R^2\times{\mathscr D}\times\R$
\\[1.75ex]
\mbox{}\hfill                                        
$
\begin{array}{c}
F^{}_1\colon(t,x)\to\, t^{}_1-x^{}_1$
\ for all $(t,x)\in\Pi^{\;\!\prime},
\quad
F^{}_2\colon(t,x)\to\, t^{}_2-x^{}_2$
\ for all $(t,x)\in\Pi^{\;\!\prime},
\\[2ex]
F^{}_3\colon(t,x)\to\,
g(x^{}_1,x^{}_2)-x^{}_3$
\ for all $(t,x)\in\Pi^{\;\!\prime}.
\end{array}
$
\hfill (1.16)
\\[2ex]
\indent
Indeed, by Definition 1.1, the functions (1.16) are first integrals on the domain $\Pi^\prime$
of system (1.15). The Jacobi matrix $J(F^{}_1,F^{}_2,F^{}_3\;\!; t,x)$ has 
${\rm rank}\,J(F^{}_1,F^{}_2,F^{}_3;t,x)=3$ for all $(t,x)\in\Pi^\prime.$
Thus the first integrals (1.16) are functionally independent on the domain $\Pi^\prime.$

Therefore, by Theorem 1.3, the first integrals (1.16) are an
integral basis on the domain $\Pi^\prime$ of 
the completely solvable autonomous total differential system (1.15).
\\[5.5ex]
\centerline{
\large\bf  
2. First integrals for linear homogeneous system 
}
\\[0.5ex]
\centerline{
\large\bf  
of partial differential equations}
\\[1.5ex]
\indent
{\bf  2.1. Basis of first integrals}
\\[0.5ex]
\indent
{\bf Definition 2.1.}
{\it
A scalar function $F\in C^1(G^{\;\!\prime})$
is said to be a 
\textit{\textbf{first integral}} on a domain $G^{\;\!\prime}\subset G$ 
of system {\rm (\!$\partial$\!)} with
${\frak L}_j^{}\in C(G),\, j=1,\ldots,m,$ if
\\[1.5ex]
\mbox{}\hfill                                   
$
{\frak L}_j^{}\;\! F(x) = 0$
\ for all $x \in G^{\;\!\prime},
\ \ j =1,\ldots,m.
$
\hfill {\rm (2.1)}
\\[2ex]
}
\indent
Using $k$ scalar functions 
\\[1.5ex]
\mbox{}\hfill         
$
F_{s}^{} \colon x \to F_{s}^{}(x)$
\ for all $x \in G^{\;\!\prime},
\ \ s = 1,\ldots, k,
\hfill 
$
\\[-0.25ex]
\mbox{}\hfill (2.2)
\\[-0.25ex]
\mbox{}\hfill
$
F^{}_s\in C^1(G^{\;\!\prime}),
\ \ s = 1,\ldots, k,
\ \ \ G^{\;\!\prime}\subset G,
\hfill
$
\\[1.5ex]
we form the vector function
\\[1.5ex]
\mbox{}\hfill              
$
F\colon x \to (F_{1}^{}(x),\ldots,F_{k}^{}(x))$
\ for all $x \in G^{\;\!\prime}
$
\hfill (2.3)
\\[1.5ex]
with range ${\rm E}F\subset \R^k.$
\vspace{0.5ex}

{\bf Theorem 2.1.}
{\it
Let the functions {\rm (2.2)} be first integrals 
on a domain $G^{\;\!\prime}\subset G$ of system {\rm (\!$\partial$\!)} with 
${\frak L}_j^{}\in C(G),\ j=1,\ldots,m.$ 
Then the function
\\[1.5ex]
\mbox{}\hfill
$
\Psi \colon x \to \Phi(F^{}_1(x),\ldots, F^{}_k(x))$
\ for all $x \in G^{\;\!\prime},
$
\hfill {\rm (2.4)}
\\[1.75ex]
where arbitrary function $\Phi\in C^1_{}({\rm E}F),$ 
\vspace{0.25ex}
is a first integral on the domain $G^{\;\!\prime}$ of system {\rm (\!$\partial$\!)}.
}

{\sl Proof}. By Definition 2.1, we have
\\[1.5ex]
\mbox{}\hfill
$
{\frak L}_{_{\scriptstyle j}} F_{s}^{}(x) = 0$
\ for all $x \in G^{\;\!\prime},
\ \ j=1,\ldots,m, \ s= 1,\ldots, k.
\hfill
$
\\[1.5ex]
\indent
Then
\\[1.5ex]
\mbox{}\hfill
$
\displaystyle
{\frak L}_{{}_{\scriptstyle j}} \Phi(F(x)) =
\sum\limits_{s = 1}^k\,
\partial_{{}_{\scriptstyle F_{{}_{\scriptsize s}}}}
\Phi(F)_{\displaystyle |_{F=F(x)}}
{\frak L}_{{}_{\scriptstyle j}}
F_{{}_{\scriptstyle s}}(x) = 0$
\ for all $x\in G^{\;\!\prime},
\ \ 
j =1,\ldots, m.
\hfill
$
\\[1.5ex]
\indent
Therefore  the composite scalar function (2.4) is 
a first integral of  system (\!$\partial$\!).\ \k
\vspace{0.5ex}

This theorem expresses functional ambiguity of first integrals for a 
linear homogeneous system of partial differential equations (see Subsection 1.2).
\vspace{0.5ex}

{\bf Definition 2.2.}
{\it
A set of functionally independent first integrals on a domain 
$G^{\;\!\prime}\subset G$ of system {\rm (\!$\partial$\!)} with
${\frak L}_j^{}\in C(G),\ j=1,\ldots,m,$ is called a 
\textit{\textbf{basis of first integrals}} {\rm(}\textit{\textbf{integral basis}}{\rm)}
on the domain $G^{\;\!\prime}$ of system {\rm (\!$\partial$\!)} if
for any first integral $\Psi$ on the domain $G^{\;\!\prime}$ of system {\rm (\!$\partial$\!)},
we have $\Psi(x) = \Phi(F(x))$ for all $x \in G^{\;\!\prime},$
where $\Phi$ is some function of class $C^{1}_{}({\rm E}F),\ 
{\rm E}F$ is the range of function  {\rm (2.3)}.
The number k is said to be the \textit{\textbf{dimension}} of 
ba\-sis of first integrals on the domain $G^{\;\!\prime}$ of system {\rm (\!$\partial$\!)}.
}
\vspace{0.5ex}

Suppose $m=n.$ Then the coefficient matrix 
\\[1.5ex]
\mbox{}\hfill                                   
$
u(x) = \|u_{_{\scriptstyle ji}}(x)\|$
\ for all $x \in G
$
\hfill (2.5)
\\[1.5ex]
of system {\rm (\!$\partial$\!)} is a square matrix of order $n.$
Since the operators (0.1) are not linearly bound on the domain $G,$
we see that the matrix (2.5) nearly everywhere on the domain $G$ is nonsingular.
Let $G^{\;\!\prime}\subset G$ be a domain such that 
${\rm rank}\,u(x) = n$ for all $x\in G^{\;\!\prime}.$ 
Then the system {\rm (\!$\partial$\!)} on the domain $G^{\;\!\prime}$ 
is equivalent to the system of $n$ differential equations 
\\[1.5ex]
\mbox{}\hfill
$
\partial_{_{\scriptstyle x_i}} y = 0, \  \, i = 1,\ldots, n.
\hfill
$
\\[1.5ex]
\indent
The function 
$y\colon x\to C$ for all $x\in G^{\;\!\prime},$
where $C$ is arbitrary real constant, is a 
first integral on the domain $G^{\;\!\prime}$ of this system.

Thus we obtain the following
\vspace{0.5ex}

{\bf Property 2.1.}
{\it
If $m=n,$ then all first integrals of system {\rm (\!$\partial$\!)} are identical constants.
}
\vspace{0.5ex}

Therefore the basic object of our research is the system $(\partial)$ with $m <n.$
\\[2.25ex]
\indent
{\bf  2.2. Incomplete system}
\\[0.75ex]
\indent
Suppose the coordinate functions 
$u^{}_{ji}\,,\ j=1,\ldots,m,\ i=1,\ldots,n,$ of 
the linear differential ope\-ra\-tors (0.1) 
are a sufficient number of times continuously differentiable 
or holomorphic on the domain $G,$ i.e.,  
${\frak L}^{}_j\in C^k(G)$ or
${\frak L}^{}_j\in C^\infty(G),\ j=1,\ldots,m,$ where $k\in\N.$
\vspace{0.5ex}

{\bf Lemma 2.1.}                                             
{\it
If the function $y\in C^2(G^\prime)$ is a first integral
on a domain $G^{\;\!\prime}\subset G$ of the system of equations
\\[0.25ex]
\mbox{}\hfill
$
{\frak L}^{}_1(x)\,y=0,
\quad \ \  
{\frak L}^{}_2(x)\,y=0,
\hfill
$
\\[1.5ex]
then this function is a first integral on the domain $G^{\;\!\prime}$ of the
linear homogeneous first-order par\-ti\-al differential equation  
\\[0.5ex]
\mbox{}\hfill
$
\bigl[{\frak L}^{}_1(x),
{\frak L}^{}_2(x)\bigr]y=0.
\hfill
$
\\[2ex]
}
\indent
{\sl Proof} [92, pp. 520 -- 521].
The operations of Poisson brackets and commutators [91, p. 165]
on twice continuously differentiable functions coincide.
This implies that 
\\[1.75ex]
\mbox{}\hfill
$
\bigl[{\frak L}^{}_1\,,
{\frak L}^{}_2\bigr]y(x)=
{\frak L}^{}_1{\frak L}^{}_2\,y(x)-
{\frak L}^{}_2{\frak L}^{}_1\,y(x)=0$
\ for all $x\in G^{\;\!\prime}.\ \ \k
\hfill
$
\\[2ex]
\indent
The following lemma is an immediate consequence of Lemma 2.1.
\vspace{0.5ex}

{\bf Lemma 2.2.}                                            
{\it
If the function $F\in C^2(G^{\;\!\prime})$ is a first integral
on a domain $G^{\;\!\prime}\subset G$ of  sys\-tem $(\partial)$ with
$\!{\frak L}^{}_j\in C^1(G),\,  j=1,\ldots,m,\!$ 
then this function is a first integral on the do\-ma\-in $G^{\;\!\prime}$ of the
linear homogeneous system of par\-ti\-al differential equations}
\\[2ex]
\mbox{}\hfill                                            
$
{\frak L}^{}_{j}(x)\,y=0,
\ \
j=1,\ldots,m, 
\quad \ \,
\bigl[{\frak L}^{}_{j^{}_\nu}(x),
{\frak L}^{}_{l_\mu}(x)\bigr]\;\!y=0,
\hfill
$
\\[-0.25ex]
\mbox{}\hfill (2.6)
\\
\mbox{}\hfill
$
\nu=1,\ldots,m^{}_1\,,\ \ 
\mu=1,\ldots,m^{}_2\,,\ \ \ 
m^{}_1\leq m,\ 
m^{}_2\leq m,\ \,
j^{}_\nu\,,\ 
l^{}_\mu\in\{1,\ldots,m\}.
\hfill
$
\\[2.5ex]
\indent
Let the system $(\partial)$ be complete. 
If to add to the system $(\partial)$ at least one 
partial differential equation of the form
\\[1.75ex]
\mbox{}\hfill                                            
$
\bigl[{\frak L}^{}_{j^{}_\nu}(x),
{\frak L}^{}_{l_\mu}(x)\bigr]y=0,
\quad
j^{}_{\nu}\,,\ 
l^{}_{\mu}\in \{1,\ldots,m\},
$
\hfill(2.7)
\\[2ex]
then the system of par\-ti\-al differential equations (2.6) is built
on the basis of the linearly bound on the domain $G$ operators.

If the system $(\partial)$ is incomplete, then some operators
\\[2ex]
\mbox{}\hfill
$
{\frak L}^{}_{\tau\theta}(x)=
\bigl[{\frak L}^{}_\tau(x),
{\frak L}^{}_\theta(x)\bigr],
\  \ \,
\tau,\, 
\theta\in\{1,\ldots,m\},
\hfill
$
\\[2.25ex]
are not a linear combination of the operators
${\frak L}^{}_j\,,\  j=1,\ldots,m.$
Adding the equation ${\frak L}^{}_{\tau\theta}(x)\,y\!=0$ to 
the system $(\partial),$ we obtain 
the linear homogeneous par\-ti\-al differential system
\\[2ex]
\mbox{}\hfill                                         
$
{\frak L}^{}_s(x)\,y=0,
\quad 
s=1,\ldots,k^{}_1\,,\ \
m<k^{}_1<n,
$
\hfill (2.8)
\\[2.25ex]
where the operators ${\frak L}^{}_s,\ s=1,\ldots,k^{}_1\;\!,$ 
are not linearly bound on the domain $G.$  

It is important to underline, that

1) $k^{}_1>m,$ i.e., 
to system $(\partial)$ one equation is added at least;

2) the supplement of system $(\partial)$ to the system (2.8) is 
made by adding of the equations of the form (2.7);

3) the systems $(\partial)$ and (2.8) are
integrally equivalent on a domain $G^{\;\!\prime}\subset G,$ i.e.,
each first integral on the domain $G^{\;\!\prime}$ of system $(\partial)$
is a first integral on the domain $G^{\;\!\prime}$ of system (2.8)
and on the contrary each first integral on the domain $G^{\;\!\prime}$ of system (2.8)
is a first integral on the domain $G^{\;\!\prime}$ of system $(\partial).$

If the system (2.8) is complete, then the process is completed.

If the system (2.8) is incomplete, then 
we do the similar procedure to the system (2.8) and 
obtain one more system.

Note that on each step of this procedure
the number of equations of an incomplete system increases at least by one. 
This argument shows that 
after a finite number of steps we get either the  complete system or the  system of $n$ equations.

Note also such conditions.

First note that the system $(\partial)$ with $m=n$ is a complete system.
Indeed, this follows from that 
the set of  $n$ not linearly bound on the domain $G$ from the 
$n\!$-dimensional space $\R^n$ first order linear differential operators of  
$n$ variables is a basis of linear differential operators on the domain $G.$

Therefore we may say that any incomplete system can be reduced to the  complete system by 
a finite number of steps of the described procedure.

Secondly note that the incomplete system $(\partial)$ is reduced  
to the integrally equivalent complete system by adding 
the linear homogeneous par\-ti\-al differential equations of the forms 
\\[2ex]
\mbox{}\hfill                                        
$
\bigl[{\frak L}^{}_{j^{}_\nu}(x),
{\frak L}^{}_{l_\mu}(x)\bigr]y=0,
\qquad
\left[{\frak L}^{}_{\alpha^{}_\xi}(x),
\bigl[{\frak L}^{}_{j^{}_\nu}(x),
{\frak L}^{}_{l^{}_\mu}(x)\bigr]\right]y=0,
\hfill
$
\\[2.5ex]
\mbox{}\hfill
$
\Bigl[{\frak L}^{}_{\beta^{}_\zeta}(x),
\Bigl[{\frak L}^{}_{\alpha^{}_\xi}(x),
\bigl[{\frak L}^{}_{j^{}_\nu}(x), 
{\frak L}_{l^{}_\mu}(x)\bigr]
\Bigr]\Bigr]y=0,\,\ldots\,,
\hfill
$
\\[0.15ex]
\mbox{}\hfill (2.9)
\\[0.15ex]
\mbox{}\hfill
$
\nu=1,\ldots,m^{}_1\,,\ \
\mu=1,\ldots,m^{}_2\,,\ \
\xi=1,\ldots,m^{}_3\,, \ \ 
\zeta=1,\ldots,m^{}_4\,,\,\ldots\,, 
\hfill
$
\\[2ex]
\mbox{}\hfill
$
m^{}_s\leq m,\ \
s=1,2,\ldots, \ \ \
\{1,\ldots,m\}\ni 
j^{}_\nu\,,\ 
l^{}_\mu\,,\ 
\alpha^{}_\xi\,,\ 
\beta^{}_\zeta\,,\ 
\ldots\,.
\hfill
$
\\[1.75ex]
\indent
Using these notations, we can state the following
\vspace{0.5ex}

{\bf Theorem 2.2.}                                          
{\it
Any incomplete system $(\partial)$ with ${\frak L}^{}_j\in C^k(G),\ j=1,\ldots,m,\ k\in\N,$ 
can be reduced to the integrally equivalent on some domain $G^{\;\!\prime}\subset G$ 
complete system by adding the equations of the form {\rm (2.9)} to the system $(\partial).$
}
\vspace{0.25ex}

At this point, we may give
\vspace{0.25ex}

{\bf Definition 2.3.}                                      
{\it
We'll say that a number $\delta$ is the \textit{\textbf{defect}} of the 
incomplete system $(\partial)$ if 
this system can be reduced to the integrally equivalent on some domain $G^{\;\!\prime}\subset G$ 
complete system by adding $\delta$ equations of the form {\rm (2.9)}.
}
\vspace{0.25ex}

In this definition we mean that the corresponding complete system to the incomplete system $(\partial)$
is constructed on the base of not linearly bound on the domain  $G$ operators.  

It is obvious that the incomplete system $(\partial)$ has the defect
$\delta$ such that $0<\delta\leq n-m.$

We may assume that a complete system has the defect $\delta=0.$ 
Then any system $(\partial)$ has the defect $\delta$ such that $0\leq\delta\leq n-m.$
\vspace{0.75ex}

{\bf Example 2.1.}                                           
The normal linear homogeneous 
system of partial differential equations
\\[2ex]
\mbox{}\hfill                                              
$
{\frak L}^{}_1(x)\,y=0,\qquad 
{\frak L}^{}_2(x)\,y=0,
$
\hfill (2.10)
\\[2ex]
with the linear differential operators
\\[2ex]
\mbox{}\hfill
$
{\frak L}^{}_1(x)=
\partial^{}_{x^{}_1}+\,
x^{}_5\,\partial^{}_{x^{}_4}-
x^{}_4\,\partial^{}_{x^{}_5}$
\ for all $x\in\R^5,
\hfill
$
\\[2ex]
\mbox{}\hfill
$
{\frak L}^{}_2(x)=
\partial^{}_{x^{}_2}+\,
2x^{}_3 x^{}_5\,\partial^{}_{x^{}_3}+\,
2x^{}_4 x^{}_5\,\partial^{}_{x^{}_4}+
(1-x_3^2-x_4^2+x_5^2)\,
\partial^{}_{x^{}_5}$
\ for all $x\in\R^5
\hfill
$
\\[2.25ex]
is incomplete as the Poisson bracket
\\[2ex]
\mbox{}\hfill
$
{\frak L}^{}_{21}(x)=
\bigl[{\frak L}^{}_2(x),
{\frak L}^{}_1(x)\bigr]=
2x^{}_3x^{}_4\,
\partial^{}_{x^{}_3}+
(1-x_3^2+x_4^2-x_5^2)\,
\partial^{}_{x^{}_4}+
2x^{}_4x^{}_5\,
\partial^{}_{x^{}_5}$
\ for all $x\in\R^5
\hfill
$
\\[2.25ex]
is not the null operator.

The system (2.10) is integrally equivalent on some domain $G^{\;\!\prime}\subset\R^5$
to the linear ho\-mo\-ge\-ne\-ous system of partial differential equations
\\[2ex]
\mbox{}\hfill                                             
$
{\frak L}^{}_1(x)\,y=0,\qquad 
{\frak L}^{}_2(x)\,y=0,\qquad 
{\frak L}^{}_{21}(x)\,y=0.
$
\hfill(2.11)
\\[2ex]
\indent
The Poisson bracket
\\[2ex]
\mbox{}\hfill
$
{\frak L}^{}_{1;\;\! 21}(x)=
\bigl[{\frak L}^{}_1(x),
{\frak L}^{}_{21}(x)\bigr]=
2x^{}_3x^{}_5\,\partial^{}_{x^{}_3}+
2x^{}_4x^{}_5\,\partial^{}_{x^{}_4}+
(1-x_3^2-x_4^2+x_5^2)\,
\partial^{}_{x^{}_5}$
\ for all $x\in\R^5
\hfill
$
\\[2.25ex]
is not a linear combination of the operators 
${\frak L}^{}_1,\ {\frak L}^{}_2,\ {\frak L}^{}_{21}.$ 
Therefore the partial differential system (2.11) is incomplete.

The system (2.10) is integrally equivalent on some domain $G^{\;\!\prime}\subset\R^5$
to the linear ho\-mo\-ge\-ne\-ous system of partial differential equations
\\[2ex]
\mbox{}\hfill                                             
$
{\frak L}^{}_1(x)\,y=0,\quad 
{\frak L}^{}_2(x)\,y=0,\quad 
{\frak L}^{}_{21}(x)\,y=0,\quad 
{\frak L}^{}_{1;\;\! 21}(x)\,y=0.
$
\hfill(2.12)
\\[2ex]
\indent
The Poisson brackets
\\[2ex]
\mbox{}\hfill
$
\bigl[{\frak L}^{}_1(x),
{\frak L}^{}_{1;\;\! 21}(x)\bigr]=
{}-{\frak L}^{}_{21}(x)$
\ for all $x\in\R^5,
\hfill
$
\\[2.5ex]
\mbox{}\hfill
$
\bigl[{\frak L}^{}_{21}(x),
{\frak L}^{}_{1;\;\! 21}(x)\bigr]=
4x^{}_5\,\partial^{}_{x^{}_4}-
4x^{}_4\,\partial^{}_{x^{}_5}={}
\hfill
$
\\[2ex]
\mbox{}\hfill
$
=\ \dfrac{4x^{}_5}{1-x_3^2-x_4^2-x_5^2}\ 
{\frak L}^{}_{21}(x) \ - \ 
\dfrac{4x^{}_4}
{1-x_3^2-x_4^2-x_5^2}\ 
{\frak L}^{}_{1;\;\! 21}(x)$
\ \ for all $x\in G,
\hfill
$
\\[2.75ex]
\mbox{}\hfill
$
\bigl[{\frak L}^{}_2(x),
{\frak L}^{}_{21}(x)\bigr]=
\bigl[\partial^{}_{x^{}_2}+
{\frak L}^{}_{1;\;\! 21}(x),
{\frak L}^{}_{21}(x)\bigr]=
{}-\bigl[{\frak L}^{}_{21}(x),
{\frak L}^{}_{1;\;\! 21}(x)\bigr]$
\ for all $x\in\R^5,
\hfill
$
\\[2.5ex]
\mbox{}\hfill
$
\bigl[{\frak L}^{}_2(x),
{\frak L}^{}_{1;\;\! 21}(x)\bigr]=
\bigl[\partial^{}_{x^{}_2}+
{\frak L}^{}_{1;\;\! 21}(x),
{\frak L}^{}_{1;\;\! 21}(x)\bigr]= {\frak O}$
\ for all $x\in\R^5,
\hfill
$
\\[2.5ex]
where $G$ is any domain from the set 
$D=\bigl\{x\colon\ x_3^2+x_4^2+x_5^2\ne1\bigr\}$
of the space $\R^5.$

Thus the Poisson brackets of the operators
${\frak L}^{}_1,\ {\frak L}^{}_2,\ {\frak L}^{}_{21},\ {\frak L}^{}_{1;\;\! 21}$
are the linear combinations of these operators on the domain $\!\!G.\!\!$
Therefore the system (2.12) is complete on the domain $\!\!G\!.\!\!$

The system (2.12) is obtained by adding of two equations to the system (2.10).
The incomplete system (2.10) is integrally equivalent on some domain 
$G^{\,\prime}$ from the set $D$ to the complete system (2.12). 
Thus the incomplete system (2.10) has the defect $\delta=2.$ 
\\[2ex]
\indent
{\bf
2.3. Complete system} 
\\[1ex]
\indent
{\bf Property 2.2.}                                    
{\it
A complete linear homogeneous system of  partial differential equations 
is invariant under a holomorphism.
}
\vspace{0.25ex}

{\sl Proof}. Let the map
\\[1.5ex]
\mbox{}\hfill                                         
$
x\colon\xi\to\varphi(\xi)$
\ for all $\xi\in\Omega\subset\R^n
$
\hfill (2.13)
\\[1.75ex]
be a holomorphism between the domain $\Omega$ and the domain $G$ of the space $\R^n.$
\vspace{0.25ex}

The expression ${\frak L}^{}_j(x)\,y(x)$ is invariant under the holomorphism (2.13):
\\[2.25ex]
\mbox{}\hfill                                          
$
{\frak L}^{}_j(x)\,y(x)_{\displaystyle 
\mbox{}|_{x=\varphi(\xi)}}=\,
\widetilde{\frak L}^{}_j(\xi)\,z(\xi)$
\ for all $\xi\in\Omega,$
\ for all $x\in G,\ \
j=1,\ldots,m,
$
\hfill (2.14)
\\[2.5ex]
where $z(\xi)=y(\varphi(\xi))$ for all $\xi\in\Omega.$
Using the transformation (2.13), we have  
the system $(\partial)$ is reduced to the system
\\[2ex]
\mbox{}\hfill                                           
$
\widetilde{\frak L}^{}_j(\xi)\,z=0,
\quad 
j=1,\ldots,m,
$
\hfill (2.15)
\\[2ex]
with the not linearly bound on the domain $\Omega$ 
linear differential operators of first order $\widetilde{\frak L}^{}_j\,,\ j=1,\ldots,m$ 
(because a holomorphism is bijective). 

Let us prove that the system (2.15) is complete on the domain $\Omega$ 
under the condition the system $(\partial)$ is complete on the domain $G.$

Taking into account the identities (2.14), we have 
\\[2ex]
\mbox{}\hfill
$
{\frak L}^{}_j(x)\,
{\frak L}^{}_l(x)\,y(x)_{\displaystyle \mbox{}|_{x=\varphi(\xi)}}=\,
{\frak L}^{}_j(x)\, v^{}_l(x)_{\displaystyle \mbox{}|_{x=\varphi(\xi)}}=\,
\widetilde{\frak L}^{}_j(\xi)\, v^{}_l(\varphi(\xi))=\,
\widetilde{\frak L}^{}_j(\xi)\,
\widetilde{\frak L}^{}_l(\xi)\,z(\xi) 
\hfill
$
\\[3.25ex]
\mbox{}\hfill
for all $\xi\in\Omega,$ \ \
for all $x\in G, 
\quad 
j=1,\ldots,m,\ \
l=1,\ldots,m,
\hfill
$
\\[2.25ex]
where $v^{}_l(x)={\frak L}^{}_l\, y(x)$ for all $x\in G,\ l=1,\ldots,m.$
Hence,
\\[2.5ex]
\mbox{}\hfill                                     
$
\Bigl[{\frak L}^{}_j(x),
{\frak L}^{}_l(x)\Bigr]y(x)_{\displaystyle 
\mbox{}|_{x=\varphi(\xi)}}=
\Bigl[\widetilde{\frak L}^{}_j(\xi),
\widetilde{\frak L}^{}_l(\xi)
\Bigr] z(\xi)$
\ \ 
for all $\xi\in\Omega,$\ \
for all $x\in G, 
\hfill
$ 
\\
\mbox{}\hfill  (2.16)
\\
\mbox{}\hfill
$
j=1,\ldots,m,\ \ \ 
l=1,\ldots,m.
\hfill
$
\\[2ex]
\indent
Since the system $(\partial)$ is complete, we see that 
the identities (0.12) are fulfilled. Therefore,
\\[2.25ex]
\mbox{}\hfill
$
\left[\widetilde{\frak L}^{}_j(\xi),
\widetilde{\frak L}^{}_l(\xi)
\right]z(\xi)=
\left[{\frak L}^{}_j(x),
{\frak L}^{}_l(x)\right]
y(x)_{\displaystyle
\mbox{}|_{x=\varphi(\xi)}}=
\sum\limits_{\nu=1}^m\, 
A^{}_{jl\nu}(x)\,
{\frak L}^{}_\nu(x)\,
y(x)_{\displaystyle
\mbox{}|_{x=\varphi(\xi)}}=
\hfill
$
\\[2.5ex]
\mbox{}\hfill
$
\displaystyle
=\ \sum\limits_{\nu=1}^m\, 
\widetilde{A}^{}_{jl\nu}(\xi)\,
\widetilde{\frak L}^{}_\nu(\xi)\,z(\xi)$
\ for all $\xi\in\Omega,$
\ for all $x\in G, 
\quad 
j=1,\ldots,m,\ \
l=1,\ldots,m.\ \ \k
\hfill
$
\\[2.25ex]
\indent
{\bf Property 2.3.}                                    
{\it
A jacobian linear homogeneous system of  partial differential equations 
is invariant under a holomorphism.
}
\vspace{0.25ex}

The {\sl proof} is analogous to the proof of Property 2.2 if we take into account that 
the identities (0.14) for the jacobian system $(\partial)$ are valid.\ \k
\vspace{0.5ex}

{\bf Property 2.4.}                                  
\vspace{0.35ex}
{\it
The complete system $(\partial)$ in a neighbourhood of any point $x\in G$  that sa\-tis\-fies 
\vspace{0.35ex}
${\rm det}\bigl\|\psi^{}_{jl}(x)\bigr\|\ne 0$ can be reduced to an integrally equivalent 
on some do\-main $G^{\;\!\prime}\subset G$ 
complete system by the nonsingular on the domain $G$ linear transformation 
of the operators}
\\[2ex]
\mbox{}\hfill                                        
$
\displaystyle
{\frak L}^{}_j(x)=
\sum\limits_{l=1}^m\, \psi^{}_{jl}(x)\,{\frak N}^{}_l(x)$
\ for all $x\in G,
\quad 
j=1,\ldots,m,
$
\hfill(2.17)
\\[2ex]
{\it
where the linear differential operators of first order
\vspace{0.35ex}
${\frak N}^{}_l\,,\ l=1,\ldots,m,$ and the scalar functions 
$\psi^{}_{jl}\colon G\to\R,\ j=1,\ldots,m,\ l=1,\ldots,m,$
are holomorphic on the domain $G.$ 
}
\vspace{0.25ex}

{\sl Proof}. 
The linear transformation (2.17) is nonsingular.
Hence the linear differential operators ${\frak N}^{}_l,\, l=1,\ldots,m,$ 
can be presented as the linear combinations of operators (0.1):
\\[1.75ex]
\mbox{}\hfill                                         
$
\displaystyle
{\frak N}^{}_l(x)=
\sum\limits_{j=1}^m\,
\theta^{}_{lj}(x)\, {\frak L}^{}_j(x)$
\ for all $x\in\widetilde{G}\subset G,
\quad 
l=1,\ldots,m,
$
\hfill(2.18)
\\[1.75ex]
and from the system $(\partial),$ we obtain the system 
\\[1.75ex]
\mbox{}\hfill                                          
$
\displaystyle
\sum\limits_{l=1}^m\,
\psi^{}_{jl}(x)\, {\frak N}^{}_l(x)\,y=0,
\quad 
j=1,\ldots,m.
$
\hfill(2.19)
\\[1.75ex]
\indent
The domain $\widetilde{G}$ in the expression (2.18) such that 
${\rm det}\bigl\|\psi^{}_{jl}(x)\bigr\|\ne 0$ for all $x\in\widetilde{G}.$
\vspace{0.5ex}

The system (2.19) disintegrates on the system of equations 
\\[2ex]
\mbox{}\hfill                                          
$
{\frak N}^{}_l(x)\,y=0,
\quad 
l=1,\ldots,m, 
$
\hfill(2.20)
\\[2ex]
where the linear differential operators of first order
${\frak N}^{}_l\,,\ l=1,\ldots,m,$
are  not linearly bound on the domain $\widetilde{G}.$

From the notion (2.19) and that the matrix $\bigl\|\psi^{}_{jl}(x)\bigr\|$ of order $m$
is nonsingular on the domain $\widetilde{G}\subset G$ 
it follows that
the system (2.20) is integrally equivalent on some domain $G^\prime\subset\widetilde{G}$
to the system $(\partial).$

Let us show that the system (2.20) is complete.

Using the operator identities
\\[2ex]
\mbox{}\hfill
$
\bigl[f(x)\,{\frak L}^{}_j(x), g(x)\,{\frak L}^{}_l(x)\bigr]=
f(x)\,g(x)\, \bigl[{\frak L}^{}_j(x), {\frak L}^{}_l(x)\bigr] +
f(x)\, {\frak L}^{}_j g(x)\, {\frak L}^{}_l(x) -
g(x)\, {\frak L}^{}_l f(x)\,{\frak L}^{}_j(x)
\hfill
$
\\[2.5ex]
\mbox{}\hfill
for all $x\in G,
\quad
f\in C^1(G),\ \
g\in C^1(G),\ \ \ 
j=1,\ldots,m,\ \
l=1,\ldots,m, 
\hfill
$
\\[1.25ex]
and
\\[1ex]
\mbox{}\hfill
$
\bigl[{\frak L}^{}_j(x)+
{\frak L}^{}_l(x),
{\frak L}^{}_\zeta(x)\bigr]=
\bigl[{\frak L}^{}_j(x),
{\frak L}^{}_\zeta(x)\bigr]+
\bigl[{\frak L}^{}_l(x),
{\frak L}^{}_\zeta(x)\bigr]$
\ for all $x\in G,
\hfill
$
\\[2.25ex]
\mbox{}\hfill
$
j=1,\ldots,m,\ \
l=1,\ldots,m,\ \
\zeta=1,\ldots,m, 
\hfill
$
\\[2ex]
and the representations (2.18), we get
\\[1.75ex]
\mbox{}\hfill
$
\displaystyle
\bigl[{\frak N}^{}_\mu(x),
{\frak N}^{}_\nu(x)\bigr] \ = \
\sum\limits_{j=1}^m 
\sum\limits_{l=1}^m\,
A^{}_{jl\mu\nu}(x)\bigl[{\frak L}^{}_j(x),
{\frak L}^{}_l(x)\bigr] \ + \
\sum\limits_{s=1}^m\, B^{}_{s\mu\nu}(x)\,{\frak L}^{}_s(x)
\hfill
$
\\[2.25ex]
\mbox{}\hfill
for all $x\in\widetilde{G}, 
\quad
\mu=1,\ldots,m,\ \ 
\nu=1,\ldots,m.
\hfill
$
\\[2ex]
\indent
From here using the decompositions (0.12) and the transformation (2.17), 
\vspace{0.35ex}
we obtain the Poisson brackets
$\bigl[{\frak N}^{}_\mu(x), {\frak N}^{}_\nu(x)\bigr]$ for all $x\in\widetilde{G},\  
\mu=1,\ldots,m,\ \nu=1,\ldots,m,$
\vspace{0.35ex}
are linear combinations on the domain $\widetilde{G}$ of the operators
${\frak N}^{}_j\,,\ j=1,\ldots,m.$
This means that the system (2.20) is complete. \k
\vspace{0.35ex}

From Property 2.4, we have the following
\vspace{0.35ex}

{\bf Property 2.5.}                                
{\it
If the system $(\partial)$ is complete, then the system}
\\[2ex]
\mbox{}\hfill                                      
$
{\frak D}^{}_j(x)\,y=0,
\quad 
j=1,\ldots,m,
$
\hfill(2.21)
\\[2ex]
{\it
where the linear differential operators of first order}
\\[2ex]
\mbox{}\hfill                                       
$
\displaystyle
{\frak D}^{}_j(x)=
\sum\limits_{l=1}^m\, 
v^{}_{jl}(x)\,{\frak L}^{}_l(x)$
\ for all $x\in G,
\quad 
j=1,\ldots,m,
$
\hfill(2.22)
\\[2ex]
{\it
the functional square matrix $v(x)=\bigl\|v^{}_{jl}(x)\bigr\|$ of order $m$
\vspace{0.35ex}
is nonsingular on the domain $G,$
is also complete and integrally equivalent to the system $(\partial)$
on a neighbourhood of any point $x\in G$ that satisfies ${\rm det}\,v(x)\ne 0.$
}
\vspace{0.5ex}

{\bf Theorem 2.3.}                                     
{\it
The complete system $(\partial)$ can be reduced to  
an integrally equivalent on some do\-main $G^{\;\!\prime}\subset G$ 
complete normal system by the nonsingular on the domain $G$ linear tran\-s\-for\-ma\-tion 
of operators {\rm(0.1)} 
{\rm (}under this transformation we may have an restriction of the domain $G).$
}
\vspace{0.5ex}

{\sl Proof}. 
Let the system $(\partial)$ be complete.
\vspace{0.5ex}
Then the square matrix $\widehat{u}(x)=\bigl\|u^{}_{ji}(x)\bigr\|$ for all $x\in G$ of order $\!m\!$
\vspace{0.5ex}
(this matrix is obtain from the $m\times n$ matrix  $u(x)=\bigl\|u^{}_{ji}(x)\bigr\|$ for all $x\!\in\! G$ 
by taking out the first  $m$ columns)  is nonsingular on the domain $G$
(since ${\rm rank}\,u(x)=m$ almost everywhere on $G,$ we see that 
it always can be received by renumbering variables).
Therefore there exists a linear nonsingular transformation of operators (0.1)
such that the system $(\partial)$ can be reduced to the normal system $(N\partial).$

By Property 2.5, this normal system is complete.\ \k
\vspace{0.5ex}

From Theorem 2.2, we get the following
\vspace{0.5ex}

{\bf Theorem 2.4.}\!                                       
{\it
The complete system $\!(\partial)\!$ can be reduced to  
an integrally equivalent on some do\-main $\!G^{\;\!\prime}\!\subset\! G\!$
jacobian system by the nonsingular on the domain $\!G\!$ linear transformation 
of operators\! {\rm(0.1)}\!
{\rm (}\!under this transformation we may have an restriction 
\vspace{0.5ex}
of the domain $\!G).\!\!$
}

Using Property 2.4 (or Property 2.5) and the process of building of complete normal system 
(see the proof of Theorem 2.3), we get the following theorem  
about integral equivalence of the complete system $(\partial)$ and a complete normal system. 
\vspace{0.5ex}

{\bf Theorem 2.5.}                                       
{\it
\vspace{0.5ex}
Suppose the complete system $(\partial)$ has 
the nonsingular in the domain $G$ square matrix $\widehat{u}$ of order $m$
{\rm(}this matrix is obtain from the $m\times n$ matrix  $u(x)=\bigl\|u^{}_{ji}(x)\bigr\|$ for all $x\!\in\! G$ 
by taking out the first  $m$ columns{\rm)}.
Then the complete system $(\partial)$ can be reduced to the complete normal system $(N\partial)$ 
and these systems in a neighbourhood of any point $x\in G$ that sa\-tis\-fies 
${\rm det}\,{\widehat{u}}(x)\ne 0$ are integrally equivalent.
}
\vspace{0.5ex}

Theorems 2.2 and 2.5 are the substantiation of the following notion  [1, p. 48].
\vspace{0.5ex}

{\bf Definition 2.4.}                                   
{\it
We'll say that a domain $H\subset G$ is a \textit{\textbf{normalization domain}} of system $(\partial)$ if
the system $(\partial)$ in a neighbourhood of any point of the domain $H$ can be reduced 
to an integrally equivalent complete normal system.
}
\vspace{0.5ex}

Notice that a normalization domain of an incomplete system is 
a normalization domain of an integrally equivalent complete system (see Theorem 2.2).

In general case a normalization domain is ambiguous.
This normalization domain depends on zeroes of determinants of square matrices of order $m$ 
(these matrices are obtain from the matrix  $u(x)=\bigl\|u^{}_{ji}(x)\bigr\|_{m\times n}$ for all $x\in G$ 
by taking out $m$ columns).
\vspace{0.5ex}

{\bf Example 2.2.}                                        
Consider the linear ho\-mo\-ge\-ne\-ous system of partial differential equations
\\[2ex]
\mbox{}\hfill                                            
$
{\frak L}^{}_1(x)\,y=0,
\qquad
{\frak L}^{}_2(x)\,y=0,
$
\hfill(2.23)
\\[1.75ex]
where the linear differential operators of the first order
\\[2ex]
\mbox{}\hfill                                          
$
{\frak L}^{}_1(x)=
x^{}_1\,\partial^{}_{x^{}_1}+
x^{}_2\,\partial^{}_{x^{}_2}+
x^{}_3\,\partial^{}_{x^{}_3}+
x^{}_4\,\partial^{}_{x^{}_4}+
x^{}_5\,\partial^{}_{x^{}_5}$
\ for all $x\in\R^5, 
\hfill
$
\\[2.5ex]
\mbox{}\hfill
$
{\frak L}^{}_2(x)=
x^{}_1\,\partial^{}_{x^{}_1}+
x^{}_2\,\partial^{}_{x^{}_2}+
x^{}_3\,\partial^{}_{x^{}_3}+
x^{2}_4\,\partial^{}_{x^{}_4}+
x^{2}_5\,\partial^{}_{x^{}_5}$
for all $x\in\R^5.
\hfill
$
\\[2ex]
\indent
The Poisson bracket
\\[2ex]
\mbox{}\hfill
$
\bigl[{\frak L}^{}_1(x),
{\frak L}^{}_2(x)\bigr]=
x_4^2\,\partial^{}_{x^{}_4}+
x_5^2\,\partial^{}_{x^{}_5}=
{\frak L}^{}_{12}(x)$
\ for all $x\in\R^5
\hfill
$
\\[2.25ex]
is not a linear combination of the operators ${\frak L}^{}_1$ and ${\frak L}^{}_2.$
\vspace{0.35ex}
Hence the system (2.23) is incomplete.

Using the operator ${\frak L}^{}_{12},$ we get 
the system (2.23) is reduced to the system
\\[2ex]
\mbox{}\hfill                                          
$
{\frak L}^{}_1(x)\,y=0,
\quad \
{\frak L}^{}_2(x)\,y=0,
\quad \
{\frak L}^{}_{12}(x)\,y=0. 
$
\hfill(2.24)
\\[2ex]
\indent
Since the Poisson brackets
\\[2ex]
\mbox{}\hfill
$
\bigl[{\frak L}^{}_1(x),
{\frak L}^{}_{12}(x)\bigr]=
{\frak L}^{}_{12}(x)$
\ for all $x\in\R^5,
\qquad
\bigl[{\frak L}^{}_2(x),
{\frak L}^{}_{12}(x)\bigr]=
{\frak O}$
\ for all $x\in\R^5,
\hfill
$
\\[2ex]
we see that the system (2.24) is complete.

Therefore the incomplete system (2.23) has the defect $\delta=1.$

From the second equation of system (2.24) by virtue of the third equation of this system
it follows that
\\[0.75ex]
\mbox{}\hfill                                          
$
x^{}_1\,\partial^{}_{x^{}_1}y+
x^{}_2\,\partial^{}_{x^{}_2}y+
x^{}_3\,\partial^{}_{x^{}_3}y=0.
$
\hfill(2.25)
\\[2.25ex]
\indent
Then from the first equation of system (2.23), we have
\\[2ex]
\mbox{}\hfill
$
x^{}_4\partial^{}_{x^{}_4}y+
x^{}_5\partial^{}_{x^{}_5}y=0.
\hfill
$
\\[2.25ex]
\indent
From this equation and the third equation of system (2.24), we obtain the equalities
\\[2ex]
\mbox{}\hfill
$
\partial^{}_{x^{}_4}y=0,
\quad \  
\partial^{}_{x^{}_5}y=0.
\hfill
$
\\[2.25ex]
\indent
By solving the equation (2.25) for $\partial^{}_{x^{}_1}y$ we reduce the system (2.24) 
 to the normal system
\\[2ex]
\mbox{}\hfill                                        
$
\partial^{}_{x^{}_1}y=
{}-x^{}_2 x_1^{-1}\,\partial^{}_{x^{}_2}y-
x^{}_3 x_1^{-1}\,\partial^{}_{x^{}_3}y,
\qquad
\partial^{}_{x^{}_4}y=0,
\qquad
\partial^{}_{x^{}_5}y=0.
$
\hfill(2.26)
\\[2.5ex]
\indent
The complete normal system (2.26) is 
integrally equivalent to the complete system (2.24) and to the 
incomplete system (2.23) on a normalization domain  
$H^{}_1\subset\{x\colon x^{}_1\ne 0\}\subset\R^5.$ 
\vspace{0.35ex}

It is readily seen that the systems (2.23) and (2.24) have 
else two normal forms, which can be obtain by solving the equation (2.25) 
for $\partial^{}_{x^{}_2}y$ and $\partial^{}_{x^{}_3}y\colon$
\\[2ex]
\mbox{}\hfill                                        
$
\partial^{}_{x^{}_2}y=
{}-x^{}_1 x_2^{-1}\,\partial^{}_{x^{}_1}y-
x^{}_3 x_2^{-1}\,\partial^{}_{x^{}_3}y,
\qquad
\partial^{}_{x^{}_4}y=0,
\qquad
\partial^{}_{x^{}_5}y=0
$
\hfill(2.27)
\\[1ex]
and
\\[1ex]
\mbox{}\hfill                                        
$
\partial^{}_{x^{}_3}y=
{}-x^{}_1 x_3^{-1}\,\partial^{}_{x^{}_1}y-
x^{}_2 x_3^{-1}\,\partial^{}_{x^{}_2}y,\quad
\partial^{}_{x^{}_4}y=0,\quad
\partial^{}_{x^{}_5}y=0.
$
\hfill(2.28)
\\[2.25ex]
\indent
The complete normal system (2.27) is 
integrally equivalent to the complete system (2.24) and to the 
incomplete system (2.23) on a normalization domain  
$H^{}_2\subset\{x\colon x^{}_2\ne 0\}\subset\R^5.$
\vspace{0.35ex}

The complete normal system (2.28) is 
integrally equivalent to the complete system (2.24) and to the 
incomplete system (2.23) on a normalization domain  
\vspace{2.75ex}
$H^{}_3\subset\{x\colon x^{}_3\ne 0\}\subset\R^5.$

{\bf
2.4. Dimension of integral basis}          
\\[1ex]
\indent
The system of total differential equations
\\[2ex]
\mbox{}\hfill                                          
$
\displaystyle
dx^{}_s={}-\sum\limits_{j=1}^m\,
u^{}_{js}(x)\,dx^{}_j\,,
\quad 
s=m+1,\ldots,n,
$
\hfill(2.29)
\\[2ex]
is associated to the normal linear homogeneous system of partial differential equations $(N\partial).$ 
For the system (2.29), as well as for the system {\rm (TD)},
the linear differential operators ${\frak X}^{}_j\;\!,
\linebreak 
j=1,\ldots,m,$ have the form
\\[2ex]
\mbox{}\hfill                                           
$
{\frak X}^{}_j(x)=
{\frak L}^{}_j(x)=
\partial^{}_{x^{}_j}-{\frak M}^{}_j(x)$
\ for all $x\in G,
\quad 
j=1,\ldots,m,
$
\hfill(2.30)
\\[2.15ex]
where the operators ${\frak M}^{}_j\,,\ j=1,\ldots,m,$ are given by (0.12).
\vspace{0.75ex}

{\bf Definition 2.5.}
{\it
We'll say that a total differential system 
and a linear homogeneous system of partial differential equations
are \textit{\textbf{integrally equivalent}} on some domain if  on this domain
each first integral of the first system is a first integral of the second system 
and on the contrary each first integral of the second system 
is a first integral of the first system.
}
\vspace{0.35ex}

If a total differential system and a linear ho\-mo\-ge\-ne\-ous system of partial differential equations 
are integrally equivalent on the domain $\Pi^\prime,$ then these systems  
have the same integral basis on this domain (see Definitions 1.2, 2.2, and 2.5).

From the connections (2.30) it follows that the identities (1.2) for the function $F\colon G^{\;\!\prime}\to\R$
coincide with the identities (2.1) for the function $F\colon G^{\;\!\prime}\to\R.$
By Definitions 1.1 and 2.1, we get the next theorem about the integral equivalence of a 
linear ho\-mo\-ge\-ne\-ous system of partial differential equations with a system of total differential equations. 

\newpage

{\bf Theorem 2.6.}                                        
{\it
The scalar function $F\colon G^\prime\to\R$ is a first integral on the domain
$G^\prime\subset G$ of the normal linear homogeneous system of partial differential equations
if and only if this function is a first integral on the domain $G$ of the 
total differential system {\rm(2.29).}
}
\vspace{0.35ex}

Using the definitions of the complete and jacobian 
linear homogeneous systems of partial differential equations $(\partial)$ 
and the definition of  the completely solvable total differential system (TD),
we get the relation between these notions.
\vspace{0.35ex}

{\bf Theorem 2.7.}\!                                        
{\it
The normal linear homogeneous system of partial differential equations $\!\!(\!N\!\partial)\!\!$ 
is complete {\rm(}jacobian{\rm)} if and only if  
\vspace{0.35ex}
the total differential system {\rm(2.29)} is completely solvable.
}

By Theorem 2.6, using the definitions of an integral basis for a total differential system (Definition 1.2) 
and for a linear homogeneous systems of partial differential equations (De\-fi\-ni\-ti\-on 2.2), 
we have the relation between integral bases 
\vspace{0.35ex}
of the systems  $\!(N\partial)\!$ and (2.29).

{\bf Theorem 2.8.}                                        
{\it
The set of functions {\rm(2.2)} is a basis of first integrals 
on a domain $G^\prime\subset G$ of the 
normal linear homogeneous system of partial differential equations $(N\partial)$ 
if and only if  this set of functions {\rm(2.2)} is a basis of first integrals 
on the domain $G^\prime$ of the system of total differential equations {\rm(2.29)}.
}
\vspace{0.35ex}

From Theorem 1.3 it follows that the dimension of an integral basis for 
the completely solvable total differential system {\rm(2.29)} is equal $n-m.$
Then, using Theorems 2.7 and 2.8, we obtain the dimension of an integral basis for the 
complete normal linear homogeneous system of partial differential equations $(N\partial).$ 
\vspace{0.35ex}

{\bf Theorem 2.9.}                                         
{\it
The complete {\rm(}jacobian{\rm)}
normal linear homogeneous system of partial dif\-fe\-ren\-ti\-al equations $(N\partial)$ 
with ${\frak M}^{}_j\in C^\infty(G)$ on a neighbourhood of any point from the domain $G$ 
has an integral basis of dimension $n-m.$
}
\vspace{0.75ex}

{\bf Example 2.3.}                                           
The complete {\rm(}jacobian{\rm)} normal linear homogeneous system of partial 
dif\-fe\-ren\-ti\-al equations 
\\[2ex]
\mbox{}\hfill                                              
$
\partial^{}_{t^{}_1}y=
{}-\partial^{}_{x^{}_1}y-
\partial^{}_{x^{}_1}
g(x^{}_1\,,x^{}_2)\,
\partial^{}_{x^{}_3}y,
\qquad
\partial^{}_{t^{}_2}y=
{}-\partial^{}_{x^{}_2}y-
\partial^{}_{x^{}_2}
g(x^{}_1\,,x^{}_2)\,
\partial^{}_{x^{}_3}y,
$
\hfill (2.31)
\\[2.25ex]
where the scalar function $g\in C^{\;\!\infty}(D),\ D\subset\R^2,$ is 
associated to the completely solvable system of total differential equations  (1.15). 
\vspace{0.25ex}

The system (2.31) has the form $(N\partial)$ with $m=2,\,n=5.$
\vspace{0.25ex}

In Example 1.3 an integral basis on the domain $\Pi^\prime=\R^2\times D\times\R$ for system (1.15)
was constructed. This integral basis is three functionally independent on the domain $\Pi^\prime$ 
first integrals (1.16) of system (1.15).

By Theorem 2.6, the functions (1.16) are first integrals on the domain $\Pi^\prime$ 
of the partial differential system (2.31).

From Theorems 2.8 and 2.9 it follows that 
the complete normal system (2.31) has an integral basis on the domain $\Pi^\prime$ of dimension $n-m=5-2=3.$
This integral basis is the functionally independent on the domain $\Pi^\prime$ 
first integrals (1.16).
\vspace{0.35ex}

Using the definition of the normalization domain (Definition 2.4), Theorems 2.5 
(to the effect that the complete system $(\partial)$ can be reduced to an 
integrally equivalent complete normal system) and 
2.9 (about the dimension of an integral basis of the complete normal system $(N\partial)),$ 
we obtain the dimension of an integral basis for the complete system $(\partial).$
\vspace{0.35ex}

{\bf Theorem 2.10.}                                       
{\it
The complete linear homogeneous system of partial differential equations $(\partial)$
with ${\frak L}^{}_j\in C^\infty(G),\ j=1,\ldots,m,$ on 
a neighbourhood of any point of a normalization domain 
has a basis of first integrals of dimension $n-m.$
}
\vspace{0.35ex}

{\bf Example 2.4.}                                        
We consider the linear homogeneous system of partial differential equations
\\[2ex]
\mbox{}\hfill                                            
$
{\frak L}^{}_1(x)\,y=0,
\qquad
{\frak L}^{}_2(x)\,y=0
$
\hfill(2.32)
\\[2ex]
with the linear differential operators of first order
\\[2ex]
\mbox{}\hfill                                          
$
{\frak L}^{}_1(x)=
\partial^{}_{x^{}_1}+
\partial^{}_{x^{}_2}+
\partial^{}_{x^{}_3}$
\ for all $x\in\R^3, 
\quad \
{\frak L}^{}_2(x)=
x^{}_1\,\partial^{}_{x^{}_1}+
x^{}_2\,\partial^{}_{x^{}_2}+
x^{}_3\,\partial^{}_{x^{}_3}$
\ for all $x\in\R^3.
\hfill
$
\\[2ex]
\indent
Since the Poisson bracket
\\[1.5ex]
\mbox{}\hfill
$
\bigl[{\frak L}^{}_1(x),
{\frak L}^{}_2(x)\bigr]=
{\frak L}^{}_1(x)$
\ for all $x\in\R^3,
\hfill
$
\\[2ex]
we see that the system (2.32) is complete.

By Theorem 2.10, the system (2.32) has an integral basis of dimension $n-m=3-2=1.$

Using the definition of a first integral (Definition 2.1), we obtain 
a first integral of system (2.32) and consequently an integral basis of system  (2.32).

Thus the function
\\[1.75ex]
\mbox{}\hfill
$
F^{}_{12}\colon x\to
(x^{}_2-x^{}_3)(x^{}_1-x^{}_2)^{-1}$
\ for all $x\in D^{}_{12}\,,
\quad 
D^{}_{12}=\{x\colon x^{}_2\ne x^{}_1\},
\hfill
$
\\[2ex]
is a basis of first integrals for the system (2.32) on any domain $G^{}_{12}$ 
from the set $D^{}_{12}.$

Similarly the function
\\[2ex]
\mbox{}\hfill
$
F^{}_{13}\colon x\to
(x^{}_3-x^{}_2)
(x^{}_1-x^{}_3)^{-1}$
\ for all $x\in D^{}_{13}\,,
\quad 
D^{}_{13}=\{x\colon x^{}_3\ne x^{}_1\},
\hfill
$
\\[2.25ex]
is a basis of first integrals for the system (2.32) on any domain $G^{}_{13}$ 
from the set $D^{}_{13}.$

In the same way the function
\\[1.75ex]
\mbox{}\hfill
$
F^{}_{23}\colon x\to
(x^{}_3-x^{}_1)
(x^{}_2-x^{}_3)^{-1}$
\ for all $x\in D^{}_{23}\,,
\quad 
D^{}_{23}=\{x\colon x^{}_3\ne x^{}_2\},
\hfill
$
\\[2.25ex]
is a basis of first integrals for the system (2.32) on any domain $G^{}_{23}$ 
from the set $D^{}_{23}.$

Each of the integral bases $F^{}_{12}\,,\ F^{}_{13},$ and $F^{}_{23}$ 
of the linear homogeneous partial differential system (2.32) defines the family of integral surfaces of this system
which consists of the planes $C^{}_1 x^{}_1+C^{}_2 x^{}_2-(C^{}_1+C^{}_2) x^{}_3=0,$
where $C^{}_1$ and $C^{}_2$ are arbitrary real constants.
\vspace{0.5ex}

Using the notion of defect for a system (Definition 2.3), 
the theorem of the reduction of an incomplete system to an 
integrally equivalent complete system (Theorem 2.2), and the theorem of 
dimension for an integral basis of a complete system (Theorem 2.10), we get 
the theorem of dimension for an integral basis of an incomplete system. 
\vspace{0.35ex}

{\bf Theorem 2.11.}                                      
{\it
Suppose the incomplete linear homogeneous system of partial differential equations $(\partial)$
with ${\frak L}^{}_j\!\in\! C^\infty(G),\, j\!=\!1,\ldots,m,$ has the defect $\!\delta.\!$ Then 
this system on a neig\-h\-bo\-ur\-ho\-od of any point of a normalization domain 
\vspace{0.35ex}
has an integral basis of dimension $n\!-\!m\!-\!\delta.\!$
}

Using Property 2.1 and the procedure of the reduction of an incomplete system to a complete system 
(see Subsection 2.2), we obtain
\vspace{0.35ex}

{\bf Corollary 2.1.}                                      
{\it
Suppose the incomplete linear homogeneous system of partial differential equations $(\partial)$
has $n-1$ equations with $n$ unknowns.
Then first integrals of this system are only arbitrary constants.
}
\vspace{0.35ex}

Recall that the complete system $(\partial)$ has the defect $\delta=0.$ 
Using this notation, Theorems 2.10 and 2.11, we can state the following
\vspace{0.35ex}

{\bf Theorem 2.12.}                                        
{\it
Suppose the linear homogeneous system of partial differential equations $(\partial)$
with ${\frak L}^{}_j\in C^\infty(G),\ j=1,\ldots,m,$ has the defect $\delta,\ 0\leq\delta\leq n-m.$ Then 
this system on a neig\-h\-bo\-ur\-ho\-od of any point of a normalization domain 
has a basis of first integrals of dimension $n-m-\delta.$
}
\vspace{0.35ex}

From Theorem 2.12, we obtain {\sl the completeness criterion for 
linear homogeneous system of partial differential equations}.
\vspace{0.35ex}

{\bf Theorem 2.13.}                                     
{\it
The linear homogeneous system of partial differential equations $(\partial)$
with ${\frak L}^{}_j\in C^\infty(G),\ j=1,\ldots,m,$ is complete if and only if this system 
on a neig\-h\-bo\-ur\-ho\-od of any point of a normalization domain 
has a basis of first integrals of dimension $n-m.$
}
\vspace{0.5ex}

The following agreement is needed for the sequel.
\vspace{0.5ex}

{\bf Agreement 2.1.}                                       
{\it
By ${\frak L}^{{}^{\scriptstyle *}}_{\nu}(x)$ for all $x\in G,\ \nu=1,\ldots,p,$ denote
linear differential ope\-ra\-tors of first order, which constructed on the base of the operators {\rm(0.1)},
such that the ope\-ra\-tors ${\frak L}^{}_1\,,\ldots, {\frak L}^{}_m\,,
{\frak L}^{{}^{\scriptstyle *}}_1\,,\ldots, {\frak L}^{{}^{\scriptstyle *}}_p$
are not linearly bound on the domain $G$ and the equations 
${\frak L}^{{}^{\scriptstyle *}}_{\nu}(x)\,y=0,\ \nu=1,\ldots,p,$ have the forms {\rm(2.9)}.
}
\vspace{0.5ex}

Along with the system $(\partial)$ we'll consider 
the linear homogeneous system of partial differential equations
\\[1.75ex]
\mbox{}\hfill                                              
$
{\frak L}^{}_j(x)\,y=0,
\quad 
j=1,\ldots,m,
\qquad
{\frak L}^{{}^{\scriptstyle *}}_{\nu}(x)\,y=0,
\quad 
\nu=1,\ldots,p.
$
\hfill(2.33)
\\[2ex]
The system (2.33) is constructed on the base of  system $(\partial)$ according to Agreement 2.1.

By Lemma 2.2 and in accordance with Agreement 2.1, we obviously have
\vspace{0.25ex}

{\bf Theorem 2.14.}                                        
{\it
The system $(\partial)$ and the system {\rm(2.33)}, 
which constructed on the base of  system $(\partial)$ according to Agreement {\rm 2.1},
are integrally equivalent on some domain $G^\prime\subset G.$
}
\vspace{0.25ex}

From Theorem 2.14, we get the following assertion for an integral basis.
\vspace{0.25ex}

{\bf Corollary 2.2.}                                         
{\it
A set of scalar functions is an integral basis of dimension $r$  on a domain $G^\prime\subset G$ 
for the partial differential system $(\partial)$ if and only if this set is an integral basis of dimension $r$  on the domain 
$G^\prime$ for the partial differential system {\rm(2.33)}, 
which constructed on the base of  system $(\partial)$ according to Agreement {\rm 2.1}.
}
\vspace{0.25ex}

The system (2.33) under the condition $p=0$ is the system $(\partial).$

If $p=\delta,$ where $\delta$ is a defect of system $(\partial),$
then the system (2.33) is complete.

Using the pocedure of the reduction of an incomplete system to a complete system (see Subsection 2.2), 
Agreement {\rm 2.1}, Theorem 2.10, and Corollary 2.2, we clearly have
\vspace{0.35ex}

{\bf Theorem 2.15.}                                        
{\it
The scalar functions
$F^{}_\tau\colon G^\prime\to\R,\ \tau=1,\ldots,m-n-\delta,$
are a basis of first integrals on the domain $G^\prime\subset G$ for system $(\partial)$ 
with ${\frak L}^{}_j\in C^\infty(G),\ j=1,\ldots,m,$ and with the defect $\delta,\ 0\leq\delta\leq n-m,$ 
if and only if this functions are a basis of first integrals on the domain $G^\prime$ 
for the complete system {\rm(2.33)}, where $p=\delta.$
}
\vspace{0.5ex}

{\bf Example 2.5.}                                       
In accordance with the definition of a first integral (Definition 2.1), 
we obtain the scalar function 
\\[2ex]
\mbox{}\hfill                                              
$
F\colon x\to\,
x^{}_3(1+x_3^2+x_4^2+x_5^2)^{{}-1}
$ \
for all 
$
x\in\R^5
$
\hfill(2.34)
\\[2ex]
is a first integral of system (2.10).

In Example 2.1 we proved that the system (2.10) is incomplete and has the defect $\delta=2.$
By Theorem 2.11, an integral basis of system (2.10) has the dimension  $n-m-\delta=5-2-2=1.$

Therefore the first integral (2.34) is an integral basis on the space $\R^5$ 
of the incomplete system (2.10).

By Corollary 2.2 and Theorem 2.15, it follows that 
the scalar function (2.34) is an integral basis on the space $\R^5$ 
both the incomplete system (2.11) and the complete system (2.12).
\vspace{0.5ex}

{\bf Example 2.6.}                                       
In Example 2.2 it has been shown that the system (2.23) is 
incomplete and has the defect $\delta=1.$

By Theorem 2.12, the incomplete system (2.23) 
on a neighbourhood of any point of a normalization domain 
has a basis of first integrals of dimension $n-m-\delta=5-2-1=2.$

From the definition of a first integral (Definition 2.1) it follows that  
the complete normal system (2.26) has the first integrals 
\\[1.75ex]
\mbox{}\hfill                                           
$
F^{}_{21}\colon x\to \,
x^{}_2\;\! x_1^{{}-1}$
\ for all $x\in H^{}_1\ \ \ 
\text{and}\ \ \ 
F^{}_{31}\colon x\to \,
x^{}_3\;\! x_1^{{}-1}$ 
\ for all $x\in H^{}_1\,,
$
\hfill(2.35)
\\[2ex]
where $H^{}_1$ is any domain from the set $\{x\colon x^{}_1\ne 0\}$ of the space $\R^5.$
\vspace{0.35ex}

By Theorem 2.15, the scalar functions (2.35) are 
a basis of first integrals on any domain $H^{}_1\subset\{x\colon x^{}_1\ne 0\}$
of the complete normal system (2.26), of the complete system (2.24), and of the incomplete systems (2.23).

Similarly, the scalar functions
\\[1.75ex]
\mbox{}\hfill                                           
$
F^{}_{12}\colon x\to \,
x^{}_1\;\! x_2^{{}-1}$
\ for all $x\in H^{}_2\ \ \ 
\text{and}\ \ \ 
F^{}_{32}\colon x\to \,
x^{}_3\;\! x_2^{{}-1}$
\ for all $x\in H^{}_2
$
\hfill(2.36)
\\[2ex]
are a basis of first integrals on any domain $H^{}_2\subset \{x\colon\! x^{}_2\ne 0\}$
of the complete normal system (2.27), of the complete system (2.24), and of the incomplete systems (2.23).
\vspace{0.25ex}

The scalar functions
\\[1.75ex]
\mbox{}\hfill                                           
$
F^{}_{13}\colon x\to \,
x^{}_1\;\! x_3^{{}-1}$
\ for all $x\in H^{}_3\ \ \ 
\text{and}\ \ \ 
F^{}_{23}\colon x\to \,
x^{}_2\;\! x_3^{{}-1}$
\ for all $x\in H^{}_3
$
\hfill(2.37)
\\[2ex]
are a basis of first integrals on any domain 
$H^{}_3\subset \{x\colon\! x^{}_3\ne 0\}$
of the complete normal system (2.28), of the complete system (2.24), and of the incomplete systems (2.23).

Each of the integral bases (2.35), (2.36), (2.37) for the incomplete system (2.23) (the com\-p\-le\-te system (2.24)) 
defines the two families of the integral surfaces,
which consistes of planes, for this system respectively 
\\[1ex]
\mbox{}\hfill
$
p^{}_{21}=
\{x\colon C^{}_1 x^{}_1+
C^{}_2 x_2=0\}
\, \ \ \quad \text{and}\quad  \ \ \,
p^{}_{31}=
\{x\colon C^{}_3 x^{}_3+
C^{}_4 x_1=0\},
\hfill
$
\\[2ex]
\mbox{}\hfill
$
p^{}_{12}=
\{x\colon C^{}_1 x^{}_1+
C^{}_2 x_2=0\}
\, \ \ \quad \text{and}\quad  \ \ \,
p^{}_{32}=
\{x\colon C^{}_3 x^{}_3+
C^{}_4 x_2=0\},
\hfill
$
\\[2ex]
\mbox{}\hfill
$
p^{}_{13}=
\{x\colon C^{}_1 x^{}_1+
C^{}_3 x_3=0\}
\, \ \ \quad \text{and}\quad  \ \ \,
p^{}_{23}=
\{x\colon C^{}_2 x^{}_2+
C^{}_4 x_3=0\},
\hfill
$
\\[1.75ex]
where $C^{}_1,\ldots,C^{}_4$ are arbitrary real constants.
\\[3.75ex]
\centerline{
\large\bf  
3. Dimension of  integral basis for not completely solvable 
}
\\[0.5ex]
\centerline{
\large\bf  
total differential system
}
\\[1.5ex]
\indent
The normal on the domain $\Pi\subset\R^{m+n}$ linear homogeneous partial differential system 
\\[1.75ex]
\mbox{}\hfill                                         
$
\displaystyle
\partial_{t_j^{}}^{}y=
{}-\sum\limits_{i=1}^{n} X_{ij}^{}(t,x)\, \partial_{x_i^{}}^{}y,
\quad 
j=1,\ldots,m,
$
\hfill (3.1)
\\[1.75ex]
is associated to the system of total differential equations {\rm(TD)}.
\vspace{0.35ex}

By Definitions 1.1, 2.1, and 2.5, we obtain
\vspace{0.5ex}

{\bf Theorem 3.1.}                                    
{\it
The total differential system {\rm(TD)}
is integrally equivalent on some domain $\Pi^\prime\subset\Pi$ to 
the normal linear homogeneous system of partial differential equations {\rm(3.1)}.
}
\vspace{0.5ex}

From Theorems 2.7 and 3.1, we have 
\vspace{0.5ex}

{\bf Theorem 3.2.}\!\!                                    
{\it
The total differential system \!{\rm(TD)}\! is completely solvable 
if and only if 
the normal linear homogeneous system of partial differential equations\! {\rm(3.1)}\! 
\vspace{0.5ex}
is complete\! {\rm(}\!jacobian{\rm)}\!.
}

Using Theorems 1.3, 2.9, 3.1, and 3.2, we can state the following
\vspace{0.5ex}

{\bf Theorem 3.3.}                                    
{\it
If the total differential system {\rm(TD)} with $X\in C^\infty(\Pi)$ is completely solvable 
{\rm(}the normal linear homogeneous system of partial differential equations {\rm(3.1)} 
with $X_{ij}^{}\in C^\infty(\Pi),\ i=1,\ldots,n,\ j=1,\ldots,m,$ is jacobian{\rm)}, 
then the systems {\rm(TD)} and {\rm(3.1)} have the same integral basis of dimension $n$ 
on some domain $\Pi^\prime\subset\Pi.$
}
\vspace{0.5ex}

The total differential system (1.15) is completely solvable. 
The jacobian normal linear homogeneous system of partial differential equations (2.31) is 
associated to the system (1.15). 
The integral basis of these systems was built in Examples 1.3 and 2.3. 
\vspace{0.75ex}

{\bf Example 3.1.}                                      
Let us consider the total differential system
\\[2ex]
\mbox{}\hfill                                          
$
dx_1^{}=
x_1^{}(x_1^{}+1)(dt_1^{}+dt_2^{}),\qquad
dx_2^{}=
x_2^{}(x_1^{}+2)(dt_1^{}+dt_2^{}),
\hfill
$
\\
\mbox{}\hfill(3.2)
\\
\mbox{}\hfill
$
dx_3^{}=
x_3^{}(x_1^{}+3)\,dt_1^{}+ x_3^{}(x_1^{}+5)\,dt_2^{}.
\hfill
$
\\[2ex]
\indent
The normal linear homogeneous system of partial differential equations 
\\[1.75ex]
\mbox{}\hfill                                          
$
{\frak X}_1^{}(t,x)\,y=0,
\qquad
{\frak X}_2^{}(t,x)\,y=0, 
$
\hfill(3.3)
\\[1.5ex]
where the linear differential operators of first order
\\[2ex]
\mbox{}\hfill
$
{\frak X}_1^{}(t,x)=
\partial_{t_1^{}}^{}+ x_1^{}(x_1^{}+1)\,\partial_{x_1^{}}^{}+
x_2^{}(x_1^{}+2)\,\partial_{x_2^{}}^{}+
x_3^{}(x_1^{}+3)\,\partial_{x_3^{}}^{}$
\ for all $(t,x)\in\R^5, 
\hfill
$
\\[2.5ex]
\mbox{}\hfill
$
{\frak X}_2^{}(t,x)=
\partial_{t_2^{}}^{}+ x_1^{}(x_1^{}+1)\,\partial_{x_1^{}}^{}+
x_2^{}(x_1^{}+2)\, \partial_{x_2^{}}^{}+
x_3^{}(x_1^{}+5)\,\partial_{x_3^{}}^{}$
\ for all $(t,x)\in\R^5,
\hfill
$
\\[2.15ex]
is associated to the total differential system (3.2). 

\newpage

Since the Poisson bracket 
\vspace{0.35ex}
$\bigl[ {\frak X}_1^{}(t,x), {\frak X}_2^{}(t,x) \bigr]= {\frak O}$
for all $(t,x)\in\R^5,$ we see that 
the system (3.2) is completely solvable and 
the system (3.2) is complete. 

The systems (3.2) and (3.3) have the same integral basis of dimension $r=3.$  
\vspace{0.25ex}

The functionally independent on any domain 
$H_1^{}\subset \Xi_1^{}=\{(t,x)\colon x_1^{}\ne 0\}$ scalar functions
\\[1.75ex]
\mbox{}\hfill
$
F_1^{}\colon(t,x)\to\,
x_2^{}(x_1^{}+1)\,x_1^{-2}$
\ for all $(t,x)\in\Xi_1^{}\,,
\hfill
$
\\[2.5ex]
\mbox{}\hfill
$
F_2^{}\colon(t,x)\to\,
(x_1^{}+1)\,x_1^{-1}\exp(t_1^{}+t_2^{})$
\ for all $(t,x)\in\Xi_1^{}\,,
$
\hfill (3.4)
\\[2.5ex]
\mbox{}\hfill
$
F_3^{}\colon(t,x)\to\,
x_3^{}x_1^{-1}
\exp({}-2(t_1^{}+2t_2^{}))$
\ for all $(t,x)\in\Xi_1^{}
\hfill
$
\\[2ex]
are a basis of first integrals on the domain $H_1^{}$ 
both for the system (3.2) and the system (3.3).

By the assumption of functional ambiguity of first integrals (Theorems 1.2 and 2.1),
we can build integral bases of the systems (3.2) and (3.3)
which are different from the integral basis (3.4).
Note also that domains of definition 
for these integral bases can contain the set of points 
$(t_1^{}\,,t_2^{}\,,0, x_2^{}\,,x_3^{}).$
\vspace{0.35ex}

For example, the scalar functions $F_1^{{}-1},\ F_2^{{}-1},$ and $F_3^{{}-1}$ on any domain
\vspace{0.25ex}
$H_2^{}$ from the set $\Xi_2^{}=\{(t,x)\colon x_1^{}\ne {}-1,\, x_2^{}\ne 0,\, x_3^{}\ne 0\}$
\vspace{0.25ex}
are a basis of first integrals on the domain $H_2^{}$ 
both for the system (3.2) and the system (3.3).
\vspace{0.5ex}

If the system {\rm(TD)} is not completely solvable, then the system (3.1) is incomplete. 
In this case, we reduce the system (3.1)  to the complete system.
We obtain the defect $\delta,\ 0<\delta\leq n,$ and a normalization domain of system (3.1).
By Theorem 2.12, we have the following
\vspace{0.25ex}

{\bf Theorem 3.4.}\!\!                                     
{\it
The not completely solvable total differential system {\rm(TD)} with $\!X\!\in\C^\infty(\Pi)\!\!$
on a neig\-h\-bo\-ur\-ho\-od of any point of a normalization domain for the
linear homogeneous {\rm(}incomplete{\rm)} system of partial differential equations {\rm(3.1)}
has an integral basis of dimension $n-\delta,$ where $\delta$ is the defect of system {\rm(3.1).}
This basis of first integrals for system {\rm(TD)} is also a basis of first integrals for system {\rm(3.1)}. 
}
\vspace{0.75ex}

{\bf Example 3.2.}                                      
The normal linear homogeneous system of partial differential equations 
\\[1.75ex]
\mbox{}\hfill                                          
$
{\frak X}_1^{}(t,x)\,y=0,
\qquad
{\frak X}_2^{}(t,x)\,y=0, 
$
\hfill(3.5)
\\[1.75ex]
where ${\frak X}_1^{}$ and ${\frak X}_1^{}$ are the 
linear differential operators of first order (1.4) and (1.5) respectively, 
is associated to the total differential system (1.3). 

The Poisson bracket
\\[1.75ex]
\mbox{}\hfill
$
{\frak X}_{12}^{}(t,x)=
\bigl[{\frak X}_1^{}(t,x), {\frak X}_2^{}(t,x) \bigr]=\,
t_1^{}x_2^{}x_3^{}\, \partial_{x_1^{}}^{} +\,
x_3^{}(1+x_1^{}t_1^{-1}- x_2^2\,t_1^{-1})\, \partial_{x_2^{}}^{} \ +
\hfill
$
\\[2.5ex]
\mbox{}\hfill
$
+ \ x_2^{}({}-2-2x_1^{}t_1^{-1}+
2x_1^2\, t_1^{-2}+x_2^2+x_3^2)\, \partial_{x_3^{}}^{}$
\ for all $(t,x)\in D
\hfill
$
\\[2ex]
is not the null operator. Therefore the system (1.3) is not completely solvable 
and the system (3.5) is incomplete.

Since the Poisson bracket
\\[1.75ex]
\mbox{}\hfill
$
{\frak X}_{1;\;\! 12}^{}(t,x)=
\bigl[
{\frak X}_1^{}(t,x),
{\frak X}_{12}^{}(t,x)
\bigr]=
\hfill
$
\\[2ex]
\mbox{}\hfill
$
={}-2x_3^{} (t_1^{}-t_1^{}x_2^2- x_1^2\, t_1^{-1}+x_1^{})\,
\partial_{x_1^{}}^{}+\,
x_2^{} x_3^{} (1-5x_1^{} t_1^{-1}+ x_1^2\, t_1^{-2})\, \partial_{x_2^{}}^{}\ +
\hfill
$
\\[2ex]
\mbox{}\hfill
$
+ \  (2+4x_1^{}t_1^{-1}-2x_1^2\, t_1^{-2}-
4x_1^3 \,t_1^{-3}+ 2x_1^4\, t_1^{-4}+
x_1^{}x_2^2\, t_1^{-1}+
3x_1^2 x_2^2\, t_1^{-2}+
2x_2^4-5x_2^2 \ -
\hfill
$
\\[2ex]
\mbox{}\hfill
$
-\ 2x_1^{}x_3^2\, t_1^{-1}+ 2x_1^2 x_3^2\, t_1^{-2}+
2x_2^2 x_3^2-2x_3^2)\,\partial_{x_3^{}}^{}$
\ for all $(t,x)\in D
\hfill
$
\\[2ex]
is not a linear combination of the operators ${\frak X}_1^{},\, {\frak X}_2^{},\, {\frak X}_{12}^{},$
we see that the linear homogeneous system of partial differential equations 
\\[2ex]
\mbox{}\hfill                                          
$
{\frak X}_1^{}(t,x)\,y=0,
\qquad
{\frak X}_2^{}(t,x)\,y=0,
\qquad
{\frak X}_{12}^{}(t,x)\,y=0 
$
\hfill(3.6)
\\[1.25ex]
is incomplete. 

In Example 1.1, we proved that the function (1.6) is a first integral on a domain $\Pi\subset D$ of system (1.3). 
Therefore the incomplete system  (3.6) has the defect $\delta=1$ and 
the function (1.6) is an integral basis on a domain $\Pi\subset D$ of system (3.6).

Thus the function (1.6) is an integral basis on a domain $\Pi\subset D$
both for the not completely solvable total differential system (1.3) and 
the incomplete normal system of partial differential equations (3.5) with the defect $\delta=2.$
\vspace{0.5ex}

Recall that a complete linear homogeneous system of partial differential equations has 
the defect $\delta=0.$ Hence for the completely solvable total differential system {\rm(TD)}
and for the not completely solvable total differential system {\rm(TD)} we have the following assertion
\vspace{0.5ex}

{\bf Theorem 3.5.}\!\!                                   
{\it
The system {\rm(TD)} with $\!X\!\in\! C^\infty(\Pi)\!$ and the 
associated system {\rm(3.1)} to the system {\rm(TD)} have the same integral basis of 
dimension $n-\delta$  on a neig\-h\-bo\-ur\-ho\-od of any point of a normalization domain for 
system {\rm(3.1)}, where $\delta,\, 0\leq\delta\leq n,\!$ is the defect of system {\rm (3.1)}. 
}

Using this notation, we can state the definitions of defect and of normalization domain for 
completely solvable system {\rm(TD)} and for not completely solvable  system {\rm(TD)}.
\vspace{0.35ex}

{\bf Definition 3.1.}\!                               
{\it
The total differential system {\rm(TD)} has the defect $\delta,\,  0\leq\delta\leq n,\!$ 
where $\delta\!$ is the defect of the associated linear homogeneous system of partial differential equations {\rm(3.1)}.
Thus a normalization domain of system {\rm(3.1)} is called  
\vspace{0.35ex}
a normalization domain of system\! {\rm(TD\!)}\!.
}

If the system (TD) is completely solvable, then a normalization domain of this system is 
the domain of complete solvability for system (TD).

By Definition 3.1 and Theorem 3.5, we obtain
\vspace{0.35ex}

{\bf Theorem 3.6.}                                   
{\it
Suppose the total differential system {\rm(TD)} with $X\in C^\infty(\Pi)$ 
has the defect $\delta,\ 0\leq\delta\leq n.$ Then this system has 
an integral basis of dimension $n-\delta$ on a neig\-h\-bo\-ur\-ho\-od of any point of a normalization domain. 
}
\vspace{0.5ex}

{\bf Example 3.3.}                                      
The total differential system 
\\[2ex]
\mbox{}\hfill                                          
$
dx_1^{}=x_1^{}\, d t_1^{}+ x_1^2\, d t_2^{}\,,
\qquad
dx_2^{}=x_2^2\, d t_1^{}+ x_2^3\, d t_2^{}
$
\hfill(3.7)
\\[2ex]
induces the linear differential ope\-ra\-tors of first order
\\[2ex]
\mbox{}\hfill
$
{\frak X}_1^{}(t,x)=
\partial_{t_1^{}}^{}+ x_1^{}\,\partial_{x_1^{}}^{}+x_2^2\,\partial_{x_2^{}}^{}
$
\ for all $(t,x)\in\R^4,
\hfill
$
\\[2.5ex]
\mbox{}\hfill
$
{\frak X}_2^{}(t,x)=
\partial_{t_2^{}}^{}+x_1^2\,\partial_{x_1^{}}^{}+x_2^3\,\partial_{x_2^{}}^{}
$
for all $(t,x)\in\R^4.
\hfill
$
\\[2ex]
\indent
Since the Poisson bracket
\\[2ex]
\mbox{}\hfill
$
{\frak X}_{12}^{}(t,x)=
\bigl[ {\frak X}_1^{}(t,x), {\frak X}_2^{}(t,x)\bigr]=
x_1^2\,\partial_{x_1^{}}^{}+ x_2^4\,\partial_{x_2^{}}^{}
$
\ for all $(t,x)\in\R^4
\hfill
$
\\[2.25ex]
is not the null operator, we see that the system (3.7) is not completely solvable. 

The linear differential ope\-ra\-tors of first order
${\frak X}_1^{}\,,\,
{\frak X}_2^{}\,,\,
{\frak X}_{12}^{}\,,\,
{\frak X}_{2;\;\! 12}^{}\,,$
where 
\\[2ex]
\mbox{}\hfill
$
{\frak X}_{2;\;\! 12}^{}(t,x)=
\bigl[
{\frak X}_2^{}(t,x),
{\frak X}_{12}^{}(t,x)
\bigr]=\,
x_2^6\,\partial_{x_2^{}}^{}
$
\ for all $(t,x)\in\R^4, 
\hfill
$
\\[2.25ex]
are not linearly bound on $\R^4.$ 
Therefore the associated incomplete normal linear homogeneous system of partial differential equations
\\[1.75ex]
\mbox{}\hfill
$
{\frak X}_1^{}(t,x)\,y=0,
\qquad
{\frak X}_2^{}(t,x)\,y=0
\hfill
$
\\[1.5ex]
to the system (3.7) has the defect $\delta=2.$

Thus the system (3.7) has the defect $\delta=2.$ 
From $n-\delta=2-2=0$ it follows that the system (3.7) has no first integrals.
\vspace{0.5ex}

{\bf Example 3.4.}                                      
The normal linear homogeneous system of partial differential equations
\\[2ex]
\mbox{}\hfill                                          
$
{\frak X}_1^{}(t,x)\,y \, \equiv\,
\partial_{t_1^{}}^{}y+ x_1^{}\,\partial_{x_1^{}}^{}y+
(1+x_1^{}+2x_2)\, \partial_{x_2^{}}^{}y=0,
\hfill
$
\\
\mbox{}\hfill (3.8)
\\
\mbox{}\hfill
$
{\frak X}_2^{}(t,x)\,y \, \equiv\,
\partial_{t_2^{}}^{}y+3x_1^{}\,\partial_{x_1^{}}^{}y+
(x_1^{}+3x_2)\,\partial_{x_2^{}}^{}y=0
\hfill
$
\\[2.25ex]
is associated to the not completely solvable total differential system (1.7). 
Therefore the system (3.8) is incomplete.

The Poisson brackets
\\[2ex]
\mbox{}\hfill
$
{\frak X}_{12}^{}(t,x)=
\bigl[
{\frak X}_1^{}(t,x),
{\frak X}_2^{}(t,x)
\bigr]=
(3-x_1^{})\,\partial_{x_2^{}}^{}
$
\ for all $(t,x)\in\R^4,
\hfill
$
\\[2.5ex]
\mbox{}\hfill
$
\bigl[
{\frak X}_1^{}(t,x),
{\frak X}_{12}^{}(t,x)
\bigr]=
(x_1^{}-6)(3-x_1^{})^{-1}\,
{\frak X}_{12}^{}(t,x)
$
\ for all $(t,x)\in
\{(t,x)\colon x_1^{}\ne 3\},
\hfill
$
\\[2.5ex]
\mbox{}\hfill
$
\bigl[
{\frak X}_2^{}(t,x),
{\frak X}_{12}^{}(t,x)
\bigr]=
9(x_1^{}-3)^{-1}\,
{\frak X}_{12}^{}(t,x)
$
\ for all $(t,x)\in
\{(t,x)\colon x_1^{}\ne 3\}.
\hfill
$
\\[2ex]
Thus the system (3.8) has the defect $\delta=1.$
\vspace{0.35ex}

By Definition 3.1, the system (1.7) has the defect $\delta=1.$ 
By Theorem 3.6, the dimension of an integral basis of system (1.7) is $n-1=2-1=1.$
Therefore the function (1.8) is an integral basis on the space $\R^4$ of system (1.7).

By Theorem 3.5, the function (1.8) is an integral basis on the space $\R^4$ 
of the system of partial differential equations (3.8).
\vspace{0.75ex}

{\bf Example 3.5.}                                      
The normal linear homogeneous system of partial differential equations (2.10) has 
the defect $\delta=2$ and the integral basis (2.34). 
This system is associated to the total differential system 
\\[1.5ex]
\mbox{}\hfill                                          
$
dx_3^{}=
2 x_3^{} x_5^{}\, d x_2^{}\,,
\qquad
dx_4^{}=x_5^{}\, d x_1^{}+
2 x_4^{} x_5^{}\, d x_2^{}\,,
\hfill
$
\\[0.1ex]
\mbox{}\hfill (3.9)
\\[0.1ex]
\mbox{}\hfill
$
dx_5^{}={}-x_4^{}\, d x_1^{}+
(1-x_3^2-x_4^2+x_5^2)\,d x_2^{}\,.
\hfill
$
\\[2ex]
\indent
Therefore the autonomous system (3.9) is not completely solvable and has the defect $\delta=2.$
The function (2.34) is an autonomous first integral of system (3.9) and this function is 
an integral basis on the space $\R^5$ of system (3.9).
\vspace{0.5ex}

Suppose the system {\rm(TD)} has the defect $\delta,\ 0\leq\delta<n.$ 
Then the associated system (3.1) to the system {\rm(TD)} we reduced to  
the integrally equivalent complete system
\\[2ex]
\mbox{}\hfill                                         
$
\displaystyle
\partial_{t_j^{}}^{}y+
\sum\limits_{i=1}^n\,
X_{ij}^{}(t,x)\,\partial_{x_i^{}}^{}y=0,
\quad \, 
j=1,\ldots,m,
\hfill
$
\\
\mbox{}\hfill(3.10)
\\
\mbox{}\hfill
$
\displaystyle
\sum\limits_{i=1}^n\,
{X}_{i\;\!\nu}^{{}^{\scriptstyle *}}(t,x)\,
\partial_{x_i^{}}^{}y=0,
\quad \, 
\nu=1,\ldots,\delta,
\hfill
$
\\[2ex]
where the functions 
\vspace{0.35ex}
${X}_{i\;\!\nu}^{{}^{\scriptstyle *}}\colon\Pi\to\R,\ 
i=1,\ldots,n,\ \nu=1,\ldots,\delta,$ are constructed on the base of the functions 
$X_{ij}^{},\, i=1,\ldots,n,\, j=1,\ldots,m,$ by the rule (2.9).
\vspace{0.5ex}

We reduced the system (3.10) to a normal system and then for this normal system we build 
the associated total differential system
\\[2ex]
\mbox{}\hfill                                          
$
\displaystyle
dx_{k_\gamma^{}}^{}\, =\,
\sum\limits_{j=1}^m\,
G_{k_\gamma^{}\,j}^{}(t,x)\, dt_j^{}\ +\,
\sum\limits_{\mu=n-\delta+1}^n
G_{k_\gamma^{}\,k_\mu^{}}^{}(t,x)\, dx_{k_\mu^{}}^{}\,,
\hfill
$
\\
\mbox{}\hfill (3.11)
\\[0.25ex]
\mbox{}\ \ \,
$
\gamma=1,\ldots,n-\delta,\ \ 
k_\gamma^{},\,
k_\mu^{}\in\{1,\ldots,n\},\ \,
k_i^{}\ne k_\xi^{}\,,
\ \ \ 
i=1,\ldots,n,\ \,
\xi=1,\ldots,n,\ \,
i\ne\xi.
\hfill
$
\\[2.25ex]
\indent
The system (3.11) is completely solvable on a normalization domain of system (3.10)
(this nor\-ma\-li\-za\-tion domain is a normalization domain both for 
the partial differential system (3.1) 
\vspace{0.25ex}
and the total differential system (TD)).

The system (3.11) in relation to the system (TD) has 
the extended coordinate space $Ot$ on $\delta$ coordinates 
at the expense of $\delta$ coordinates of the coordinate space $Ox.$
\vspace{0.25ex}

Using the rearrangement of the dependent and independent variables in the system (TD), we get 
the system (3.11) has the form
\\[2ex]
\mbox{}\hfill                                         
$
\displaystyle
dx_i^{}=
\sum\limits_{j=1}^{m+\delta}\,
H_{ij}^{}(t_1^{}\,,\ldots,
t_{m+\delta}\,,x_1^{}\,,\ldots,
x_{n-\delta})\,dt_j^{}\,,\ 
\quad
i=1,\ldots,n-\delta,
$
\hfill (3.12)
\\[2ex]
where $t_{m+\nu}^{}=x_{n-\delta+\nu}^{}\,,\ \nu=1,\ldots,\delta.$
\vspace{0.5ex}

Note also that the system (TD) and the system (3.12) are integrally equivalent 
(they have the same integral basis) on 
a normalization domain (the system (TD) on this domain is reduced to the system (3.12)) of system (TD).
\vspace{0.75ex}

{\bf Theorem 3.7.}                                    
{\it
The total differential system {\rm(TD)} with the coordinate spaces $Ot$ and $Ox$ and 
with the defect $\delta,\ 0\leq\delta<n,$ on a normalization domain is integrally equivalent 
to the completely solvable total differential system with the coordinate spaces 
$Ot_1^{}\,,\ldots,t_{m+\delta}^{}$ and 
$Ox_1^{}\,,\ldots,x_{n-\delta}^{}\,,$ 
\vspace{0.5ex}
where $t_{m+\nu}^{}=x_{n-\delta+\nu}^{}\,,\ \nu=1,\ldots,\delta$ 
{\rm(}accurate to numbering of the dependent variables $x_1^{}\,,\ldots,x_n^{}$ in the system {\rm(TD))}.
}
\vspace{0.75ex}

{\bf Example 3.6.}                                       
\vspace{0.5ex}
Adding the equation ${\frak X}_{12}^{}(t,x)\,y\equiv
(3-x_1^{})\,\partial_{x_2^{}}^{}y=0$ to the incomplete system  (3.8), we get 
the integrally equivalent complete normal system 
\\[2ex]
\mbox{}\hfill                                          
$
\partial_{t_1^{}}^{}y=
{}-x_1^{}\,\partial_{x_1^{}}^{}y,
\qquad
\partial_{t_2^{}}^{}y={}-3x_1^{}\,\partial_{x_1^{}}^{}y,
\qquad
\partial_{x_2^{}}^{}y=0.
$
\hfill(3.13)
\\[2.25ex]
\indent
The total differential equation
\\[1.5ex]
\mbox{}\hfill                                          
$
d x_1^{}=x_1^{}\,d t_1^{}+ 3x_1^{}\,d t_2^{}+0\,d x_2^{}
$
\hfill(3.14)
\\[1.75ex]
is associated to the system (3.13).

Therefore the not completely solvable total differential system (1.7) is 
integrally equivalent to the completely solvable total differential equation (3.14).
Moreover, the system (1.7) and the equation (3.14) have the same integral basis, which is 
the first integral (1.8).

In Example 1.2, we proved that the system (1.8) has no solutions.
At the same time the integrally equivalent equation (3.14) to the system (1.8) has the general solution
\\[1.75ex]
\mbox{}\hfill
$
x_1^{}\colon
(t_1^{}\,,t_2^{}\,,x_2^{})\to
C\exp(t_1^{}+3t_2^{})
$
\ for all $(t_1^{}\,,t_2^{}\,,x_2^{})\in\R^3.
\hfill
$
\\[2ex]
\indent
Thus integrally equivalent total differential systems have the same first integrals.
But this statement is not true for solutions of integrally equivalent systems.
\vspace{0.5ex}

{\bf Example 3.7.}                                       
The complete system of partial differential equations (2.12) is reduced to the normal systems:
\\[2ex]
\mbox{}\hfill
$
\partial_{x_1^{}}^{}y=0,
\qquad
\partial_{x_2^{}}^{}y=0,
\qquad
\partial_{x_3^{}}^{}y=
{}-\dfrac{1-x_3^2+x_4^2+x_5^2}
{2x_3^{}x_5^{}}\ 
\partial_{x_5^{}}^{}y,
\qquad
\partial_{x_4^{}}^{}y=
\dfrac{\,x_4^{}}{\,x_5^{}}\ \partial_{x_5^{}}^{}y\,;
\hfill
$
\\[2.5ex]
\mbox{}\hfill
$
\partial_{x_1^{}}^{}y=0,
\qquad
\partial_{x_2^{}}^{}y=0,
\qquad
\partial_{x_3^{}}^{}y=
{}-\dfrac{1-x_3^2+x_4^2+x_5^2}
{2x_3^{}x_4^{}}\ 
\partial_{x_4^{}}^{}y,
\qquad
\partial_{x_5^{}}^{}y=
\dfrac{\,x_5^{}}{\,x_4^{}}\ 
\partial_{x_4^{}}^{}y\,;
\hfill
$
\\[2.5ex]
\mbox{}\hfill
$
\partial_{x_1^{}}^{}y=0,
\quad 
\partial_{x_2^{}}^{}y=0,
\quad 
\partial_{x_4^{}}^{}y=
{}-\dfrac{2x_3^{}x_4^{}}
{1-x_3^2+x_4^2+x_5^2}\ 
\partial_{x_3^{}}^{}y,
\quad 
\partial_{x_5^{}}^{}y=
{}-\dfrac{2x_3^{}x_5^{}}
{1-x_3^2+x_4^2+x_5^2}\ 
\partial_{x_3^{}}^{}y.
\hfill
$
\\[2ex]
\indent
The completely solvable total differential equations
\\[1.75ex]
\mbox{}\hfill                                     
$
d x_5^{}=0\,d x_1^{}+0\,d x_2^{}-
\dfrac{1-x_3^2+x_4^2+x_5^2}
{2x_3^{}x_5^{}}\ d x_3^{}+
\dfrac{\,x_4^{}}{\,x_5^{}}\ 
d x_4^{}\,,
$
\hfill(3.15)
\\[2.25ex]
\mbox{}\hfill                                     
$
d x_4^{}=0\,d x_1^{}+0\,d x_2^{}-
\dfrac{1-x_3^2+x_4^2+x_5^2}
{2x_3^{}x_4^{}}\ d x_3^{}+
\dfrac{\,x_5^{}}{\,x_4^{}}\ 
d x_5^{}\,,
$
\hfill(3.16)
\\[2.25ex]
\mbox{}\hfill                                     
$
d x_3^{}=0\,d x_1^{}+0\,d x_2^{}-
\dfrac{2x_3^{}x_4^{}}
{1-x_3^2+x_4^2+x_5^2}\ d x_4^{}-
\dfrac{2x_3^{}x_5^{}}
{1-x_3^2+x_4^2+x_5^2}\ d x_5^{}
$
\hfill (3.17)
\\[1.5ex]
are associated to these normal systems respectively.

Thus the not completely solvable total differential system (3.9) is:
\vspace{0.25ex}

a)\, integrally equivalent on any domain
\vspace{0.35ex}
$H_5^{}\subset \{x\colon x_3^{}\ne 0,\, x_5^{}\ne 0\}$
to the completely solvable total differential equation (3.15) and they have 
the same integral basis on the domain $H_5^{}\,,$ which is the first integral (2.34);
\vspace{0.25ex}

b)\, integrally equivalent on any domain
\vspace{0.35ex}
$H_4^{}\subset \{x\colon x_3^{}\ne 0,\, x_4^{}\ne 0\}$
to the completely solvable total differential equation (3.16) and they have 
the same integral basis on the domain $H_4^{}\,,$ which is the first integral (2.34);
\vspace{0.25ex}

c)\, integrally equivalent on any domain
\vspace{0.35ex}
$H_3^{}\subset \{x\colon 1-x_3^2+x_4^2+x_5^2\ne 0\}$
to the completely solvable total differential equation (3.17) and they have 
the same integral basis on the domain $H_3^{}\,,$ which is the first integral (2.34).
\\[5.75ex]
\centerline{
\large\bf  
4. First integrals for Pfaff system of equations
}
\\[1.75ex]
\indent
{\bf  4.1. Integrally equivalent Pfaff systems of equations }
\\[1.25ex]
\indent
{\bf Definition 4.1.}
{\it
A scalar function $F\in C^1(G^{\;\!\prime})$
is said to be a 
\textit{\textbf{first integral}} on a domain $G^{\;\!\prime}\subset G$ 
of system ${\rm (Pf)}$ with $\omega^{}_j\in C(G),\  j=1,\ldots,m,$ if 
there exist the scalar functions $a^{}_j\in C(G^\prime),\ j=1,\ldots,m,$ such that 
the total differential
}
\\[1ex]
\mbox{}\hfill                                           
$
\displaystyle
d F(x)=\sum\limits_{j=1}^m\,
a_j^{}(x)\;\!\omega_j^{}(x)
$
\ for all $x\in G^\prime.
$
\hfill(4.1)
\\[2ex]
\indent
Let us introduce the equivalence relation on a set of Pfaff systems of equations.
\vspace{0.5ex}

{\bf Definition 4.2.}
{\it
We'll say that two Pfaff systems of equations are \textit{\textbf{integrally equivalent}} 
on some domain if  on this domain
each first integral of the first system is a first integral of the second system 
and on the contrary each first integral of the second system 
is a first integral of the first system.
}
\vspace{0.5ex}

We claim that the demand that the 1-forms (0.2) are not linearly bound on the domain $G$ 
is not narrow the set of all possible Pfaff systems of equations (Pf) 
(from the point of view of the integral equivalence).
Indeed, let the 1-forms (0.2) be linearly bound on the domain $G.$
Then the functional matrix
\\[1.5ex]
\mbox{}\hfill                                           
$
w(x)= \bigl\|
w_{ji}^{}(x)\bigr\|_{m\times n}$
\ for all $x\in G
$
\hfill(4.2)
\\[1.25ex]
has the rank
\\[1.25ex]
\mbox{}\hfill
$
{\rm rank}\,w(x)=s(x),\ 
1\leq s(x)<\min\,\{m,n\}$
\ for all $x\in G.
\hfill
$
\\[2.25ex]
\indent
Take $s=\min\{s(x)\colon x\in G\}$ not linearly bound on a domain $\Omega\subset G$ 1-forms
\\[1.75ex]
\mbox{}\hfill                                            
$
\omega_{j_l^{}}^{}\,,
\ \ 
j_l^{}\in\{1,\ldots,m\},
\ \ \ 
l=1,\ldots,s.
$
\hfill(4.3)
\\[2ex]
The domain $\Omega$ is such that the complement $\Omega$ on  $G$
has the null measure: $\mu{\rm C}_{G}^{}\Omega=0.$
\vspace{0.25ex}

Using the 1-forms (4.3), we get the new Pfaff system of equations
\\[1.75ex]
\mbox{}\hfill                                            
$
\omega_{j_l^{}}^{}(x)=0,
\ \
j_l^{}\in\{1,\ldots,m\},
\ \ \ 
l=1,\ldots,s.
$
\hfill(4.4)
\\[2ex]
\indent
Since the 1-forms (4.3) are not linearly bound and $\!s\!=\!\min\{s(x)\colon x\!\in\! G\},\!$ 
we see that the Pfaff systems of equations\! (Pf)\!  and\! (4.4) have the same 
first integrals on the domain $\!\Omega\!$ (by De\-fi\-ni\-ti\-on 4.1). 
Thus the systems (Pf) and (4.4) are integrally equivalent on the domain $\!\Omega.\!\!$ 

Since ${\rm C}_{G}^{}\Omega$ has the null measure, we see that a class of systems (Pf) 
doesn't restrict.
\\[2ex]
\indent
{\bf
4.2. Integral basis}
\\[0.75ex]
\indent
Suppose the 1-forms (0.2) are not linearly bound on the domain $G.$ 
Then the fun\-c\-tional matrix (4.2) has ${\rm rank}\,w(x)=m$ almost everywhere on the domain $G.$
To be definite, assume that $m\leq n.$
\vspace{0.25ex}

If $m=n,$ then, since the linear differential forms (0.2) are not linearly bound on the domain $G,$ 
it follows that the square matrix (4.2) of order $n$ is nonsingular almost everywhere on the domain $G.$

In this case the system (Pf) on a domain $\Omega\subset G$ 
by a nonsingular algebraic transformation can be reduced to the differential system
\\[1.75ex]
\mbox{}\hfill
$
dx_i^{}=0,
\quad 
i=1,\ldots,n,
\hfill
$
\\[2ex]
where the domain $\Omega$ is such that $\mu{\rm C}_{G}^{}\Omega=0.$

Whence, we obtain $x_i^{}=C_i^{}\,,\ i=1,\ldots,n,$
where $C_1^{}\,,\ldots, C_n^{}$ are arbitrary real constants.
Hence the functions 
\\[2ex]
\mbox{}\hfill                                           
$
F_i^{}\colon x\to x_i^{}$
\ for all $x\in\Omega,
\ \ \  
i=1,\ldots,n,
$
\hfill(4.5)
\\[2.25ex]
are first integrals on the domain $\Omega$ of system (Pf).

Thus the case $m=n$ is singular and here we have
\vspace{0.5ex}

{\bf Property 4.1.}                                     
{\it 
The system {\rm(Pf)} with $m=n$ has $n$ functionally independent first integrals {\rm(4.5)}
on such a subdomain $\Omega$  of the domain $G$ that $\mu{\rm C}_{G}^{}\Omega=0.$
}
\vspace{0.75ex}

{\bf Theorem 4.1.}                                     
{\it 
Let the functions {\rm(2.2)} be first integrals on a domain $G^\prime\subset G$ of 
system {\rm(Pf)} with $\omega_j^{}\in C(G).$ Then the function {\rm(2.4)} is also 
a first integral on the domain $G^\prime$ of system} (Pf).
\vspace{0.35ex}

{\sl Proof}.
By Definition 4.1, the functions (2.2) are first integrals on the domain $G^\prime$ of system (Pf)
if and only if 
there exist the scalar functions $a_{\xi j}^{}\in C(G^\prime),\ \xi=1,\ldots,k,\ j=1,\ldots,m,$
such that the total differentials of functions (2.2) have the forms
\\[1.5ex]
\mbox{}\hfill                                          
$
\displaystyle
dF_\xi^{}(x)=
\sum\limits_{j=1}^m\,
a_{\xi j}^{}(x)\;\!
\omega_j^{}(x)$
\ for all $x\in G^\prime,
\ \
\xi=1,\ldots,k.
$
\hfill(4.6)
\\[1.5ex]
\indent
Suppose $\!\Phi\!$ is arbitrary scalar function from the space $\!C^1({\rm E}F),\!$
where $\!F\!$ is the vector fun\-c\-ti\-on\! (2.3).
Then, using  the identities (4.6), we get the total differential of function\! (2.4)\! is 
\\[2ex]
\mbox{}\hfill
$
\displaystyle
d\Psi(x)=d\Phi
\bigl(F_1^{}(x),\ldots,F_k^{}(x)\bigr)=
\sum\limits_{\xi=1}^k \,
\partial_{F_\xi^{}}^{}
\Phi\bigl(F_1^{}\,,\ldots,F_k^{}\bigr)
_{\displaystyle 
\mbox{}|_{F=F(x)}}\,
dF_\xi^{}(x)=
\hfill
$
\\[2ex]
\mbox{}\hfill
$
\displaystyle
=\ \sum\limits_{\xi=1}^k\,
\sum\limits_{j=1}^m\,
\partial_{F_\xi^{}}^{}
\Phi\bigl(F_1^{}\,,\ldots,F_k^{}\bigr)
_{\displaystyle 
\mbox{}|_{F=F(x)}}\,
a_{\xi j}^{}(x)\;\!
\omega_j^{}(x)$
\ for all $x\in G^\prime.
\hfill
$
\\[2.25ex]
\indent
By this identity and Definition 4.1, it follows that 
the function (2.4) is a first integral on the domain $G^\prime$  
of the Pfaff system of equations (Pf). \k
\vspace{0.5ex}

It was shown in Theorem 4.1 that first integrals for a Pfaff system of equations are functional ambiguous.
Thus the priority of functionally independent  first integrals is installed.

The same property of functional ambiguous of first integrals we have 
for  systems of ordinary differential equations [95, pp. 262 -- 263],
for total differential systems (see Subsection~1.2),
for linear homogeneous partial differential equations [52, p. 16], 
and for linear homogeneous systems of partial differential equations (Theorem 2.1).
\vspace{0.5ex}

{\bf Example 4.1.}                                    
Consider the Pfaff system of equations
\\[1.75ex]
\mbox{}\hfill                                        
$
\omega_1^{}(x)=0,
\quad
\omega_2^{}(x)=0,
$
\hfill(4.7)
\\[1.75ex]
where the linear differential forms
\\[2ex]
\mbox{}\hfill
$
\omega_1^{}(x)=
x_1^{}(1+x_2^{})\;dx_1^{}+ x_2^{}(1-x_2^{})\;dx_2^{}+
(x_3^{}+x_2^{}x_4^{})\;dx_3^{}+ (x_4^{}+x_2^{}x_3^{})\;dx_4^{}$
\ for all $x\in\R^4,
\hfill
$
\\[2.25ex]
\mbox{}\hfill
$
\omega_2^{}(x)=
x_1^{}\;dx_1^{}-
x_2^{}\;dx_2^{}+
x_4^{}\;dx_3^{}+
x_3^{}\;dx_4^{}$
\ for all $x\in\R^4.
\hfill
$
\\[1.25ex]
\indent
We have
\\[1.25ex]
\mbox{}\hfill
$
2\;\!\omega_1^{}(x)+2(1-x_2^{})\;\!\omega_2^{}(x)=
d\bigl(2x_1^2+(x_3^{}+x_4^{})^2\bigr)$
\ for all $x\in\R^4,
\hfill
$
\\[2.25ex]
\mbox{}\hfill
$
2\;\!\omega_1^{}(x)- 2(1+x_2^{})\;\!\omega_2^{}(x)=
d\bigl(2x_2^2+(x_3^{}-x_4^{})^2\bigr)$
\ for all $x\in\R^4.
\hfill
$
\\[2ex]
\indent
Therefore, by Definition 4.1, the functions
\\[1.75ex]
\mbox{}\hfill                                        
$
F_1^{}\colon x\to 
2x_1^2+(x_3^{}+x_4^{})^2$
\ for all $x\in\R^4,
$
\hfill(4.8)
\\[2.5ex]
\mbox{}\hfill                                        
$
F_2^{}\colon x\to 
2x_2^2+(x_3^{}-x_4^{})^2$
\ for all $x\in\R^4
$
\hfill(4.9)
\\[2.25ex]
are first integrals on the space $\R^4$ of system (4.7). The first integrals (4.8) and (4.9) are 
functionally independent on the space $\R^4.$

By Theorem 4.1, the function 
\\[2ex]
\mbox{}\hfill                                        
$
F_3^{}\colon x\to 
x_1^2-x_2^2+2x_3^{}x_4^{}$
\ for all $x\in\R^4
$
\hfill(4.10)
\\[2.25ex]
is a first integral on the space $\R^4$ of the Pfaff system of equations (4.7).
Indeed, the function 
$F_3^{}(x)=\bigl(F_1^{}(x)-F_2^{}(x)\bigr)\slash 2$ for all $x\in\R^4.$
\vspace{0.75ex}

{\bf Definition 4.3.}
{\it
A set of functionally independent first integrals on the domain 
$G^{\;\!\prime}\subset G$ of system {\rm(Pf)} with $\omega_j^{}\in C(G),\ j=1,\ldots,m,$ is called a 
\textit{\textbf{basis of first integrals}} {\rm(}\textit{\textbf{integral basis}}{\rm)}
on the domain $G^{\;\!\prime}$ of system {\rm(Pf)} if
for any first integral $\Psi$ on the domain $G^{\;\!\prime}$ of system {\rm(Pf)},
we have $\Psi(x) = \Phi(F(x))$ for all $x \in G^{\;\!\prime},$
where $\Phi$ is some function of class $C^{1}_{}({\rm E}F),\ 
{\rm E}F$ is the range of the vector function  {\rm (2.3)}.
The number k is said to be the \textit{\textbf{dimension}} of 
ba\-sis of first integrals on the domain $G^{\;\!\prime}$ of system {\rm(Pf)}.
}
\vspace{0.75ex}

From Definition 4.3 and Property 4.1, we get the following
\vspace{0.5ex}

{\bf Property 4.2.}                                  
{\it 
The scalar functions {\rm(4.5)} are a basis of first integrals on a domain $\Omega\subset G$
for the Pfaff system of equations {\rm(Pf)} with $m=n.$
}
\vspace{0.75ex}

{\bf Example 4.2.}                                    
The Pfaff system of equations
\\[2ex]
\mbox{}\hfill                                        
$
\omega_1^{}(x)=0,
\quad 
\omega_2^{}(x)=0,
\quad 
\omega_3^{}(x)=0
$
\hfill(4.11)
\\[2.25ex]
with the linear differential forms 
\\[2ex]
\mbox{}\hfill
$
\omega_1^{}(x)=
dx_1^{}+dx_2^{}+2\,dx_3^{}$
\ for all $x\in\R^3,
\hfill
$
\\[2ex]
\mbox{}\hfill
$
\omega_2^{}(x)=
dx_1^{}+2\,dx_2^{}+2\,dx_3^{}$
\ for all $x\in\R^3,
\hfill
$
\\[2ex]
\mbox{}\hfill
$
\omega_3^{}(x)=
dx_1^{}+dx_2^{}+(2+x_2^{})\,dx_3^{}$
\ for all $x\in\R^3
\hfill
$
\\[2.25ex]
has the nonsingular matrix (4.2) on the set $\Xi=\{x\colon x_2^{}\ne0\}$ 
\vspace{0.35ex}
(the determinant of this matrix is ${\rm det}\,w(x)=x_2^{}\ne 0$ for all $x\in \Xi).$
\vspace{0.35ex}
By Property 4.2, restrictions of the functions 
$F_\xi^{}\colon x\to x_\xi^{}$ for all $x\in\R^3,\ \xi=1,2,3,$
are a basis of first integrals on any domain $\Omega\subset\Xi$ of the 
Pfaff system of equations (4.11).

Under the condition $x_2^{}=0$ 
the Pfaff system of equations (4.11) is 
the first-order ordinary differential equation $dx_1^{}+2\,dx_3^{}=0.$ 
\vspace{0.5ex}
This differential equation has the general integral
$
F\colon (x_1^{}\,,x_3^{})\to  x_1^{}+2x_3^{}$
for all $(x_1^{}\,,x_3^{})\in\R^2.$
\\[2.5ex]
\indent
{\bf
4.3. Existence criterion of first integral}
\\[0.75ex]
\indent
The Pfaff system of equations {\rm(Pf)} induces the linear differential forms (0.2). We add $n-m$ linear differential forms 
\\[2ex]
\mbox{}\hfill                                         
$
\displaystyle
\omega_\zeta^{}(x)=
\sum_{i=1}^n\,
w_{\zeta i}^{}(x)\,dx_i^{}$
\ for all $x\in G,
\ \ \ 
\zeta=m+1,\ldots,n,
$
\hfill(4.12)
\\[2.25ex]
with coefficients 
\vspace{0.5ex}
$w_{\zeta i}^{}\in C(G),\ 
\zeta=m+1,\ldots,n,\ 
i=1,\ldots,n,$ to the 1-forms (0.2) such that the 
set of the linear differential forms (4.12) and (0.2) 
\\[2ex]
\mbox{}\hfill                                          
$
\displaystyle
\omega_\xi^{}(x)=
\sum_{i=1}^n\,
w_{\xi i}^{}(x)\,dx_i^{}$
\ for all $x\in G,
\ \ \ 
\xi=1,\ldots,n,
$
\hfill(4.13)
\\[2.25ex]
are  not linearly bound on the domain $G.$ 
We form the square matrix of order $n$
\\[2ex]
\mbox{}\hfill                                           
$
\displaystyle
\widetilde{w}(x)=
\bigl\|w_{\xi i}^{}(x)\bigr\|$
\ for all $x\in G.
$
\hfill(4.14)
\\[2.25ex]
The matrix (4.14) is nonsingular on a domain $\Omega\subset G$ with $\mu {\rm C}_G\Omega=0.$
Then the matrix (4.14) on the domain $\Omega$ has the inverse matrix
\\[2ex]
\mbox{}\hfill                                            
$
\displaystyle
g(x)=\bigl\|g_{i\xi}^{}(x)\bigr\|$
\ for all $x\in\Omega.
$
\hfill (4.15)
\\[2.25ex]
The matrix (4.15) is a  nonsingular on the domain $\Omega$ square matrix of order $n$ and 
\\[2ex]
\mbox{}\hfill                                            
$
\displaystyle
\widetilde{w}(x) g(x)=E$
\ for all $x\in\Omega,
$
\hfill(4.16)
\\[2.25ex]
where $E$ is the identity matrix of order $n.$

We build the $n$ linear differential ope\-ra\-tors of first order
\\[2ex]
\mbox{}\hfill                                             
$
\displaystyle
{\frak G}_i^{}(x)=
\sum_{\xi=1}^n\,
g_{\xi i}^{}(x)
\partial_{x_\xi^{}}^{}$
\ for all $x\in \Omega,
\ \ \ 
i=1,\ldots,n,
$
\hfill(4.17)
\\[2ex]
which are not linearly bound on the domain $\Omega$ 
(because the matrix (4.15) is nonsingular on the domain $\Omega).$

The ope\-ra\-tors (4.17) and 1-forms (4.13) are called {\it contragredient} if 
the coordinate relations (4.16) are hold.

By the contragredient ope\-ra\-tors (4.17) and 1-forms (4.13), 
using the identity (4.16),  we get 
the total differential of any scalar function $F\in C^1(\Omega)$ have the form 
\\[1.5ex]
\mbox{}\hfill                                            
$
\displaystyle 
d F(x)=
\sum_{i=1}^n\, 
{\frak G}_i^{} F(x)\;\! 
\omega_i^{}(x)$
\ for all $x\in\Omega.
$
\hfill(4.18)
\\[2ex]
\indent
By virtue of  (4.18) and Definition 4.1, we obtain an {\sl existence criterion of a first integral 
for a Pfaff system of equations}.
\vspace{0.35ex}

{\bf Theorem 4.2.}                                        
{\it
A scalar function $F\in C^1(\Omega)$ is a first integral on a domain $\Omega\subset G$
of  the Pfaff system of equations {\rm (Pf)} with $\omega_j^{}\in C(G),\ j=1,\ldots,m,$ 
if and only if the following conditions hold
}
\\[1.25ex]
\mbox{}\hfill                                              
$
\displaystyle
{\frak G}_\zeta^{} F(x)=0$
\ for all $x\in\Omega,
\ \ \
\zeta=m+1,\ldots,n.
$
\hfill(4.19)
\\[2.25ex]
\indent
{\bf                                                   
4.4. Integral equivalence with linear homogeneous system 
of partial differential equations                                                                          
}
\\[0.75ex]
\indent
By the existence criterion of a first integral for the Pfaff system of equations (Pf) (Theorem 4.2), 
using a linear homogeneous system of partial differential equations, we can build 
an integral basis of  the Pfaff system of equations (Pf).
\vspace{0.35ex}

{\bf Theorem 4.3.}\!                                       
{\it
The scalar functions {\rm(2.2)} are a basis of first integrals on a domain  $\!G^\prime\!\subset\! G\!$
for the Pfaff system of equations {\rm (Pf)} with $\omega_j^{}\in C(G),\ j=1,\ldots,m,$ 
if and only if 
the functions {\rm(2.2)} are a basis of first integrals on the  domain $G^\prime\subset\Omega\subset G$ 
for the linear ho\-mo\-ge\-ne\-ous system of partial differential equations}
\\[2ex]
\mbox{}\hfill                                             
$
{\frak G}_\zeta^{}(x)\;\!y=0,
\quad 
\zeta=m+1,\ldots,n,
$
\hfill(4.20)
\\[2ex]
{\it
induced by the operators} (4.17).

{\sl Proof}.\! From the system of identities\! (4.19),\! we get the function 
$\!F\colon\Omega\to\R\!$ is a first in\-teg\-ral on a domain $\Omega\subset G$ of 
the linear ho\-mo\-ge\-ne\-ous system of partial differential equations (4.20) (by Definition 2.1).
Then, by the definitions of integral bases for the Pfaff system of equ\-a\-ti\-ons (Definition 4.3) and 
for the linear ho\-mo\-ge\-ne\-ous system of partial differential equ\-a\-ti\-ons (Definition 2.2), from 
Theorem 4.2, we obtain the criterion formulated in Theorem 4.3. \k
\vspace{0.5ex}

The systems (Pf) and (4.20) are called {\it contragredient}.
\vspace{0.5ex}

{\bf Definition 4.4.}
{\it
We'll say that a Pfaff system of equations and 
a linear ho\-mo\-ge\-ne\-ous system of partial differential equ\-a\-ti\-ons
are \textit{\textbf{integrally equivalent}} 
on some domain if  on this domain
each first integral of the first system is a first integral of the second system 
and on the contrary each first integral of the second system 
is a first integral of the first system.
}
\vspace{0.5ex}

By Theorem 4.3, the linear ho\-mo\-ge\-ne\-ous system of partial differential equ\-a\-ti\-ons (4.20)
is integrally equivalent on the domain $\Omega$ to the contragredient
Pfaff system of equations (Pf).

We supplement the 1-forms (0.2)  to the 1-forms (4.13) with only one condition to the  
linear differential forms (4.12): the 1-forms (4.13) is not linearly bound on the domain $\Omega.$ 
At this viewpoint, the contragredient linear ho\-mo\-ge\-ne\-ous system of partial differential equ\-a\-ti\-ons (4.20) 
to the Pfaff system of equations (Pf) is constracted ambiguously.
At the same time it does not influence (accurate within  functional expression of basis integrals)
on an integral basis of the Pfaff system of equations (Pf) and is regulated by Theorem 4.1.

In regard to the domain $\Omega,$ we have this domain is established by the domain of definition $G$ of 
the Pfaff system of equations (Pf) and 
corrected by the possibility of construction of the inverse matrix (4.15) to the matrix (4.14).
\vspace{0.5ex}

{\bf Example 4.3.}                                          
Let us consider the Pfaff system of equations
\\[2ex]
\mbox{}\hfill                                              
$
\omega_1^{}(x)=0,\qquad
\omega_2^{}(x)=0
$
\hfill(4.21)
\\[1.25ex]
with the 1-forms
\\[1.5ex]
\mbox{}\hfill
$
\omega_1^{}(x)=
dx_1^{}-dx_2^{}-
(x_1^{}x_2^{}+x_2^2-
2x_3^2-2x_3x_4)\;\!
x_2^{-1}(x_3^{}-x_4^{})^{-1}\,dx_3^{}+{}
\hfill
$
\\[2ex]
\mbox{}\hfill
$
{}+(x_1^{}x_2^{}+x_2^2-
2x_3x_4-2x_4^2)\;\!
x_2^{{}-1}(x_3^{}-x_4^{})^{-1}\,dx_4^{}$
\ for all $x\in\Xi,
\hfill
$
\\[2.5ex]
\mbox{}\hfill
$
\omega_2^{}(x)=
dx_1^{}+dx_2^{}-
(x_1^{}x_2^{}-x_2^2+
2x_3^2+2x_3x_4)\;\!
x_2^{-1}(x_3^{}-x_4^{})^{-1}\,dx_3^{}+{}
\hfill
$
\\[2ex]
\mbox{}\hfill
$
{}+(x_1^{}x_2^{}-x_2^2+
2x_3x_4+2x_4^2)\;\!x_2^{{}-1}(x_3^{}-x_4^{})^{-1}\,dx_4^{}$
\ for all $x\in\Xi,
\hfill
$
\\[2ex]
where $\Xi=\{x\colon x_2^{}\ne 0,\ x_4^{}\ne x_3\}.$
\vspace{0.35ex}

We add two 1-forms
\\[2ex]
\mbox{}\hfill
$
\omega_3^{}(x)=
x_3^{} x_2^{-1}
(x_3^{}-x_4^{})^{-1}\,dx_3^{}-
x_4^{} x_2^{-1}
(x_3^{}-x_4^{})^{-1}\,dx_4^{}$
\ for all $x\in\Xi,
\hfill
$
\\[2.75ex]
\mbox{}\hfill
$
\omega_4^{}(x)=
{}-(x_3^{}-x_4^{})^{-1}\,dx_3^{}+
(x_3^{}-x_4^{})^{-1}\,dx_4^{}$
\ for all $x\in\Xi
\hfill
$
\\[2.25ex]
to the linear differential forms $\omega_1^{}$ and $\omega_2^{}.$ 
\vspace{0.25ex}

The square matrix of fourth order (4.14) is 
\vspace{0.35ex}
generated by the coefficients of the 1-forms $\omega_i^{}\,,\ i=1,\ldots,4.$
\vspace{0.35ex}
Since the determinant ${\rm det}\, \widetilde{w}(x)=2x_2^{{}-1}(x_3^{}-x_4^{})^{-1}\ne 0$
for all $x\in\Xi,$ we see that the matrix $\widetilde{w}$ is nonsingular on the set $\Xi.$
Therefore the 1-forms $\omega_i^{}\,,\ i=1,\ldots,4,$ are not linearly bound on any domain $\Omega\subset \Xi.$ 

We introduce the linear differential operators
\\[1.75ex]
\mbox{}\hfill
$
{\frak G}_1^{}(x)=
0{,}5\, \partial_{x_1^{}}^{}-
0{,}5\, \partial_{x_2^{}}^{}$
\ for all $x\in\Xi,
\qquad
{\frak G}_2^{}(x)=
0{,}5\,\partial_{x_1^{}}^{}+\,
0{,}5\,\partial_{x_2^{}}^{}$
\ for all $x\in\Xi,
\hfill
$
\\[2ex]
\mbox{}\hfill
$
{\frak G}_3^{}(x)=
2(x_3^{}+x_4^{})\,\partial_{x_2^{}}^{}+
x_2^{}\,\partial_{x_3^{}}^{}+
x_2^{}\,\partial_{x_4^{}}^{}$
\ for all $x\in\Xi,
\hfill
$
\\[2ex]
\mbox{}\hfill
$
{\frak G}_4^{}(x)=
{}-x_1^{}\,\partial_{x_1^{}}^{}+
x_2^{}\,\partial_{x_2^{}}^{}+
x_4^{}\,\partial_{x_3^{}}^{}+
x_3^{}\,\partial_{x_4^{}}^{}$
\ for all $x\in\Xi,
\hfill
$
\\[2.25ex]
which are contragredient to the 1-forms $\omega_i^{}\,,\ i=1,\ldots,4.$

The linear ho\-mo\-ge\-ne\-ous system of partial differential equ\-a\-ti\-ons
\\[2ex]
\mbox{}\hfill                                         
$
{\frak G}_3^{}(x)\;\!y=0,
\quad
{\frak G}_4^{}(x)\;\!y=0
$
\hfill(4.22)
\\[2ex]
is contragredient to the Pfaff system of equations (4.21).

The restrictions of the functions
\\[1.75ex]
\mbox{}\hfill                                          
$
F_1^{}\colon x\to
x_1^{}(x_3^{}-x_4^{})^{{}-1}$
\ for all $x\in
\{x\colon x_4^{}\ne x_3^{}\},
\hfill
$
\\
\mbox{}\hfill(4.23)
\\[0.5ex]
\mbox{}\hfill
$
F_2^{}\colon x\to 
x_1^2\bigl(x_2^2-
(x_3^{}+x_4^{})^2\bigr)$
\ for all $x\in\R^4
\hfill
$
\\[2.25ex]
are a basis of first integrals [96, p. 200; 34; 41] on the domain $\Omega$ of system (4.22).

By Theorem 4.3, the restrictions of the functions (4.23) 
are an integral basis on the domain $\Omega$ of 
the Pfaff system of equations (4.21).
\\[2ex]
\indent
{\bf                                              
4.5. Transformation of a Pfaff system of equations by known first integrals
}
\\[1ex]
\indent
{\bf Theorem 4.4.}                                    
{\it
If the Pfaff system of equations {\rm(Pf)} has $k,\ 1\leq k\leq m,$ 
functionally independent on a domain $G^\prime\subset G$ first integrals {\rm(2.2)},
then this system on a domain $\Omega\subset G^\prime$ with $\mu{\rm C}_{G^\prime}\Omega=0$
can be reduced to the form 
\\[1ex]
\mbox{}\hfill                                          
$
dF_\xi^{}(x)=0,
\ \
\xi=1,\ldots,k,
\hfill
$
\\[-0.1ex]
\mbox{}\hfill {\rm(4.24)}
\\[-0.1ex]
\mbox{}\hfill
$
\omega_{j_\lambda^{}}^{}(x)=0,
\ \
j_\lambda^{}\in\{1,\ldots,m\},
\ \ \ 
\lambda=1,\ldots,m-k.
\hfill
$
\\[2.5ex]
by a nonsingular linear transformation of the {\rm 1}-forms {\rm(0.2)}.
}
\indent

{\sl Proof}.
If the functions (2.2) are first integrals on a domain $G^\prime\subset G$ of system (Pf),
then, by Definition 4.1, 
there exist the scalar functions $b_{\xi j}^{}\in C^1(G),\ \xi=1,\ldots,k,\ j=1,\ldots,m,$
such that the total differentials 
\\[2ex]
\mbox{}\hfill                                          
$
\displaystyle
dF_\xi^{}(x)=
\sum\limits_{j=1}^{m}\,
b_{\xi j}^{}(x)\;\! \omega_j^{}(x)$
\ for all $x\in G^\prime,
\ \ \
\xi=1,\ldots,k.
$
\hfill(4.25)
\\[2ex]
\indent
Since the first integrals (2.2) are functionally independent on a domain 
\vspace{0.5ex}
$G^\prime$ and the 1-forms (0.2) are not linearly bound on the domain $G,$ we see that 
\vspace{0.5ex}
the matrix $b(x)=\bigl\|b_{\xi j}^{}(x)\bigr\|_{\textstyle{}_{k \times m}}$ for all $x\in G^\prime$ 
(this matrix is induced by the coefficients of expansion (4.25))
has ${\rm rank}\,b(x)=k$ for all $x\in\Omega,$
where $\Omega$ is a domain such that $\Omega\subset G^\prime$ and $\mu{\rm C}_{G^\prime}\Omega=0.$
\vspace{0.25ex}

Without loss of generality it can be assumed that 
the square matrix  $\widehat{b}(x)=
\bigl\|b_{\xi j}^{}(x)\bigr\|$ for all $x \in\Omega$ of order $k$ 
(we obtain the matrix  $\widehat{b}$ from the restriction on the domain $\Omega$ of the matrix $b$ 
by deletion of the last $m-k$ columns)
is nonsingular on the domain $\Omega$
(it always can be received by renumbering the 1-forms $\omega_j^{}\,,\ j=1,\ldots,m).$
Then the system (Pf) is transformed into the system (4.24) under 
the nonsingular on the domain $\Omega$ linear transformation of the 1-forms
\\[2ex]
\mbox{}\hfill
$
\displaystyle
l_\xi^{}=
\sum\limits_{j=1}^{m}\,
b_{\xi j}^{}(x)\;\!
\omega_j^{}\,,
\ \
\xi=1,\ldots,k,
\qquad
l_\theta^{}=
\omega_\theta^{}\,,
\ \ 
\theta=k+1,\ldots,m.
\ \ \ \k
\hfill
$
\\[2ex]
\indent
The differential system (4.24) is constructed with the help of 
the not linearly bound on the domain $G$ differential forms
\\[2ex]
\mbox{}\hfill
$
\displaystyle
\widehat{\omega}_\xi^{}(x)=
\sum_{i=1}^{n}\,
\partial_{x_i^{}}^{}
F_\xi^{}(x)\,dx_i^{}$
\ for all $x\in\Omega,
\ \ \ 
\xi=1,\ldots,k,
\hfill
$
\\[2.5ex]
\mbox{}\hfill
$
\omega_{j_\lambda^{}}^{}(x)$
\ for all $x\in\Omega,\ \
j_\lambda^{}\in\{1,\ldots,m\},
\ \ \ 
\lambda=1,\ldots,m-k.
\hfill
$
\\[2.25ex]
\indent
In this connection, we have
\vspace{0.35ex}

{\bf Theorem 4.5.}                                         
{\it
The restrictions on a domain $\Omega\subset G^\prime,\ \mu\, {\rm C}_{G^{\;\!\prime}}\Omega=0,$
of the functions {\rm(2.2)} are $k,\ 1\leq k\leq m,$ functionally independent first integrals 
on the domain $\Omega$ of the Pfaff system of equations {\rm(4.24)}.
} 
\vspace{0.35ex}

By Theorems 4.4 and 4.5, we obtain an {\sl existence criterion of 
functionally independent first integrals for a Pfaff system of equations}.
\vspace{0.35ex}

{\bf Theorem 4.6.}                                         
{\it
The Pfaff system of equations {\rm(Pf)} has $k,\ 1\leq k\leq m,$
functionally in\-de\-pen\-dent first integrals on a domain $\Omega\subset G$ 
if and only if this system on the domain $\Omega$ can be reduced to the form {\rm (4.24)}
by the nonsingular linear transformation of the {\rm 1}-forms {\rm(0.2)}.
}
\\[2ex]
\indent
{\bf                                                    
4.6. Closed systems
}
\\[1ex]
\indent
{\bf Definition 4.5.}
{\it
The Pfaff system of equations {\rm(Pf)} is called \textit{\textbf{closed}} on a domain $\Omega\subset G$ if
a contragredient linear ho\-mo\-ge\-ne\-ous system of partial differential equ\-a\-ti\-ons {\rm(4.20)} is 
complete on the domain $\Omega.$
}
\vspace{0.35ex}

For example, the contragredient linear ho\-mo\-ge\-ne\-ous system of partial differential equ\-a\-ti\-ons (4.22) to 
the Pfaff system of equations (4.21) is complete on a domain $\Omega\subset\Xi.$  Indeed,
since the basis of first integrals (4.23) has the dimension $n-m=4-2=2,$ we see that 
the contragredient system (4.22) is complete (by Theorem 2.13).
Thus  the Pfaff system of equations (4.21) is closed on the domain $\Omega.$ 
\vspace{0.5ex}

{\bf Theorem 4.7.}                                          
{\it
The Pfaff system of equations {\rm(Pf)} is closed on a domain $H\subset G$ 
if and only if
this system on the domain $H$ has an integral basis of dimension $m.$
}

{\sl Proof}.
By Theorem 2.10, 
a basis of first integrals of the complete system (4.20) on a neig\-h\-bo\-ur\-ho\-od of 
any point of its normalization domain $H$
has the dimension $m.$ Taking into account Theorem 4.4, we obtain 
the closed Pfaff system of equations {\rm(Pf)} and the 
contragredient complete linear ho\-mo\-ge\-ne\-ous system of partial differential equ\-a\-ti\-ons (4.20)
have the same dimensions of integral bases. 
Therefore these dimensions are equal $m.$

Thus, by Theorem 2.13, the system  {\rm(Pf)} is closed on a domain $H\subset G$
if and only if an integral basis of this system on the domain $H$ has the dimension $m.$ \k
\vspace{0.35ex}

Theorem 4.7 is a {\sl criterion of closure for a Pfaff system of equations 
in the terms of the dimension of an integral basis.}
\vspace{0.5ex}

{\bf Example  4.4.}                                          
The Pfaff system of equations (4.7) has two equations $(m=2)$ and 
two functionally independent on the space $\R^4$ first integrals (4.8) and (4.9).

Therefore the system (4.7) is closed and the functions (4.8) and (4.9) are an integral basis
of system (4.7) on the space $\R^4.$

Likewise, by Theorem 4.7, we can prove that the Pfaff system of equations (4.11) (see Example 4.2) 
is closed on space $\R^3.$ Indeed, this system has an integral basis of the dimension three on space $\R^3.$
\vspace{0.5ex}

In Theorem 4.7 we can take $H$ as a normalisation domain (Definition 2.4) of 
the contragredient system (4.20) to the system (Pf).

From Theorem 4.6 under the condition $k=m$ and Theorem 4.7, we obtain the following 
{\sl criterion of closure for a Pfaff system of equations} [53, pp. 110 -- 111].
\vspace{0.5ex}

{\bf Theorem 4.8.}                                        
{\it
The Pfaff system of equations {\rm(Pf)} is closed on a domain $\Omega\subset G$ 
if and only if this system on the domain $\Omega$ can be reduced to the differential system 
\\[1.5ex]
\mbox{}\hfill                                             
$
dF_j^{}(x)=0,
\ \ \,
j=1,\ldots,m,
$
\hfill {\rm (4.26)}
\\[1.5ex]
by the nonsingular linear transformation of the {\rm 1}-forms {\rm(0.2)}. 
}
\vspace{0.5ex}

The scalar functions $F_j^{}\in C^1(\Omega),\ j=1,\ldots,m,$ are 
first integrals on the domain $\Omega$ both for the system (4.26) and the system (Pf).
The systems (Pf) and (4.26) are integrally equivalent on the domain $\Omega\subset G.$
\vspace{0.5ex}

{\bf Example 4.5.}                                          
The Pfaff system of equations (4.7) is transformed into the system 
\\[1.5ex]
\mbox{}\hfill                                               
$
\eta_1^{}(x)=0,
\quad
\eta_2^{}(x)=0,
$
\hfill(4.27)
\\[1ex]
where the 1-forms
\\[1.5ex]
\mbox{}\hfill
$
\eta_1^{}(x)=
4x_1^{}\,dx_1^{}+
2(x_3^{}+x_4^{})\,dx_3^{}+
2(x_3^{}+x_4^{})\,dx_4^{}$
\ for all $x\in\R^4,
\hfill
$
\\[2.5ex]
\mbox{}\hfill
$
\eta_2^{}(x)=
4x_2^{}\,dx_2^{}+
2(x_3^{}-x_4^{})\,dx_3^{}-
2(x_3^{}-x_4^{})\,dx_4^{}$
\ for all $x\in\R^4,
\hfill
$
\\[1.75ex]
under the nonsingular on the space $\R^4$ linear transformation of the 1-forms $\omega_1^{}$ and $\omega_2^{}$
\\[2ex]
\mbox{}\hfill
$
\eta_1^{}=
2\;\!\omega_1^{}+ 2(1-x_2^{})\;\!\omega_2^{}\,,
\qquad
\eta_2^{}=
2\;\!\omega_1^{}- 2(1+x_2^{})\;\!\omega_2^{}.
\hfill
$
\\[2ex]
\indent
Since 
\\[1.5ex]
\mbox{}\hfill
$
\eta_1^{}(x)=
d\bigl(2x_1^2+
(x_3^{}+x_4^{})^2\bigr)$
\ for all $x\in\R^4,
\hfill
$
\\[2.5ex]
\mbox{}\hfill
$
\eta_2^{}(x)=
d\bigl(2x_2^2+(x_3^{}-
x_4^{})^2\bigr)$
\ for all $x\in\R^4,
\hfill
$
\\[2ex]
we see that the system (4.27) can be reduced to the form
\\[2ex]
\mbox{}\hfill
$
d\bigl(2x_1^2+(x_3^{}+x_4^{})^2\bigr)=0,
\qquad
d\bigl(2x_2^2+(x_3^{}-x_4^{})^2\bigr)=0.
\hfill
$
\\[2.25ex]
\indent
By Theorem 4.8, the system (4.7) is closed and the first integrals (4.8) and (4.9) are 
an integral basis on the space $\R^4$ of system (4.7)
(using Theorem 4.7, the same result was obtained in Example 4.4).

Using the nonsingular on the space $\R^4$ linear transformation of the 1-forms $\omega_1^{}$ and $\omega_2^{}$
\\[1ex]
\mbox{}\hfill
$
\sigma_1^{}=
2\omega_1^{}+
2(1-x_2^{})\omega_2^{}\,,\qquad
\sigma_2^{}=
\omega_2^{}\,,
\hfill
$
\\[1.5ex]
we get the system (4.7) can be reduced to the system
\\[1.5ex]
\mbox{}\hfill                                               
$
\sigma_1^{}(x)=0,
\quad
\sigma_2^{}(x)=0,
$
\hfill(4.28)
\\[1.5ex]
where the 1-forms
\\[1.5ex]
\mbox{}\hfill
$
\sigma_1^{}(x)=
4x_1^{}\,dx_1^{}+
2(x_3^{}+x_4^{})\,dx_3^{}+
2(x_3^{}+x_4^{})\,dx_4^{}$
\ for all $x\in\R^4,
\hfill
$
\\[2.5ex]
\mbox{}\hfill
$
\sigma_2^{}(x)=
x_1^{}\,dx_1^{}-
x_2^{}\,dx_2^{}+
x_4^{}\,dx_3^{}+
x_3^{}\,dx_4^{}$
\ for all $x\in\R^4.
\hfill
$
\\[2ex]
\indent
Since 
\\[1ex]
\mbox{}\hfill
$
\sigma_1^{}(x)=
d\bigl(2x_1^2+
(x_3^{}+x_4^{})^2\bigr)$
\ for all $x\in\R^4,
\hfill
$
\\[2.5ex]
\mbox{}\hfill
$
2\sigma_2^{}(x)=
d\bigl(x_1^2-x_2^2+2x_3^{}x_4^{}\bigr)$
\ for all $x\in\R^4,
\hfill
$
\\[2.25ex]
we see that the system (4.28) can be reduced to the form 
\\[2ex]
\mbox{}\hfill
$
d\bigl(2x_1^2+(x_3^{}+
x_4^{})^2\bigr)=0,\qquad
d\bigl(x_1^2-x_2^2+
2x_3^{}x_4^{}\bigr)=0.
\hfill
$
\\[2.25ex]
\indent
By Theorem 4.8, the system (4.7) is closed and the first integrals (4.8) and (4.10) are 
an integral basis on the space $\R^4$ of system (4.7).
\\[2.5ex]
\indent
{\bf                                                   
4.7. Interpretation of closure in terms of differential forms}
\\[1ex]
\indent
In [58], the interpretation of complete solvability for the total differential system (TD) 
in terms of differential forms was given. We give the interpretation of closure 
for the Pfaff system of equations (Pf) in terms of differential forms.
\vspace{0.5ex}

{\bf Lemma 4.1.}                                          
{\it
Suppose the linear differential forms
$\omega_\rho^{}\in C^\infty(G),\ \rho=1,\ldots,s.$
Then the system of exterior differential identities
\\[2ex]
\mbox{}\hfill                                             
$
d\omega_\rho^{}(x)\wedge
\omega_1^{}(x)\wedge
\ldots\wedge
\omega_s^{}(x)=0$
\ for all $x \in G,
\ \ \ \rho=1,\ldots,s,
$
\hfill{\rm(4.29)}
\\[2.25ex]
is invariant under the nonsingular on the domain $G$ linear transformation of the {\rm 1}-forms
$\omega_\rho^{},\  \rho=1,\ldots,s.$
}
\vspace{0.35ex}

{\sl Proof.}
Let the 1-forms $l_\delta^{}\,,\ \delta=1,\ldots,s,$ be 
nonsingular on the domain $G$ linear com\-bi\-na\-tions of the 1-forms
$\omega_\rho^{}\,,\ \rho=1,\ldots,s,$ i.e.,
\\[2ex]
\mbox{}\hfill                                             
$
\displaystyle
l_\delta^{}(x)=
\sum\limits_{\rho=1}^{s}\,
\Psi_{\delta\rho}^{}(x)\;\!
\omega_\rho^{}(x)$
\ for all $x \in G,
\ \ \ 
\delta=1,\ldots,s,
$
\hfill{\rm(4.30)}
\\[2.25ex]
where the scalar functions 
\vspace{0.5ex}
$\Psi_{\delta\rho}^{}\in C^\infty(G),\ 
\delta=1,\ldots,s,\ \rho=1,\ldots,s,$
the square matrix $\Psi(x)=\|\Psi_{\delta\rho}^{}(x)\|$
for all $x\in G$ of order $s$ is nonsingular on the domain $G.$
Then the exterior product
\\[1ex]
\mbox{}\hfill
$
\mathop{\wedge}\limits_{\rho=1}^{s} 
l_{\rho}^{}(x)=
{\rm det}\,\Psi(x)
\Bigl(\,\mathop{\wedge}
\limits_{\rho=1}^{s} 
\omega_{\rho}^{}(x)\Bigr)$
\ for all $x\in G.
\hfill
$
\\[2.25ex]
\indent
Therefore the system of identities
\\[2ex]
\mbox{}\hfill                                              
$
d\,l_{\delta}^{}(x)\wedge
l_1^{}(x)\wedge\ldots\wedge l_s^{}(x)=0$
\ for all $x \in G,
\ \
\delta=1,\ldots,s,
$
\hfill (4.31)
\\[2.25ex]
is valid if and only if 
\\[2ex]
\mbox{}\hfill                                              
$
d\,l_{\delta}^{}(x)\wedge
\omega_1^{}(x)\wedge\ldots\wedge
\omega_s^{}(x)=0$
\ for all $x \in G,
\ \
\delta=1,\ldots,s.
$
\hfill(4.32)
\\[2.25ex]
\indent
By the representations (4.30), the exterior differentials
\\[2ex]
\mbox{}\hfill
$
\displaystyle
d\,l_\delta^{}(x)=
\sum\limits_{\rho=1}^{s}\,
\Psi_{\delta\rho}^{}(x)\,
d\;\!\omega_\rho^{}(x)+
\sum\limits_{\rho=1}^{s}\,
d\Psi_{\delta\rho}^{}(x)\wedge
\omega_\rho^{}(x)$
\ for all $x \in G,
\ \
\delta=1,\ldots,s.
\hfill
$
\\[2ex]
\indent
Thus the exterior products
\\[2ex]
\mbox{}\hfill                                              
$
\displaystyle
d\,l_{\delta}^{}(x)\wedge
\omega_1^{}(x)\wedge\ldots\wedge
\omega_s^{}(x)=
\sum\limits_{\rho=1}^{s}\,
\Psi_{\delta\rho}^{}(x)
d\,\omega_\rho^{}(x)\wedge
\omega_1^{}(x)\wedge\ldots\wedge
\omega_s^{}(x)
\hfill
$
\\[-0.25ex]
\mbox{}\hfill (4.33)
\\[-0.25ex]
\mbox{}\hfill
for all $x \in G,
\ \
\delta=1,\ldots,s.
\hfill
$
\\[2ex]
\indent
Using the identities (4.29), from the identities (4.33), we get the identities (4.32).
Therefore the identities (4.31) are consistent. This yields that the system of identities (4.29) is 
invariant un\-der the nonsingular on the domain $\!G\!$ transformation (4.30) of the \!{\rm 1}-forms
\vspace{0.75ex}
$\!\omega_\rho^{}, \rho\!=\!1,\ldots,s.\k\!\!$

{\bf Theorem 4.9.}                                         
{\it
\vspace{0.5ex}
If the Pfaff system of equations (Pf) with $\omega_j^{}\in C^\infty(G),\ j=1,\ldots,m,$ has 
$m$ functionally independent on a domain $G^\prime\subset G$ first integrals, then  
the exterior products}
\\[2ex]
\mbox{}\hfill                                              
$
d\omega_j^{}(x)\wedge
\omega_1^{}(x)\wedge\ldots\wedge
\omega_m^{}(x)=0$
\ for all $x\in G^\prime,\ \
j=1,\ldots,m.
$
\hfill(4.34)
\\[2.5ex]
\indent
{\sl Proof.}
Let the Pfaff system of equations (Pf) has
the $m$ functionally independent on a domain $G^\prime\subset G$ first integrals
\\[1ex]
\mbox{}\hfill                                              
$
F_j^{}\colon x\to F_j^{}(x)$
\ for all $x \in G^\prime, \ \ 
j=1,\ldots,m.
$
\hfill(4.35)
\\[2.25ex]
\indent
By Definition 4.1, the total differentials
\\[2ex]
\mbox{}\hfill                                              
$
\displaystyle
dF_j^{}(x)=
\sum\limits_{\zeta=1}^{m}\,
b_{j\zeta}^{}(x)\;\!
\omega_\zeta^{}(x)$
\ for all $x\in G^\prime,
\ \
j=1,\ldots,m,
$
\hfill(4.36)
\\[2.25ex]
where the scalar functions 
\vspace{0.5ex}
$b_{j\zeta}^{}\in C^\infty(G^\prime),\ 
j=1,\ldots,m,\ \zeta=1,\ldots,m,$
the square matrix 
$b(x) =\|b_{j\zeta}^{}(x)\|$ for all $x\in G^\prime$ 
\vspace{0.5ex}
of order $m$ is nonsingular on a domain $\Omega\subset G^\prime$ with $\mu {\rm C}_{G^\prime}\Omega=0$
(see the proof of Theorem 4.4).
\vspace{0.35ex}

By the Poincar\'{e} theorem [97, p. 111] (for any differential $q$\!-form $\Theta\in C^\infty(G),$ 
we have the identity $d(d\Theta(x))=0$ for all $x\in G),$
in view of the first integrals (4.35) of system (Pf) we obtain
\\[1.5ex]
\mbox{}\hfill                                              
$
d^2F_j^{}(x)=0$
\ for all $x\in G^\prime,
\ \
j=1,\ldots,m.
$
\hfill(4.37)
\\[1.5ex]
\indent
For the 1-forms
\\[1.75ex]
\mbox{}\hfill                                               
$
\displaystyle
l_j^{}(x)=
\sum_{\zeta=1}^{m}\,
b_{j\zeta}^{}(x)\;\!
\omega_\zeta^{}(x)$
\ for all $x\in G^\prime,
\ \
j=1,\ldots,m,
$
\hfill(4.38)
\\[1.75ex]
using the identities (4.36) and (4.37), we get  the exterior products
\\[2ex]
\mbox{}\hfill                                               
$
d\,l_j^{}(x)\wedge
l_j^{}(x)\wedge\ldots\wedge
l_m^{}(x)=0$
\ for all $x\in G^\prime,
\ \
j=1,\ldots,m.
$
\hfill(4.39)
\\[2.25ex]
\indent
Since the matrix $b$ is nonsingular  on a domain $\Omega\subset G^\prime,$ we see that 
from the rep\-re\-sen\-ta\-tions (4.38) it follows that 
the 1-forms $l_j^{}\,,\ j=1,\ldots,m,$ on the domain $\Omega$ are  the result of 
the nonsingular on the domain $\Omega$ linear transformation of the \!{\rm 1}-forms
$\omega_j^{}\,,\ j=1,\ldots,m.$
Therefore, by Lemma  4.1, the identities (4.39) are valid  if and only if the identities hold
\\[2ex]
\mbox{}\hfill
$
d\omega_j^{}(x)\wedge
\omega_j^{}(x)\wedge\ldots\wedge
\omega_m^{}(x)=0$
\ for all $x\in \Omega,
\ \
j=1,\ldots,m.
\hfill
$
\\[2.25ex]
\indent
Whence, using $\omega_j^{}\in C^\infty(G),\ j=1,\ldots,m,\ \mu {\rm C}_{G^\prime}\Omega=0,$ 
we get the identities (4.34). \k
\vspace{1ex}

From Theorems 4.7 and 4.9 we obtain the following 
\vspace{0.75ex}

{\bf Theorem 4.10.}                                      
{\it
\vspace{0.25ex}
If the Pfaff system of equations {\rm(Pf)} with $\omega_j^{}\in C^\infty(G),\ j=1,\ldots,m,$ 
is closed on a domain $H\subset G,$ then the exterior products}
\\[2ex]
\mbox{}\hfill
$
d\omega_j^{}(x)\wedge
\omega_1^{}(x)\wedge\ldots\wedge
\omega_m^{}(x)=0$
\ for all $x\in H,\ \
j=1,\ldots,m.
\hfill
$
\\[2.5ex]
\indent
{\bf Theorem 4.11.}                                      
{\it
If the system of the exterior identities
\\[2ex]
\mbox{}\hfill                                            
$
d\omega_j^{}(x)\wedge
\omega_1^{}(x)\wedge\ldots\wedge
\omega_m^{}(x)=0$
\ for all $x\in G,\ \
j=1,\ldots,m, 
$
\hfill{\rm(4.40)}
\\[2.25ex]
is valid, then the Pfaff system of equations {\rm(Pf)} with 
\vspace{0.35ex}
$\omega_j^{}\in C^\infty(G),\ j=1,\ldots,m,$ is closed on a domain $\Omega\subset G,$
where $\Omega$ is a domain such that $\mu\,{\rm C}_{G}\Omega=0.$
}
\vspace{0.5ex}

{\sl Proof.} 
Under the condition $n-m=0,$ we have the exterior products
\\[2ex]
\mbox{}\hfill
$
d\omega_j^{}(x)\wedge
\omega_1^{}(x)\wedge\ldots\wedge
\omega_m^{}(x)$
\ for all $x\in G,
\ \
j=1,\ldots,m, 
\hfill
$
\\[2.25ex]
are $(n+2)\!$-forms of  $n$ variables and 
under the condition $n-m=1,$ we get these exterior products
are $(n+1)\!$-forms of  $n$ variables.
In these cases the identities (4.40) are valid on the domain $G.$
Now let us prove that 
the system {\rm(Pf)} with $n-m=0$ and the system {\rm(Pf)} with $n-m=1$ are closed 
on a domain $\Omega\subset G,\ \mu\,{\rm C}_{G}\Omega=0.$
\vspace{0.25ex}

If $n-m=0,$ then the Pfaff system of equations {\rm(Pf)} on the domain $\Omega$ has 
$n$ first in\-teg\-rals (4.5) (by Pro\-per\-ty 4.1). These first integrals are an integral basis of system {\rm(Pf)} 
on the domain $\Omega$ (by Pro\-per\-ty 4.2). 
Taking into account Theorem 4.7, we obtain 
the Pfaff system of equations {\rm(Pf)} with $n-m=0$ is closed on the domain $\Omega.$

Suppose $n-m=1.$  
\vspace{0.25ex}
Then the matrix (4.2) of system {\rm(Pf)} is an $(n-1)\times n$ matrix and 
has ${\rm rank}\,w(x)=n-1$ for all $x\in \widetilde{\Omega},$
\vspace{0.25ex}
where $\widetilde{\Omega}\subset G$ is a domain such that $\mu\,{\rm C}_G\widetilde{\Omega}=0.$
Therefore the system {\rm(Pf)} on some domain $\Omega\subset\widetilde{\Omega}$ with  
\vspace{0.25ex}
$\mu\,{\rm C}_{\widetilde{\Omega}}\Omega=0$
can be solved for $m=n-1$ differentials.
For example, if the system {\rm(Pf)} is solved for $d\,x_1^{}\,,\ldots,d\,x_{n-1}^{}\,,$ 
then the system {\rm(Pf)} can be reduced to the 
system of  $n-1$ ordinary differential equations  
\\[1.75ex]
\mbox{}\hfill                                                
$
\dfrac{dx_\tau^{}}{dx_n^{}}=
P_\tau^{}(x_1^{},\ldots,x_n^{}),
\quad 
\tau=1,\ldots,n-1,
$
\hfill(4.41)
\\[2.25ex]
with the right hand sides $P_\tau^{}\in C^\infty(\Omega),\ \tau=1,\ldots,n-1.$
\vspace{0.25ex}

The system (4.41) has $n-1$ functionally independent on the domain $\Omega$ first integrals 
(by Theorem 1.3 with $m=1).$ Hence the system {\rm(Pf)} with $n-m=1$ has an integral basis 
of dimension $m=n-1$ on the domain $\Omega,$ where $\Omega$ is a domain from $G$ 
such that $\mu\,{\rm C}_G\Omega=0$ 
(because $\mu\,{\rm C}_{\widetilde{\Omega}}\Omega=0$ and 
$\mu\,{\rm C}_G\widetilde{\Omega}=0).$
By Theorem 4.7, the Pfaff system of equations {\rm(Pf)} with $n-m=1$ is closed on the domain $\Omega.$

Thus the assertion of Theorem 4.11 is valid for the Pfaff system of equations {\rm(Pf)} 
of codimension null $(n-m=0)$ and codimension one $(n-m=1).$

The proof of Theorem 4.11 for the Pfaff system of equations {\rm(Pf)} of codimension $s>1$ $(n-m=s)$
is by mathematical induction on $s.$

We assume that Theorem  4.11 is true for $s>1,\ s=n-m,$ i.e., if 
the system of exterior differential identities (4.40) with $m=n-s,\ s>1,$ is valid, 
then  the Pfaff system of equations {\rm(Pf)} of codimension $s$ is closed on a domain
$\Omega\subset G,\ \mu\,{\rm C}_G\Omega=0.$

By Theorem 4.8, we assume that if the system of exterior differential identities (4.40) with $m=n-s,\ s>0,$ is valid, 
then there exists a nonsingular on the domain $\Omega\subset G,$ $\mu\,{\rm C}_G\Omega=0,$
linear transformation of the 1-forms $\omega_j^{}\,,\ j=1,\ldots,m,$  such that 
the Pfaff system of equations {\rm(Pf)} of codimension  $s$ can be reduced 
to the differential system
\\[2ex]
\mbox{}\hfill
$
dF_j^{}(x)=0,
\ \
j=1,\ldots,m,
\ \ \ 
m=n-s,\ s>1.
\hfill
$
\\[2.25ex]
\indent
Let us consider the  Pfaff system of equations {\rm(Pf)} of codimension $s+1.$

Let the system of exterior differential identities (4.40) with 
$m=n-(s+1),\ s>1,$ be valid. 
If we fix $x_n^{}\,,$ then the identities (4.40) with $m=n-(s+1),\ s>1,$ are 
corresponding to the identities (4.40) with $m=n-s,\ s>1.$ 
Therefore, by the inductive assumption, 
there exists a linear transformation of the 1-forms $\omega_j^{}\,,\ j=1,\ldots,m,$ 
on the domain $\Omega\subset G,\ \mu\,{\rm C}_G\Omega=0,$ such that 
the Pfaff system of equations {\rm(Pf)} with $m=n-(s+1),\ s>1,$  can be reduced 
to the Pfaff system of equations
\\[2ex]
\mbox{}\hfill                                             
$
l_j^{}(x)=0,
\ \
j=1,\ldots,m,
\ \ \ 
m=n-(s+1),\ s>1,
$
\hfill(4.42)
\\[1.75ex]
with the 1-forms
\\[2ex]
\mbox{}\hfill                                              
$
l_j^{}(x)=
dF_j^{}(x)+g_j^{}(x)\,dx_n^{}$
\ for all $x\in\Omega,
\ \ \
j=1,\ldots,n-s-1,\ s>1.
$
\hfill (4.43)
\\[2.25ex]
\indent
Using the Poincar\'{e} identities (4.37) with 
\vspace{0.25ex}
$m=n-(s+1),\ s>1,$ we obtain 
the exterior differentials of the 1-forms (4.43)
\\[2ex]
\mbox{}\hfill                                                
$
d\,l_j^{}(x)=
dg_j^{}(x)\wedge dx_n^{}$
\ for all $x \in\Omega,
\ \ \ 
j=1,\ldots,n-s-1,\ s>1.
$
\hfill (4.44)
\\[2.25ex]
\indent
On the other hand, the 1-forms $l_j^{}\,,\ j=1,\ldots,n-s-1,\ s>1,$ 
\vspace{0.35ex}
are  the result of the nonsingular on the domain $\Omega\subset G,\ \mu\,{\rm C}_G\Omega=0,$ 
linear transformation of the \!{\rm 1}-forms $\omega_j^{},\ j=1,\ldots,n-s-1,\ s>1.$
Then, by Lemma 4.1, the  exterior products
\\[2ex]
\mbox{}\hfill
$
d\,l_j^{}(x)\wedge
l_1^{}(x)\wedge\ldots\wedge
l_m^{}(x)=0$
\ for all $x\in\Omega,
\ \
j=1,\ldots,n-s-1,\ s>1.
\hfill
$
\\[2.25ex]
\indent
Therefore the exterior differentials
\\[2ex]
\mbox{}\hfill                                                
$
\displaystyle
d\,l_j^{}(x)=
\sum\limits_{\zeta=1}^{m}\,
Q_{j\zeta}^{}(x)\wedge 
l_{\zeta}^{}(x)$
\ for all $x\in\Omega,
\ \ \ 
j=1,\ldots,n-s-1,\ s>1.
$
\hfill (4.45)
\\[2ex]
\indent
Combining the identities (4.44) and (4.45), we obtain the total differentials
\\[2ex]
\mbox{}\hfill
$
\displaystyle
dg_j^{}(x)=
\sum\limits_{\zeta=1}^{m}\,
h_{j\zeta}^{}(x)
\bigl(dF_\zeta^{}(x)+
g_\zeta^{}(x)\,dx_n^{}\bigr)=
\sum_{\zeta=1}^{m}\,
h_{j\zeta}^{}(x)\,
dF_\zeta^{}(x)+
h_j^{}(x)\,dx_n^{}
\hfill
$
\\[2ex]
\mbox{}\hfill
for all $x\in\Omega,
\quad
j=1,\ldots,n-s-1,
\ \,
s>1,
\hfill
$
\\[1ex]
where 
\\[1ex]
\mbox{}\hfill
$
\displaystyle
h_j^{}(x)=
\sum\limits_{\zeta=1}^{m}\,
h_{j\zeta}^{}(x) g_\zeta^{}(x)$
\ for all$ x\in\Omega,
\ \ \ j=1,\ldots,n-s-1, \ \ s>1.
\hfill
$
\\[2ex]
From these identities it follows that
\\[2ex]
\mbox{}\hfill
$
g_j^{}(x)=
\widehat{g}_j^{}
\bigl(F_1^{}(x),\ldots,
F_m^{}(x),x_n^{}\bigr)$
\ for all $x\in\Omega,
\ \ \ 
j=1,\ldots,n-s-1,\ s>1,
\hfill
$
\\[2.25ex]
where $\widehat{g}_j^{}\,,\ j=1,\ldots,n-s-1,\ s>1,$ 
\vspace{0.5ex}
are holomorphic scalar functions of $m+1$ variables $F_1^{}\,,\ldots, F_m^{}\,,d\,x_n^{}\,.$
\vspace{0.5ex}

Then the equations (4.42) with the 1-forms (4.43) are
\\[2ex]
\mbox{}\hfill
$
dF_j^{}+
\widehat{g}_j^{}
\bigl(F_1^{}\,,\ldots,
F_m^{}\,,x_n^{}\bigr)\, dx_n^{}=0,
\quad
j=1,\ldots,n-s-1,
\ \ s>1.
\hfill
$
\\[2ex]
\indent
This system is a system of  $m$ ordinary differential equations  and has 
a basis of first integrals of dimension $m,\ m=n-(s+1),\ s>1$ (by Theorem 1.3 with $m=1).$
Moreover, this basis is an integral basis on the domain $\Omega$ of system 
(Pf) with $n-m=s+1,\ s>1.$

By Theorem 4.7, the system (Pf) with $\!n-m\!=\!s+1,\, s>1,\!$ 
\vspace{0.75ex}
is closed on the domain $\Omega.\!\!$ \k

From the proof of Theorem 4.11, we get the following statements.
\vspace{0.75ex}

{\bf Corollary 4.1.}                                            
\vspace{0.25ex}
{\it 
The Pfaff system of equations {\rm(Pf)} with 
$\omega_j^{}\in C^\infty(G),\ j=1,\ldots,m,$ and $m=n-1$ 
is closed on a domain $\Omega\subset G$ with $\mu {\rm C}_G\Omega=0.$
}
\vspace{0.75ex}

{\bf Corollary  4.2.}                                            
\vspace{0.25ex}
{\it 
The Pfaff system of equations {\rm(Pf)} with 
$\omega_j^{}\in C^\infty(G),\ j=1,\ldots,m,$ and  $m=n-1$ 
can be reduced to the integrally equivalent on a domain $\Omega\subset G,\ \mu\,{\rm C}_G\Omega=0,$ 
system of  $n-1$ ordinary differential equations {\rm(4.41)} by 
the nonsingular on the domain $\Omega$ 
linear transformation of the {\rm 1}-forms {\rm (0.2)}.
}
\vspace{0.5ex}

Using Theorems 4.10 and 4.11, we obtain the Frobenius theorem [53, pp. 110 -- 112; 97, pp. 131 -- 136],
which is a {\sl criterion of closure for a Pfaff system of equations with the help of 
exterior products of differential forms.}
\vspace{0.5ex}

{\bf Theorem 4.12.}                                              
{\it
The Pfaff system of equations {\rm(Pf)} with 
\vspace{0.35ex}
$\omega_j^{}\in C^\infty(G),\ j=1,\ldots,m,$ 
is closed on a domain $\Omega\subset G,\ \mu {\rm C}_G\Omega=0,$
\vspace{0.35ex}
if and only if 
the system of exterior differential identities {\rm (4.40)} is valid.
}
\vspace{0.5ex}

The system of exterior differential identities {\rm (4.40)} is called [58, p. 302]
{\it the Frobenius conditions} 
of closure for the  Pfaff system of equations  (Pf).
\vspace{0.5ex}

{\bf Example 4.6.}                                                
The Pfaffian differential  equation
\\[1.75ex]
\mbox{}\hfill
$
yz\,dx+2xz\,dy+3xy\,dz=0
$
\hfill (4.46)
\\[2ex]
induces the vector field
$
A\colon (x,y,z)\to (yz,2xz,3xy)$
for all $(x,y,z)\in\R^3
$
with the  rotor
\\[1.5ex]
\mbox{}\hfill
$
{\rm rot}\,A\colon (x,y,z)\to (x,{}-2y,z)$
\ for all $(x,y,z)\in\R^3.
\hfill
$
\\[1.75ex]
\indent
The scalar product
$\!A(x,y,z)\,{\rm rot}\;\!A(x,y,z)\!=\!0\!$ for all $\!\!(x,y,z)\!\in\!\R^3,\!$ i.e.,
the vector field $\!A\!$ is orthogonal to the rotor of $\!\!A.\!$
This condition is equivalent to the Frobenius condition (4.40).

Therefore the Pfaffian differential  equation (4.46) is closed on space $\R^3$ 
and this equation has an integral basis of dimension one.

The first integral 
\\[0.75ex]
\mbox{}\hfill
$
F\colon (x,y,z)\to x y^2 z^3$
\ for all $(x,y,z)\in\R^3
\hfill
$
\\[2ex]
is an integral basis of the Pfaffian differential  equation (4.46).
\\[2ex]
\indent
{\bf                                                   
4.8. Nonclosed  systems}
\\[1ex]
\indent
Let us consider the Pfaff system of equations  (Pf) such that 
the  contragredient linear ho\-mo\-ge\-ne\-ous system of partial differential equ\-a\-ti\-ons {\rm(4.20)} is 
incomplete on a domain $\Omega\subset G.$
In this case, the Pfaff system of equations  (Pf) is said to be {\it nonclosed} on the domain $\Omega.$

Further, adding the equations of the forms (2.9) to the system (4.20), we get 
a cor\-res\-pon\-ding complete system to the incomplete system (4.20) and 
the defect $\delta,\ 1<\delta\leq m,$ of system (4.20).
For this complete system we obtain a normalization domain $H\subset\Omega$ (Definition 2.4).
Then, by Theorem 4.3, we have 
\vspace{0.5ex}

{\bf Theorem 4.13.}                                      
{\it
The nonclosed on a domain $\Omega\subset G$
Pfaff system of equations  {\rm(Pf)} has 
a basis of first integrals of dimension $m-\delta$ on 
a normalization domain $H\subset\Omega$ of the 
contragredient linear ho\-mo\-ge\-ne\-ous system of partial differential equ\-a\-ti\-ons {\rm(4.20)},
where $\delta$ is the defect of system} (4.20).
\vspace{0.5ex}

Let us remember that the complete system (4.20) has the defect $\delta=0.$ 
Then, using  Theorems  4.7 and 4.13, 
we obtain the generalizing statement about a basis of first integrals for 
the closed or nonclosed Pfaff system of equations  {\rm(Pf)}.
\vspace{0.5ex}

{\bf Theorem 4.14.}                                       
{\it
The Pfaff system of equations  {\rm(Pf)} on a normalization domain $H\subset G$ of the 
contragredient linear ho\-mo\-ge\-ne\-ous system of partial differential equ\-a\-ti\-ons {\rm(4.20)} has 
a basis of first integrals of dimension $m-\delta,$  
\vspace{0.75ex}
where $\delta$ is the defect of system} (4.20), $0\leq\delta\leq m.$

{\bf Example 4.7.}
Consider  the Pfaff system of equations  
\\[2ex]
\mbox{}\hfill                                             
$
\omega_1^{}(x)=0, 
\quad
\omega_2^{}(x)=0
$
\hfill(4.47)
\\[1.5ex]
with the 1-forms
\\[1.5ex]
\mbox{}\hfill
$
\omega_1^{}(x)=
dx_1^{}+dx_2^{}+ dx_3^{}+dx_4^{}$
\ for all $x\in\R^4,
\hfill
$
\\[2.5ex]
\mbox{}\hfill
$
\omega_2^{}(x)=
dx_1^{}+2dx_2^{}+
x_4^{}dx_3^{}+dx_4^{}$
\ for all $x\in\R^4.
\hfill
$
\\[2.25ex]
\indent
We add two 1-forms 
\\[2ex]
\mbox{}\hfill
$
\omega_3^{}(x)= dx_3^{}$
\ for all $x\in\R^4,
\qquad
\omega_4^{}(x)=dx_4^{}$
\ for all $x\in\R^4.
\hfill
$
\\[2.25ex]
to the linear differential forms $\omega_1^{}$ and $\omega_2^{}.$

The linear differential forms $\omega_i^{},\, i=1,\ldots,4,$ are not linearly bound on the space $\!\R^4.\!\!$

Using the not linearly bound on the space $\R^4$ contragredient linear differential operators
\\[1.25ex]
\mbox{}\hfill
$
{\frak G}_1^{}(x)=
2\partial_{x_1^{}}^{}-\partial_{x_2^{}}^{}$
\ for all $x\in\R^4,
\qquad
{\frak G}_2^{}(x)=
{}-\partial_{x_1^{}}^{}+\partial_{x_2^{}}^{}$
\ for all $x\in\R^4,
\hfill
$
\\[2ex]
\mbox{}\hfill
$
{\frak G}_3^{}(x)=
(x_4^{}-2)\partial_{x_1^{}}^{}+
(1-x_4^{})\partial_{x_2^{}}^{}+\partial_{x_3^{}}^{}$
\ for all $x\in\R^4,
\hfill
$
\\[2ex]
\mbox{}\hfill
$
{\frak G}_4^{}(x)=
{}-\partial_{x_1^{}}^{}+\partial_{x_4^{}}^{}$
\ for all $x\in\R^4
\hfill
$
\\[2.25ex]
to the 1-forms $\omega_i^{}\,,\ i=1,\ldots,4,$ we obtain 
the contragredient linear ho\-mo\-ge\-ne\-ous system of partial differential equ\-a\-ti\-ons
\\[1.5ex]
\mbox{}\hfill                                           
$
{\frak G}_3^{}(x)\;\!y=0,
\quad
{\frak G}_4^{}(x)\;\!y=0
$
\hfill(4.48)
\\[2ex]
to the Pfaff system of equations (4.47).

Since the Poisson bracket
\\[2ex]
\mbox{}\hfill
$
{\frak G}_{34}^{}(x)=
[{\frak G}_3^{}(x),
{\frak G}_4^{}(x)]=
{}-\partial_{x_1^{}}^{}+
\partial_{x_2^{}}^{}=
{\frak G}_2^{}(x)$
\ for all $x\in \R^4,
\hfill
$
\\[2.25ex]
we see that the system (4.48) is incomplete. 
Therefore the Pfaff system of equations (4.47) is nonclosed.

The system (4.48) with the help of the operator ${\frak G}_{34}^{}$ can be reduced 
to the complete system
\\[1ex]
\mbox{}\hfill
$
{\frak G}_3^{}(x)\;\!y=0,
\qquad
{\frak G}_4^{}(x)\;\!y=0,
\qquad
{\frak G}_{34}^{}(x)\;\!y=0,
\hfill
$
\\[2ex]
\indent
The first integral
\\[1.5ex]
\mbox{}\hfill
$
F\colon x\to 
x_1^{}+x_2^{}+x_3^{}+x_4^{}$
\ for all $x\in \R^4
\hfill
$
\\[2ex]
is an integral basis on the space $\R^4$ of this complete system.

This function is a basis of first integrals on the space $\R^4$ of the  nonclosed 
Pfaff system of equations  (4.47).
\\[2.75ex]
\indent
{\bf                                                   
4.9. Integral equivalence with  total differential system    
}
\\[1ex]
\indent
One more approach for building of an integral basis of a Pfaff system of equations is based 
on a reducing this system to integrally equivalent  total differential system. 
\vspace{0.75ex}

{\bf Definition 4.6.}
{\it
We'll say that a Pfaff system of equations and a system of total differential equations   
are \textit{\textbf{integrally equivalent}} 
on some domain if on this domain
each first integral of the first system is a first integral of the second system 
and on the contrary each first integral of the second system 
is a first integral of the first system.
}
\vspace{0.75ex}

The linear differential forms (0.2) are not linearly bound on the domain $G.$ 
Therefore the $m\times n$ matrix (4.2) has 
\vspace{0.5ex}
${\rm rank}\,w(x)=m$ for all $x\in\Omega,$ where a domain $\Omega\subset G$
and $\mu{\rm C}_G\Omega=0.$
Then, the square matrix 
\vspace{0.5ex}
$\widehat{w}(x)=\bigl\|w_{ji}^{}(x)\bigr\|$ for all $x\in G$ of order $m$ 
is nonsingular on the domain $\Omega.$ 
Thus the Pfaff system of equations (Pf) can be reduced to the 
system of total differential equations   
\\[2ex]
\mbox{}\hfill                                                
$
\displaystyle
dx_j^{}=
\sum\limits_{\nu=m+1}^{n}
\widehat{a}_{j\nu}^{}(x)\,dx_\nu^{}\,,
\quad 
j=1,\ldots,m.
$
\hfill(4.49)
\\[2.25ex]
\indent
The Pfaff system of equations (Pf) and the system of total differential equations  (4.49)
are integrally equivalent on some domain $G^\prime\subset\Omega\subset G,$ i.e., 
\vspace{0.75ex}
we have the following assertions.

{\bf Theorem 4.15.}                                          
{\it
\vspace{0.25ex}
A scalar function $F\in C^1(G^\prime)$ is a first integral on a domain $G^\prime\subset G$ 
of the Pfaff system of equations {\rm(Pf)} with $\omega_j\in C(G),\  j=1,\ldots,m,$
\vspace{0.25ex}
if and only if  this function is a first integral on the domain $G^\prime$ 
of  the system of total differential equations} (4.49).
\vspace{0.75ex}

{\bf Theorem  4.16.}                                           
{\it
\vspace{0.25ex}
The scalar functions {\rm(2.2)} are a basis of first integrals on a domain $G^\prime\subset G$ 
of the Pfaff system of equations {\rm(Pf)} with $\omega_j\in C(G),\  j=1,\ldots,m,$
\vspace{0.25ex}
if and only if  these functions are a basis of first integrals on the domain $G^\prime$ of  
the system of total differential equations} (4.49).
\vspace{0.75ex}

Using Theorems  3.6, 4.7, and 4.16, we can prove the following
\vspace{0.75ex}

{\bf Theorem 4.17.}                                           
{\it
The Pfaff system of equations {\rm(Pf)} is closed on a domain $G^\prime\subset G$
if and only if  
the total differential system {\rm (4.49)} on the domain $G^\prime$ is 
completely solvable.
}
\vspace{1.25ex}

{\bf Example 4.8.}                                              
The Pfaff system of equations 
\\[2ex]
\mbox{}\hfill                                                  
$
2x_1^{}(1+x_2^{})\,dx_1^{}+
6x_2^{}\,dx_2^{}+3x_3^{}(2+x_2^{})\,dx_3^{}+
3x_4^{}(2+x_4^{})\,dx_4^{}=0,
\hfill
$
\\[-0.15ex]
\mbox{}\hfill(4.50)
\\[-0.15ex]
\mbox{}\hfill
$
4x_1^{}(1+x_1^{})\,dx_1^{}-
6x_2^{}\,dx_2^{}+3x_3^{}(1+2x_1^{})\,dx_3^{}+
3x_4^{}(1+2x_1^{})\,dx_4^{}=0
\hfill
$
\\[2.25ex]
can be reduced to the  system of total differential equations 
\\[2ex]
\mbox{}\hfill                                                
$
dx_1^{}=
{}-\dfrac{3}{2}\,
x_3^{} x_1^{-1}\,dx_3^{}-
\dfrac{3}{2}\,
x_4^{} x_1^{-1}\,dx_4^{}\,,
\hfill
$
\\[0.25ex]
\mbox{}\hfill(4.51)
\\[0.25ex]
\mbox{}\hfill
$
dx_2^{}=
{}-\dfrac{1}{2}\,
x_3^{} x_2^{-1}\,dx_3^{}-
\dfrac{1}{2}\,
x_4^{} x_2^{-1}\,dx_4^{}\,,
\hfill
$
\\[2.25ex]
which is defined on the set $\Xi=\{x\colon x_1^{}\ne 0,\ x_2^{}\ne 0, \ 3+2x_1^{}+x_2^{}\ne0\}.$
\vspace{0.5ex}

Since the Poisson bracket
\\[2ex]
\mbox{}\hfill
$
\Bigl[
\partial_{x_1^{}}^{}-
\dfrac{3}{2}\,x_3^{} x_1^{-1}\,
\partial_{x_1^{}}^{}-
\dfrac{1}{2}\,x_3^{} x_2^{-1}\,
\partial_{x_4^{}}^{}\,,\
\partial_{x_2^{}}^{}-
\dfrac{3}{2}\,x_4^{} x_1^{-1}\,
\partial_{x_3^{}}^{}-
\dfrac{1}{2}\,x_4^{} x_2^{-1}\,
\partial_{x_4^{}}^{}\Bigr]=
{\frak O}
\hfill
$
\\[2.75ex]
\mbox{}\hfill
for all $x\in\widetilde{\Xi},
\quad 
\widetilde{\Xi}=
\{x\colon x_1^{}\ne 0,\ x_2^{}\ne 0\},
\hfill
$
\\[2.25ex]
we see that the system (4.51) is completely solvable on any domain 
$\widetilde{G}\subset\widetilde{\Xi}.$
\vspace{0.25ex}

By Theorem 4.17, the Pfaff system of equations (4.50) is closed on any domain $\widetilde{G}\subset\Xi.$

\newpage

The functionally independent first integrals 
\\[2ex]
\mbox{}\hfill
$
F_1^{}\colon x\to 
2x_1^2+3x_3^2+3x_4^2$
\ for all $x\in\R^4,
\hfill
$
\\[2.5ex]
\mbox{}\hfill
$
F_2^{}\colon x\to 
2x_2^2+x_3^2+x_4^2$
\ for all $x\in\R^4
\hfill
$
\\[2.25ex]
are a basis of first integrals on any domain $\widetilde{G}\subset\widetilde{\Xi}$ of 
the system (4.51), and a basis of first integrals on any domain $G^{\,\prime}\subset \Xi$ of 
the system (4.50).

}
\end{document}